\documentclass[12pt,titlepage]{amsart}
\usepackage{latexsym,amsfonts,amssymb,amsmath,amsthm} 
\DeclareFontFamily{U}{wncy}{}
\DeclareFontShape{U}{wncy}{m}{n}{<->wncyr10}{}
\DeclareSymbolFont{mcy}{U}{wncy}{m}{n}
\DeclareMathSymbol{\Sh}{\mathord}{mcy}{"58}

\newcommand{\mz}{\ensuremath{\mathbb Z}}
\newcommand{\mr}{\ensuremath{\mathbb R}}

\newcommand{\mq}{\ensuremath{\mathbb Q}}

\newcommand{\mf}{\ensuremath{\mathbb F}}

\newcommand{\mymod}{\ensuremath{\negthickspace \negmedspace \pmod}}
\newcommand{\shortmod}{\ensuremath{\negthickspace \negthickspace \negthickspace \pmod}}

\newcommand{\onehalf}{\ensuremath{ \frac{1}{2}}}
\newcommand{\thalf}{\textstyle \frac{1}{2}}
\newcommand{\notdiv}{\ensuremath{\not \; |}}
\newcommand{\notdivtext}{\ensuremath{\negthickspace \not | \, }}

\theoremstyle{plain}		
	\newtheorem{mytheo}{Theorem}[section]
	
	\newtheorem{myprop}[mytheo]{Proposition}
	\newtheorem{mycoro}[mytheo]{Corollary}
     \newtheorem{mylemma}[mytheo]{Lemma}
	
\theoremstyle{definition}
	\newtheorem{myconj}[mytheo]{Conjecture}

\theoremstyle{remark}		
	
	\newtheorem*{myrema}{Remark}

	\newtheorem*{myclaim}{Claim}

\begin{document}

\title{Low-lying Zeros of Families of Elliptic Curves}
\author{Matthew P. Young}
\address{American Institute of Mathematics, 360 Portage Ave.,
Palo Alto, CA 94306-2244}
\email{myoung@aimath.org}

\begin{abstract}
There is a growing body of evidence giving strong evidence that zeros of families of L-functions follow distribution laws of eigenvalues of random matrices.  This philosophy is known as the random matrix model or the Katz-Sarnak philosophy.  The random matrix model makes predictions for the average distribution of zeros near the central point for families of L-functions.  We study the low-lying zeros for families of elliptic curve L-functions.  For these L-functions there is special arithmetic interest in any zeros at the central point (by the conjecture of Birch and Swinnerton-Dyer and the impressive partial results towards resolving the conjecture).

We calculate the density of the low-lying zeros for various families of elliptic curves.  Our main foci are the family of all elliptic curves and a large family with positive rank.  A main challenge has been to obtain results with test functions that are concentrated close to the origin since the central point is a location of great interest.  An application is an improvement on the upper bound of the average rank of the family of all elliptic curves (conditional on the Generalized Riemann hypothesis (GRH)).  The upper bound obtained is less than $2$, which shows that a positive proportion of curves in the family have algebraic rank equal to analytic rank and finite Tate-Shafarevich group.  We show that there is an extra contribution to the density of the low-lying zeros from the family with positive rank (presumably from the ``extra" zero at the central point).
\end{abstract}
\maketitle

\section{Introduction}
The random matrix model predicts that many statistics associated to zeros of a family of L-functions can be modeled (or predicted) by the distribution of eigenvalues of large random matrices in one of the classical linear groups.  If the statistics of a family of L-functions are modeled by the eigenvalues of the group $G$ then we say that $G$ is the symmetry group (or symmetry type) associated to the family.

The statistic of interest to us in this work is the density of zeros near the central point (also known as the 1-level density).  The random matrix model predicts that the distribution of these zeros should be modeled by the eigenenvalues nearest $1$ for one of the symmetry types $G$.  All of the different groups $G$ have distinct behavior in this regard.  Therefore, computing the 1-level density gives a theoretical way to predict the symmetry type of a family.  

It is standard to assume the Generalized Riemann Hypothesis (GRH) to study the 1-level density and we do so throughout this work.  In truth, the GRH is necessary for only a handful of our results, but it simplifies arguments in some non-essential places so we use it freely even when it could be removed with extra work.

It is especially interesting to investigate the 1-level density for families of L-functions attached to elliptic curves over the rationals since zeros at the central point have important arithmetic information (by the conjecture of Birch and Swinnerton-Dyer).  These investigations have been the main focus of this work.

\subsection{Acknowledgement}
This work constitutes a large portion of my PhD thesis.  I thank my advisor, Henryk Iwaniec, for suggesting this problem and for his support and encouragement while doing this work.

\section{Preliminaries}
\label{section:preliminaries}
We begin by collecting some facts and setting the notation we will use.  We consider an elliptic curve $E/\mq$ given in general Weierstrass form
\begin{equation}
\label{eq:ellipticcurve}
E: y^2 + a_1 xy + a_3 y = x^3 + a_2x^2 + a_4 x + a_6,
\end{equation}
where each $a_i \in \mz$.  Under a change of variables $E$ can be brought into the simpler form
\begin{equation}
\label{eq:ellipticcurveb}
y^2 = x^3 + ax + b,
\end{equation}
where $a$ and $b$ are integers.  The canonical change of variables (cf. \cite{Silverman}, 46-48) uses the parameters $b_2, b_4, b_6, c_4,$ and $c_6$, where
\begin{eqnarray*}
b_2 & = & a_1^2 + 4a_2, \\
b_4 & = & 2a_4 + a_1 a_3, \\
b_6 & = & a_3^2 + 4a_6,
\end{eqnarray*}
and
\begin{eqnarray*}
c_4 & = & b_2^2 - 24 b_4, \\
c_6 & = &  -b_2^3 + 36b_2 b_4 - 216 b_6. 
\end{eqnarray*}
The curve (\ref{eq:ellipticcurve}) is then equivalent to
\begin{equation*}
y^2 = x^3 - 27 c_4 x - 54c_6.
\end{equation*}

When given by the form (\ref{eq:ellipticcurveb}), $E$ has discriminant
\begin{equation*}
\Delta = -16(4a^3 + 27b^2),
\end{equation*}
which is necessarily non-zero for the curve $E$ to be elliptic.  

The Weierstrass equation for the elliptic curve (\ref{eq:ellipticcurve}) is not unique.  Any two Weierstrass equations for the same curve are related by the admissible change of variables 
\begin{equation}
\label{eq:admissiblecov}
\begin{array}{lcl} x & = & u^2 x' + r \\
y & = & u^3 y' + su^2 x' + t,
\end{array}
\end{equation}
where $u$, $r$, $s$, and $t$ are integers and $u \neq 0$.  Under this change of variables the discriminant transforms by $u^{12} \Delta' = \Delta$.  Likewise, $u^4 c_4' = c_4$ and $u^6 c_6' = c_6$.  

A standard technique in studying $E$ is to reduce the equation (\ref{eq:ellipticcurve}) $\mymod{p}$ for every prime $p$.  The equation (\ref{eq:ellipticcurve}) is {\em minimal} for the prime $p$ if the power of $p$ dividing $\Delta$ cannot be decreased by an admissible change of variables.  The equation (\ref{eq:ellipticcurve}) is a {\em global minimal Weierstrass equation} if it is minimal for all primes simultaneously.  For any Weierstrass equation there is an admissible change of variables placing it in global minimal Weierstrass form (cf. \cite{Silverman}, Corollary 8.3).  We record here that if the admissible change of variables (\ref{eq:admissiblecov})  places (\ref{eq:ellipticcurve}) in global minimal Weierstrass form then the only primes $p$ dividing $u$ are those for which (\ref{eq:ellipticcurve}) is not minimal.  We remark that for $p >3$ if the equation (\ref{eq:ellipticcurveb}) is not minimal at $p$ then $p^4 | a$ and $p^6 | b$ (and therefore $p^{12} | \Delta$).  

Suppose $E$ is given by a global minimal Weierstrass equation (\ref{eq:ellipticcurve}).  The conductor $N$ of $E$ is then defined by
\begin{equation*}
N = \prod_{p | \Delta} p^{f_p}
\end{equation*}
where for $p > 3$ 
\begin{equation*}
f_p = \begin{cases}
1 & \text{if $p  \notdivtext c_4$, i.e. $E$ has multiplicative reduction at $p$},
\\
		2 & \text{if $p|c_4$, i.e. $E$ has additive reduction at $p$}.
	\end{cases} 
\end{equation*}
When $p = 2$ or 3 the definition of $f_p$ is more complicated (cf. \cite{SilvermanAEC2}, IV, \S 10), but it will usually be enough for our purposes that $N | \Delta$ and $f_p \leq 8$ for all primes $p$.  If necessary, the conductor can be computed using Ogg's formula and Tate's algorithm.

Continuing to assume (\ref{eq:ellipticcurve}) is a global minimal Weierstrass equation for $E$, the L-function of $E$ is defined by
\begin{equation}
\label{eq:defoflfunction}
L(s, E) = \prod_{p \negmedspace \, \notdiv  \Delta} \left({1 - a_p p^{-s}} + p^{1 - 2s} \right)^{-1} \prod_{p | \Delta} \left({1 - a_p p^{-s}} \right)^{-1},
\end{equation}
where
\begin{equation}
\label{eq:ap}
a_p = p - \#E(\mf_p), 
\end{equation}
and $\#E(\mf_p)$ is the number of affine points on $E$, when reduced $\mymod{p}$.  The central point is at $s=1$.  Since (\ref{eq:ellipticcurve}) is minimal the change of variables taking (\ref{eq:ellipticcurve}) to (\ref{eq:ellipticcurveb}) does not alter $a_p$ for $p > 3$.  Thus, for any $p > 3$, $a_p$ is given by
\begin{equation*}
a_p = - \sum_{x \shortmod{p}} \left( \frac{x^3 + ax + b}{p} \right).
\end{equation*}
We remark that for primes $p > 3$ dividing the conductor we have $a_p = \pm 1$ if $E$ has multiplicative reduction at $p$ and $a_p = 0$ if $E$ has additive reduction at $p$.
The infinite product  (\ref{eq:defoflfunction}) converges absolutely and uniformly for $\text{Re} \: s > \frac{3}{2}$, by Hasse's estimate $a_p < 2\sqrt{p}$.  According to the Shimura-Taniyama conjecture (proved by Wiles {\em et al} \cite{Wiles}, \cite{TaylorWiles}, \cite{BCDF}) there exists a weight two primitive cusp form $f(z)$ on $\Gamma_0(N)$ such that  $L(s,E)=L(s,f)$.  Further, $L(s, E)$ has analytic continuation to the complex plane and satisfies the functional equation
\begin{equation*}
\Lambda(s, E) := \left( \frac{\sqrt{N}}{2\pi} \right)^s \Gamma(s) L(s, E) = w \Lambda(2-s, E),
\end{equation*}
where $w = \pm 1$ is the root number of $E$.  

Throughout this work we will be assuming the Generalized Riemann Hypothesis holds, namely that all the nontrivial zeros of an arithmetic L-function lie on its line of symmetry.

Using the notation of Iwaniec-Luo-Sarnak \cite{ILS} we define for an L-function $L(s, E)$ the quantity
\begin{equation*}
D(E; \phi) = \sum_{\gamma_E} \phi\left( \frac{\gamma_E}{2\pi} \log X\right)
\end{equation*}
where $\phi$ is an even Schwartz class test function whose Fourier transform\footnote{Throughout we denote $\widehat{f}(y) = \int_{-\infty}^{\infty} f(x)e(-xy) dx$, \; $e(x) = \exp(2\pi i x)$} $\widehat{\phi}$ has compact support, $\gamma_E$ runs through the imaginary parts of the nontrivial zeros $\rho_E = 1 + i \gamma_E$ of $L(s, E)$ (counted with multiplicity), and $X$ is a parameter at our disposal (generally of size $N_E$,  the conductor of $E$; allowing $X$ to be only approximately $N_E$ gives us more freedom in averaging over a family).
The scaling factor $(2\pi)^{-1} \log X$ is inserted to normalize the number of zeros counted by the test function $\phi$, so that $D(E; \phi)$ should be thought of as representing the density of zeros of $L(s, E)$ near the central point $s=1$.  We will be interested in averaging $D(E; \phi)$ over certain families of automorphic forms arising from elliptic curves.  Each family we study will be of the form
\begin{equation*}
\mathcal{F} = \{ E_{d} \}
\end{equation*}
where $d$ ranges over a set $\mathcal{A}$, which is a subset of $\mz \text{ or } \mz \times \mz$, such that each $d \in \mathcal{A}$ naturally defines an elliptic curve $E_D$ over $\mq$.  The curve $E_d$ will be defined by a Weierstrass equation whose coefficients are polynomials in $d$.  Our main family will be parameterized by $d = (a, b) \in \mz^2$ where $E_d: y^2 = x^3 + ax +b$.  It may happen that $L(s, E_d) = L(s, E_{d'})$ for $d \neq d'$ trivially because a global minimal Weierstrass equation for $E_d$ equals that of $E_d'$ or more subtly because $E_d$ is isogenous to $E_{d'}$.  We take such forms with multiplicity.  We expect it should make no statistical difference whether one takes such forms $f$ with multiplicity or not.  As a general rule, we can easily make restrictions on $d$ that force $E_dw(E)$ to be minimal, but this does not change any statistics (at least in the main term) (see Section \ref{section:minimality} for a more thorough discussion).  We will often suppress the dependence of $\mathcal{F}$ on $\mathcal{A}$.

We study the weighted average
\begin{equation*}
\mathcal{D}(\mathcal{F}; \phi, w) = \sum_{E \in \mathcal{F}} D(E; \phi) w(E),
\end{equation*}
where $w(E_d) := w(d)$ is a smooth, compactly supported function (a cutoff function) on $\mr$ or $\mr^2$, whichever is appropriate.
To avoid trivialities we assume $w$ does not have total mass zero (i.e. $\widehat{w}(0) \neq 0$), and in particular that $w$ is not identically zero.
We measure the weighted sum against the total weight
\begin{equation*}
W(\mathcal{F}) = \sum_{E \in \mathcal{F}} w(E).
\end{equation*} 
Since $D(E; \phi)$ depends on $X$ we scale our cutoff function $w$ by $X$ also, in which case we use the notation $w_X$ to denote the scaling of $w$ by $X$ and $W_X$ to be the total weight scaled by $X$.  The precise scaling depends on the family and is used to pick out curves with conductors $N$ such that $\log{N}$ is asymptotically $\log{X}$ on average.  Often we simply take curves with discriminant $|\Delta| \asymp X$.  Usually in random matrix theory one takes a family with conductors restricted by $c(f) \leq X$ and let $X$ tend to infinity; for families of elliptic curves this is not practical since it is necessary to have a concretely-given set over which to average.

Katz and Sarnak predict that for a natural family $\mathcal{F}$ the average density should satisfy
\begin{equation*}
\lim_{X \rightarrow \infty} \frac{\mathcal{D}(\mathcal{F}; \phi, w_X)}{W_X(\mathcal{F})} = \int_{-\infty}^{\infty} \phi(t) \mathcal{W}(G)(t) dt,
\end{equation*}
where $\mathcal{W}(G)$ is the 1-level scaling density of eigenvalues near 1 for a symmetry group $G$ ($\mathcal{W}$ will in general be a distribution).  Such a result is called the {\it density theorem} for the family $\mathcal{F}$.  We have
\begin{equation*}
\mathcal{W}(G)(t) = 
\begin{cases} 1 & \text{if $G = U$} \\ 
1 - \frac{\sin{2 \pi t}}{2 \pi t} & \text{if $G = Sp$} \\
1 + \onehalf \delta_0(t) & \text{if $G = O$} \\
1 + \frac{\sin{2 \pi t}}{2 \pi t} & \text{if $G = SO(\text{even})$} \\
1 + \delta_0(t) - \frac{\sin{2 \pi t}}{2 \pi t} & \text{if $G = SO(\text{odd})$},
\end{cases}
\end{equation*}
where $\delta_0$ is the Dirac distribution \cite{KSM}.
In practice it is convenient to have the Fourier transforms of these distributions, which we record here.
\begin{equation}
\label{eq:symmetrytypes}
\widehat{\mathcal{W}}(G)(t) = 
\begin{cases} \delta_0(t) & \text{if $G = U$} \\ 
\delta_0(t) - \onehalf \eta(t) & \text{if $G = Sp$} \\
\onehalf +  \delta_0(t) & \text{if $G = O$} \\
\delta_0(t) + \onehalf \eta(t) & \text{if $G = SO(\text{even})$} \\
1 + \delta_0(t) - \onehalf \eta(t) & \text{if $G = SO(\text{odd})$}, 
\end{cases}
\end{equation}
where
\begin{equation*}
\eta(t) =
\begin{cases}
1 & \text{if \; $|t| < 1$}, \\
\onehalf & \text{if \; $|t| = 1$}, \\
0 & \text{if \; $|y| > 1$}.
\end{cases}
\end{equation*}
An important feature is that the Fourier transforms of $\widehat{\mathcal{W}}(O)(t)$, $\widehat{\mathcal{W}}(SO(\text{even}))(t)$, and $\widehat{\mathcal{W}}(SO(\text{odd}))(t)$ all agree for $|t| < 1$ but are distinguishable for $|t| > 1$.  Therefore, by the Plancherel Theorem, to distinguish the densities of these three symmetry types one needs to apply test functions $\phi$ whose Fourier transforms are supported outside $[-1, 1]$.

Iwaniec, Luo, and Sarnak prove that the family $H_2^\star(N)$ of primitive cusp forms of weight 2 and level $N$ ($N$ square-free) has symmetry type $O$ for test functions $\phi$ restricted by supp$\; \widehat{\phi} \subset (-2, 2)$ (see \cite{ILS}).  They also prove the forms with root number $+1$ have symmetry type $SO(\text{even})$ and the forms with root number $-1$ have symmetry type $SO(\text{odd})$ in the same range.  This is relevant for our families because the family of all elliptic curves forms a subfamily of weight two primitive cusp forms.  In particular, we do not expect to detect statistics of the root number of the family without obtaining support past $(-1, 1)$.  This potential change in behavior at $1$ is not surprising in light of the ``approximate" functional equation, which states that
\begin{equation*}
L(1, E) = \sum_{n} a_n n^{-1/2} g\left(\frac{2\pi n}{U}\right) + w_E \sum_{n} a_n n^{-1/2} g \left(\frac{2\pi n}{V}\right) 
\end{equation*}
where $g$ is a test function of a certain kind and $UV = N$.  If we take $U > N^{1 + \varepsilon}$ (i.e. sum the Fourier coefficients of length greater than $N$) then the first sum implicitly captures the root number (since the second sum is then small) whereas if $U \leq N^{1 - \varepsilon}$ then the root number explicitly occurs in the second sum.  After developing the explicit formula we shall see the analogy to this dramatic shift in behavior with a similar sum (except over primes); going past support $(-1, 1)$ is similar to taking $U$ larger than the conductor.

It is more difficult to gain large support for families of elliptic curves than for all cusp forms of weight two because the forms coming from elliptic curves compose a small subfamily of the weight two cusp forms.  Very loosely speaking, there are probably around $X^{5/6}$ elliptic curves with conductors $N \asymp X$, whereas there are about $X^2$ weight two cusp forms of levels $N \asymp X$.  

It is an open (and very interesting) question to estimate the number of elliptic curves that have conductor $N \leq X$.  Our figure of $X^{5/6}$ arises by simply counting the number of positive integers $a$ and $b$ such that $|\Delta| = 16(4a^3 + 27b^2) \asymp X$.  Fouvry, Nair, and Tenenbaum \cite{FNT} have shown that the number of non-isogenous semi-stable elliptic curves with conductor $N \leq X$ is $\gg X^{5/6}$.  In the other direction, Duke and Kowalski \cite{DK} (building on work of Brumer and Silverman \cite{BS}) have shown that the number of elliptic curves with conductor $N \leq X$ is $\ll X^{1 + \varepsilon}$ for any $\varepsilon > 0$.  Ny improvement in the exponent in this upper bound would be very interesting because it would show that for almost all levels $N$ there is no elliptic curves with conductor $N$.  Note also that our lack of knowledge in this regard illustrates why we cannot average $D(E;\phi)$ over all elliptic curves with $N_E \leq X$.

\section{Summary of Results}
\subsection{Main results}
Our main results are given in this section.
\begin{mytheo} 
\label{thm:mainresult}
Let $\mathcal{F}$ be the family of elliptic curves given by the Weierstrass equations $E_{a,b}: y^2 = x^3 + ax + b$ with $a$ and $b$ positive integers.  Let $w \in C_0^\infty (\mr^{+} \times \mr^{+})$\footnote{For us $\mr^+ = (0, \infty)$} and set $w_X (E_{a, b}) = w \left( \frac{a}{A}, \frac{b}{B}  \right)$, where $A = X^{1/3}$, $B = X^{1/2}$ ($X$ a positive real number).  Then
\begin{equation*}  
\mathcal{D}(\mathcal{F}; \phi, w_X) \sim [\widehat{\phi}(0) + \thalf \phi(0)] W_X(\mathcal{F}) \; \; \text{as } X \rightarrow \infty,
\end{equation*}
for $\phi$ with $\text{supp } \widehat{\phi} \subset (- \frac{7}{9}, \frac{7}{9})$.
\end{mytheo}

Note that $W_X(\mathcal{F}) \sim \widehat{w}(0,0)AB$ as $X \rightarrow \infty$ so we are taking about $X^{5/6}$ curves from our family.  In the language of random matrix theory this theorem shows that the family of all elliptic curves has symmetry type $O$, inasmuch as this can be detected without having support outside $(-1, 1)$.
Further, we have
\begin{mytheo}
\label{thm:rankone}
Let $\mathcal{F}$ be the family of elliptic curves given by the Weierstrass equations $E_{a, b}: y^2 = x^3 + ax + b^2$ with $a$ and $b$ positive integers.  Let $w \in C_0^\infty (\mr^{+} \times \mr^{+})$ and set $w_X (E_{a, b}) = w \left( \frac{a}{A}, \frac{b}{B}  \right)$, where $A = X^{1/3}$, $B = X^{1/4}$ ($X$ a positive real number).  Then
\begin{equation*}  
\mathcal{D}(\mathcal{F}; \phi, w_X) \sim [\widehat{\phi}(0) + {\textstyle \frac{3}{2}} \phi(0)]W_X(\mathcal{F}) \; \; \text{as } X \rightarrow \infty, \end{equation*}
for $\phi$ with $\text{supp } \widehat{\phi} \subset (- \frac{23}{48}, \frac{23}{48})$.
\end{mytheo}

Here we are taking about $X^{7/12}$ curves from our family.  These elliptic curves generally have positive algebraic rank (if the point $(0, b)$ is torsion  the Lutz-Nagell criterion implies $b^2 | 4a^3$ so we instantly see that the number of curves in this family such that $(0, b)$ is torsion is $O(X^{1/3 + \varepsilon})$), which explains the presence of the ``extra" $\phi(0)$ contribution.

Obtaining these two theorems with the stated support crucially requires GRH for Dirichlet L-functions.  

S. J. Miller, in his doctoral thesis \cite{Miller}, has independently proved density theorems for various families of elliptic curves (along with other things), but more restricted by the support of $\widehat{\phi}$.  Brumer \cite{Brumer}, Heath-Brown \cite{Heath-Brown}, Michel \cite{Michel}, Silverman \cite{Silverman2}, and others have proved results on the average rank of certain families of elliptic curves that can now be interpreted to essentially be density theorems.  Actually, in order to obtain a density theorem one must asymptotically compute the average of $\log{N}$ over the family; in order to obtain an upper bound on the average rank only a trivial upper bound (such as $N \leq |\Delta|$) is required.

Note that if one can apply $\widehat{\phi}$ of large support then $\phi$ can be localized near the origin, so the zeros $\rho_E = 1 + i\gamma_E$ are held closer to the central point.  Therefore, one challenge for us has been to obtain Theorem \ref{thm:mainresult} with the support of $\widehat{\phi}$ as large as possible.  Breaking support $(-1, 1)$ would be of interest, since it is at that point that the Fourier transforms of the densities of the groups $O$, $SO(\text{even})$, and $SO(\text{odd})$ become distinguishable.  It appears that the current technology is incapable of producing such a result (even assuming GRH), but we have reasons to express the following 
\begin{myconj} 
Theorem \ref{thm:mainresult} holds for test functions $\phi$ with no restrictions on the support.
\label{conj:bigsupport}
\end{myconj}
We will provide some justification for this conjecture in Section \ref{section:evidence} after we prove Theorems \ref{thm:mainresult} and \ref{thm:rankone}.  This conjecture agrees with the folklore conjecture that the root number is equidistributed in the family of all elliptic curves (but see \cite{Helfgott} for an extensive treatment of the variation in sign of the root number).  Miller also predicts the symmetry type is $O$ but from an entirely different direction.  He considers the 2-level density of the family of all elliptic curves and uses conjectures implying equidistribution of root numbers to predict symmetry type $O$ (the 2-level density can distinguish between the various orthogonal symmetry types 
using test functions with arbitrarily small support but it is necessary to know the percentage of curves with given root number).  Our method is quite different and relies on sharp estimates for the three-variable character sum (\ref{eq:S}) for large values of $P$.  It is mysterious how the distribution of the root number (essentially controlled by the M\"{o}bius function of the polynomial $4a^3 + 27b^2$) is captured by such a character sum.

An easy consequence of Theorem \ref{thm:mainresult} is the following

\begin{mycoro}
\label{coro:rankbound}
The family of all elliptic curves ordered by the discriminant as in Theorem \ref{thm:mainresult} has average analytic rank $r \leq 25/14$.  Conjecture \ref{conj:bigsupport} being true implies that $r \leq 1/2$.  
\end{mycoro}

The proof shows that any family with symmetry type $O$ and support up to $\nu$  has average rank bounded by $\onehalf + \frac{1}{\nu}$.  
Applying a density theorem to obtain an upper bound on the average rank in this fashion requires the Riemann Hypothesis for all L-functions in the family.  Using zero density estimates Kowalski and Michel \cite{KM1}, \cite{KM2} obtained an upper bound on the average order of vanishing of all weight $2$ level $q$ modular L-functions, thereby removing the assumption of GRH for their family.  The bound on the average rank is significantly larger than that which is obtained on GRH though.  It would be interesting to obtain an unconditional (yet weaker) upper bound for the average rank of the family of all elliptic curves by mimicing their methods.

Brumer \cite{Brumer} showed that the average rank $r$ of all elliptic curves satisfies $r \leq 2.3$ and Heath-Brown \cite{Heath-Brown} proved $r \leq 2$, modulo a few minor differences in choice of test functions and sieving of unpleasant curves.  Brumer's and Heath-Brown's results requires the GRH for elliptic curve L-functions.  It is a standard conjecture of random matrix theory that, once the symmetry type $G$ for the family has been identified, the density theorem should hold for test function $\widehat{\phi}$ with arbitrarily large support.  Random matrix theory predicts that Conjecture \ref{conj:bigsupport} implies the equidistribution of root numbers (via the 2-level density for instance); it would be interesting to see a direct (number-theoretic) reason for this to be so.

Using that $25/14 = 1.78... < 2$ and the famous theorem that says that if the analytic rank is $\leq 1$ then the algebraic rank and analytic rank are equal and $\Sh$ is finite (due in large part to Kolyvagin \cite{Kolyvagin} and Gross-Zagier \cite{GZ}) we obtain
\begin{mytheo}
Assume GRH.  Then a positive proportion of elliptic curves ordered by height have algebraic rank equal to analytic rank and finite Tate-Shafarevich group.
\end{mytheo}

\subsection{Further results}
\label{section:furtherresults}
Besides the results recorded in the previous section, we investigate similar density results for a variety of interesting families in Sections \ref{section:Z2Z2torsion}-\ref{section:rankonequadratictwists}.  In particular we investigate a number of families with proscribed torsion.  See \cite{Kubert}, Table 3 for the parametrizations of the various torsion structures.  

These torsion families have some surprisingly nice properties.  For instance, the family of curves with torsion group $\mz/4\mz$ is parameterized by $y^2 + xy - by = x^3 - bx^2$ with $\Delta = b^4(1 + 16b)$.  As a general rule we cannot prove a density result for a family where the degree of the discriminant is larger than $3$.  Although the discriminant has degree $5$ for this family, it is easily treated because the irreducible factors $b$ and $1+ 16b$ are linear.  In addition the conductor is much smaller than the discriminant ($N \ll b^2$ whereas $|\Delta| \asymp b^5$) so there are many more curves in the family with conductor $\leq X$ than one would expect based on the degree of the discriminant.  The other torsion families have similar characteristics that make their study extremely pleasant.

We study some popular families of curves with complex multiplication in Section \ref{section:CM}.  These families are rather small yet good results (compared with the sizes of the family) are still obtained because of the simple nature of the Fourier coefficients of the L-functions.  The barrier to obtaining better results with these families is a lack of knowledge on the oscillation of a kind of twisted cubic (or quartic, depending on the family) Gauss sum.

In Section \ref{section:rankonequadratictwists} we study a thin sequence of quadratic twists with positive rank.  This family is challenging because the conductor is essentially a cubic polynomial.  Such families of quadratic twists were studied by Rubin and Silverberg \cite{RS} for example.

\subsection{Structure of the paper}
In Section \ref{section:generalmethodofproof} we set up the machinery for proving a density result for a family of elliptic curves.  We prove Theorems \ref{thm:mainresult} and \ref{thm:rankone} in Sections \ref{section:proof1} and \ref{section:proof2}; some of the more technical details are deferred to the Appendices.

We provide evidence for Conjecture \ref{conj:bigsupport} in Section \ref{section:conjecture}.  This conjecture naturally follows from an assumption that a certain character sum in three variables has square-root cancellation in each variable (one of which is a summation over primes).  We lend credence to the conjecture by studying a sum similar to the aforementioned one but where the summation is extended to integers.  We obtain a stronger result with this new sum; see Theorem \ref{thm:integersum}.

In the remaining sections we prove density results for the interesting families discussed in Section \ref{section:furtherresults}.

\section{General Method of Proof}
\label{section:generalmethodofproof}
In this section we set up some machinery to streamline the proofs of our density theorems.

Suppose we are given a Weierstrass equation (\ref{eq:ellipticcurve}) (not necessarily minimal) which defines an elliptic curve $E$ with conductor $N$ and L-function $L(s, E)$.  To analyze $D(E; \phi)$ we will employ the explicit formula for $L(s, E)$, which for cusp forms of weight two and level $N$ takes the form (see (4.25) in \cite{ILS})
\begin{equation} \label{eq:explicitformula}
D(E; \phi) = \widehat{\phi}(0) \frac{\log N}{\log X} +\onehalf \phi(0)- P(E;\phi) + O\left(\frac{\log \log |\Delta|}{\log X}\right)
\end{equation}
where
\begin{equation*}
P(E;\phi) = \sum_{p > 3} \lambda_E (p) \widehat{\phi}\left(\frac{\log p}{\log X}\right) \frac{2 \log p}{p \log X}.
\end{equation*}
Here $\Delta$ is the discriminant of the curve defined by (\ref{eq:ellipticcurve}) (again, $\Delta$ is not assumed minimal), $X$ is a scaling parameter (any number $\geq 2$ at our disposal), and
\begin{equation}
\label{eq:lambdaf}
\lambda_E (p) = - \sum_{x \shortmod{p}} \left( \frac{x^3 + ax + b}{p} \right) 
\end{equation}
if the Weierstrass equation (\ref{eq:ellipticcurve}) defining $E$ is put into the form (\ref{eq:ellipticcurveb}).  Actually, the sum $P(E;\phi)$ in \cite{ILS} is restricted by $p \notdivtext{N}$ and assumes $\lambda_E(p)$ is the coefficient (\ref{eq:ap}) of $p^{-s}$ in the Dirichlet series expansion of $L(s, E)$.  We claim the modifications in $P(E; \phi)$ are acceptable because the discrepancy in the formula (\ref{eq:explicitformula}) is absorbed by the error term.  To prove this we first note that if $p \negthickspace \not \! | \, \Delta$ and (\ref{eq:ellipticcurveb}) is in global minimal Weierstrass form then $\lambda_E(p)$ is exactly the coefficient (\ref{eq:ap}).  Further, the character sum defining the quantity $\lambda_E (p)$ is left unchanged by a change of variables placing (\ref{eq:ellipticcurveb}) in global minimal Weierstrass form.  Therefore the terms agree for $p \notdivtext{\Delta}$. On the other hand, the terms where $p| \Delta$ are absorbed by the error term. It is worthwhile to note that the error term in (\ref{eq:explicitformula}) is derived by using the Riemann Hypothesis for the symmetric-square L-function $L(s, \text{sym}^2(f)$ to handle terms of the form $\lambda_E(p^2)$.  The formula (\ref{eq:explicitformula}) holds for individual $E$.  In practice one could eliminate the use of the GRH by averaging over families; we have used the Riemann Hypothesis for simplicity and brevity since we are assuming it for other reasons anyways.

Next we sum over the family.  We compute
\begin{equation*}
\mathcal{D}(\mathcal{F}; \phi, w_X) = \widehat{\phi}(0) \sum_{E \in \mathcal{F}} \frac{\log N_E}{\log X} w_X(E) + \onehalf \phi(0) W_X\left(\mathcal{F} \right) - \mathcal{P}(\mathcal{F};\phi, w_X) + O\left(\frac{W_X(\mathcal{F}, \Delta)}{\log X} \right),
\end{equation*}
where
\begin{equation*}
\mathcal{P}(\mathcal{F}; \phi, w_X) = \sum_{E \in \mathcal{F}} P(E; \phi) w_X(E)
\end{equation*}
and
\begin{equation*}
W_X(\mathcal{F}, \Delta) = \sum_{E \in \mathcal{F}} \left| w_X(E) \right| \log \log{|\Delta|}.
\end{equation*}
In general $\log \log{|\Delta|}$ will be $\ll \log{\log{X}}$ so this will be a true error term.

With all of our families we will take $X$ and $w_X$ so that
\begin{equation}
\label{eq:approximateconductor}
\sum_{E \in \mathcal{F}} \frac{\log N_E}{\log X} w_X(E) \sim W_X\left( \mathcal{F} \right) \; \; \text{as $X \rightarrow \infty$}.
\end{equation}
It is generally highly nontrivial to prove such an asymptotic holds; often it amounts to having control on the square divisors of a polynomial of high degree.  This is a significant barrier to producing density theorems with families of high rank.  We call (\ref{eq:approximateconductor}) the {\it conductor condition} for the family $\mathcal{F}$.

The idea in most cases is to approximate $N_E$ with numbers $R_E$ which have the same prime divisors as $N_E$ but are easier to compute (often $R_E = |\Delta|$, the discriminant of the elliptic curve).  Then $X$ and $w_X$ will be chosen so that $X$ is approximately $R_E$ for $E$ in the support of $w_X$.  For an example, with the family $y^2 = x^3 + ax + b$ given in Theorem \ref{thm:mainresult} we have $w_X$ scaled so that $a \asymp X^{1/3}$, $b \asymp X^{1/2}$, $|\Delta| \asymp X^{5/6}$, and we take $R_E = |\Delta|$.

The following lemma will streamline many of our arguments to show (\ref{eq:approximateconductor}) holds.
\begin{mylemma}
\label{lem:lemma3a}
Let $\mathcal{F} = \{E_d\}$ be a family of elliptic curves with notation as in Section \ref{section:preliminaries}.  Let $\Delta(d)$ be the discriminant of the curve $E_d$. Suppose that there exists an integer polynomial $R(d)$ dividing $\Delta(d)$ such that each irreducible factor of $\Delta(d)$ divides $R(d)$.  Further suppose
\begin{equation}
\label{eq:suminconductorlemma}
 \sum_{E \in \mathcal{F}}  \mathop{\sum_{p^{\alpha} || R_E}}_{\alpha > 0} \left| w_X(E) \right| \log{p^{\alpha - 1}} + \sum_{E \in \mathcal{F}} \mathop{\mathop{ \sum_{p || R_E}}_{p^2 | N_E}}_{p > 3} \left| w_X(E) \right| \log{p}   \ll \left|W_X\left(\mathcal{F}\right) \right|
\end{equation}
uniformly in $X$, where we have defined $R_E = R(d)$ for $E = E_d$.
Suppose $R(d) \asymp X$ for all $f$ in the support of $w_X$.  Then we have
\begin{equation*}
\sum_{E \in \mathcal{F}} \frac{\log N_E}{\log X}  w_X(E) = W_X\left( \mathcal{F} \right) + O\left(\frac{ \left|W_X\left(\mathcal{F}\right) \right| }{\log X} \right).
\end{equation*}
\end{mylemma}
\begin{myrema} Both sums in (\ref{eq:suminconductorlemma}) have an interpretation.  For the first sum to be small there must not be many large square divisors of $R_E$ (i.e. $R_E$ is not too large).  For the second sum to be small it must be rare for a prime to divide $R_E$ to lower order than $N_E$ (i.e. $R_E$ is not too small).  In most applications $N_E | R_E$ so the second sum will be void.  In general we can handle the first sum as long as all of the irreducible factors of $R(d)$ are of degree $3$ or less.  
\end{myrema}

\begin{proof}  
Suppose $R(d)$ is as above. Then
\begin{equation*}
W_X \left( \mathcal{F} \right) - \sum_{E \in \mathcal{F}} \frac{\log N_E}{\log X} w_X(E) =  \sum_{E \in \mathcal{F}} \frac{\log (X/R_E)}{\log X} w_X(E) + \sum_{E \in \mathcal{F}} \frac{\log (R_E/N_E)}{\log X} w_X(E).
\end{equation*}
Since we are assuming $R_E \asymp X$ the first sum is trivially $\ll \left| W_X(\mathcal{F}) \right| (\log{X})^{-1}$.  The second sum is
\begin{equation*}
\frac{1}{\log{X}} \sum_{E \in \mathcal{F}} \mathop{\sum_{p^\alpha || R_E}}_{p^\beta || N_E} w_X(E) \log p^{\alpha - \beta},
\end{equation*}
by using the additivity of the logarithm to separate the prime factors of $R_E$ and $N_E$ (these primes $p$ have nothing to do with the explicit formula of course).
First consider the terms with $\beta > 0$.  By taking the terms with $\alpha > \beta$ and $\alpha < \beta$ separately it's clear the contribution is
\begin{equation*}
\ll \frac{1}{\log{X}} \sum_{E \in \mathcal{F}} \mathop{\sum_{p^\alpha || R_E}}_{\alpha > 0} \left|w_X(E) \right|  \log p^{\alpha - 1} + \frac{1}{\log{X}}  \sum_{E \in \mathcal{F}} \mathop{\mathop{\sum_{p || R_E}}_{p^2 | N_E}}_{p > 3} \left|w_X(E) \right| \log p \text{ }  + O\left(\frac{\left|W_X(\mathcal{F})\right|}{\log{X}} \right),
\end{equation*}
the error term coming from $p = 2$ or $3$ and $\alpha < \beta$ (in which case $\beta$ may be larger than 2).

Now consider those terms with $\beta = 0$ (and hence $\alpha > 0$).  Since $p | \Delta$ but $p \notdivtext N$ this implies the equation (\ref{eq:ellipticcurve}) defined by $d$ is not minimal at $p$ and therefore $\alpha \geq 12$, so we may absorb these terms into the first sum which consists of $p^\alpha || R_E, \alpha >0$.
\end{proof}

For any individual family we will need to estimate $\mathcal{P}(\mathcal{F} ;\phi, w_X)$ using (ad hoc) techniques applicable to the family.
In its evaluation we will often make use of the following identity (Poisson Summation $\mymod{l}$).
\begin{myprop} Let $w$ be a Schwartz-class function, $D$ a positive real number, and $a$ an integer.  Then the following holds
\begin{equation}
\label{eq:poissonsummation}
\mathop{\sum_{d \in \mz }}_{d \equiv a \shortmod{l} } w\left( \frac{d}{D} \right) = \frac{D}{l} \sum_{h \in \mz} e\left(\frac{ha}{l} \right) \widehat{w}\left( \frac{hD}{l} \right).
\end{equation}
\end{myprop}

\section{Proof of Theorem \ref{thm:mainresult}}
\label{section:proof1}
By the discussion in Section \ref{section:generalmethodofproof}, to prove Theorem \ref{thm:mainresult} we need to show
\begin{equation*}
\mathcal{P}(\mathcal{F}; \phi, w_X) \ll \frac{W_X(\mathcal{F})}{\log{X}} \asymp \frac{AB}{\log{X}}
\end{equation*}
and to show the conductor condition (\ref{eq:approximateconductor}) holds.  
\subsection{The Conductor Condition}
In this section we show  (\ref{eq:approximateconductor}) holds.  

\begin{mylemma}  
\label{lem:4conductor}
Let $\mathcal{F}, A, B$, and $w$ be as in Theorem \ref{thm:mainresult}.  Then we have
\begin{equation*}
\sum_{E \in \mathcal{F}} \frac{\log N_E}{\log X} w_X(E) = \left\{1 + O\left(\frac{1}{\log{X}} \right) \right\} W_X\left( \mathcal{F} \right).
\end{equation*}
\end{mylemma}
\begin{proof}
We apply Lemma \ref{lem:lemma3a}.  We take the polynomial $R(d)$ for $d = (a, b)$ to be given by $R(d) = 16(4a^3 + 27b^2)$. Note $R(d) =  |\Delta(d)|$.  Recall that $A = X^{1/3}$ and $B = X^{1/2}$ so that $R(d) \asymp X$ since $a \asymp A$ and $b \asymp B$.  The sum
\begin{equation*}
\sum_{E \in \mathcal{F}} \mathop{ \sum_{p || R_E}}_{p^2 | N_E} w_X(E)  \log{p}  
\end{equation*}
is empty since $N | \Delta$.  The other term is
\begin{equation*}
\sum_{E \in \mathcal{F}}  \mathop{\sum_{p^{\alpha} || R_E}}_{\alpha > 0} w_X(E) \log{p^{\alpha - 1}}  = \sum_{a} \sum_{b} \mathop{\sum_{p^{\alpha} || 16(4a^3 + 27b^2)}}_{\alpha > 0} w\left( \frac{a}{A}, \frac{b}{B} \right) \log{p^{\alpha - 1}}
\end{equation*}
by definition.
First suppose $p > 3$.   Interchange the order of summation and for each prime $p$ define $\gamma$ by $p^\gamma || a$.  We first consider the terms where $\alpha < 3 \gamma$. For these cases $p^\alpha || b^2$ (so $\alpha$ is necessarily even).  We employ the change of variables $a = p^\gamma a'$, $b = p^{\alpha/2} b'$ and obtain (we often use primes to indicate that a summation is carried out with certain coprimality conditions in place which should be apparent from context; in the next formula the restriction is $(a', p) = (b', p) = 1$)
\begin{equation*}
\mathop{\sum_{p^\alpha \ll X}}_{\alpha > 0} \sum_{\alpha/3 < \gamma \ll \log{A}} \sum_{a'}{}^{'} \sum_{b'}{}^{'}  w\left( \frac{p^\gamma a'}{A}, \frac{p^{\alpha/2} b'}{B} \right) \log{p^{\alpha - 1}}
\end{equation*}
\begin{equation*}
\ll \mathop{\sum_{p^\alpha \ll X}}_{\alpha > 0} \sum_{\alpha/3 < \gamma \ll \log{A}} \left(1 + \frac{A}{p^\gamma} \right)\left(1 + \frac{B}{p^{\alpha/2}} \right) \log{p^{\alpha - 1}} 
\end{equation*}
\begin{equation*}
\ll \mathop{\sum_{p^\alpha \ll X}}_{\alpha > 0}\left(\log{A} + \frac{A}{p^{[\alpha/3] + 1}} + \frac{B\log{A}}{p^{\alpha/2}} + \frac{AB}{p^{\alpha/2 + [\alpha/3] + 1}} \right) \log{p^{\alpha - 1}}.
\end{equation*}
This sum is
\begin{equation*}
\ll \sqrt{X} \log{A} + A \log{X} + B (\log{X})^2 + AB \ll  X^{5/6},
\end{equation*}
because
\begin{equation*}
\mathop{\sum_{p^\alpha \leq Z}}_{\alpha > 0} \log{p^{\alpha - 1}} \ll  \sqrt{Z}, \; \; \; \mbox{and}
\; \; \;
\mathop{\sum_{p^\alpha}}_{\alpha > r} \frac{\log{p^{\alpha - 1}}}{p^{\alpha/r}} \ll 1.
\end{equation*}
Now we consider the terms with $\alpha \geq 3\gamma$.  In this case $p^{3 \gamma} || b^2$ (so $\gamma$ is necessarily even).  We employ the change of variables $a \rightarrow p^\gamma a'$, $b \rightarrow p^{3\gamma/2} b'$ and obtain
\begin{equation*}
\mathop{\sum_{p^{\alpha} \ll X}}_{\alpha > 0} \sum_{0 \leq \gamma < \alpha/3} \mathop{\sum_{a'}{}^{'} \sum_{b'}{}^{'}}_{p^{\alpha - 3\gamma} || (4a'^3 + 27b'^2) } w\left( \frac{p^\gamma a'}{A}, \frac{p^{3\gamma/2} b'}{B} \right) \log{p^{\alpha - 1}} .
\end{equation*}
We split the summation over $b'$ into progressions $\mymod{p^{\alpha - 3\gamma}}$.  Since $(p, a') = (p, b') = 1$ we get that for each $a'$ the number of solutions (in $b' \pmod{p}$) to the congruence $4a'^3 + 27b'^2 \equiv 0 \pmod{p^{\alpha - 3\gamma}}$ is bounded by 2.  Therefore we obtain
\begin{equation*}
\mathop{\sum_{p^{\alpha} \ll X}}_{\alpha > 0} \sum_{0 \leq \gamma < \alpha/3}  \left( 1 + \frac{A}{p^\gamma} \right) \left( 1 + \frac{B p^{-3\gamma/2}}{p^{\alpha - 3\gamma}} \right) \log{p^{\alpha - 1}},
\end{equation*}
\begin{equation*}
\ll \mathop{\sum_{p^{\alpha} \ll X}}_{\alpha > 0} \left(\alpha + A + \frac{B}{ p^{\alpha  - \frac{3}{2}[\alpha/3]}  }  + 
\frac{AB}{p^{\alpha - \onehalf[\alpha/3] }} 
\right) \log{p^{\alpha - 1}}
\end{equation*}
which is bounded by $X^{5/6}$ by the same type of reasoning used for the case $\alpha > 3 \gamma$ (check $\alpha = 2$ and $\alpha \geq 3$ separately).  

The primes 2 and 3 are handled in the same way as above with minor changes.  For $p = 3$ we may assume $\alpha \geq 4$ by trivially estimating the terms with $\alpha \leq 3$.  The estimation for $\alpha < 3\gamma$ is the same as before after the change of variables $\alpha \rightarrow \alpha + 3$.  The estimation for $\alpha \geq 3\gamma$ is as before after the change of variables $\gamma \rightarrow \gamma + 1$.  For $p = 2$ we may assume $\alpha \geq 6$.  After cancelling $2^6$ in the discriminant congruence we are left in a very similar case to $p = 3$.  We omit the tedious yet elementary details. 
\end{proof}

\subsection{The Central Estimate}
We will have proved Theorem \ref{thm:mainresult} once we have proved
\begin{mylemma} \label{lem:5b} Set $A=X^{1/3}$ and $B = X^{1/2}$.  Then
\begin{equation*}
\sum_{a} \sum_{b} P(E;\phi) w\left(\frac{a}{A}, \frac{b}{B}\right) \ll \frac{X^{5/6}}{\log{X}} 
\end{equation*}
provided supp ${\widehat{\phi}} \subset (- \frac{7}{9}, \frac{7}{9})$.
\end{mylemma}
This lemma is the heart of the matter.
\begin{proof}  
We calculate
\begin{equation*}
\sum_{a} \sum_{b} P(E;\phi) w\left(\frac{a}{A}, \frac{b}{B}\right)
= \sum_{p > 3} \frac{2 \log p }{p \log X} \widehat{\phi}\left( \frac{\log p}{\log X} \right) \sum_a \sum_b \lambda_{a, b} (p) w\left(\frac{a}{A}, \frac{b}{B}\right).
\end{equation*}
Apply Poisson summation $\mymod{p}$ in the summation over $a$ and $b$ and obtain
\begin{equation*}
\sum_a \sum_b \lambda_{a, b} (p) w\left(\frac{a}{A}, \frac{b}{B}\right) = \frac{AB}{p^2} \sum_h \sum_k \mathop{ \mathop{\sum \sum}_{\alpha \shortmod{p} }}_{\beta \shortmod{p}}  \lambda_{\alpha, \beta} \, e \left( \frac{\alpha h + \beta k}{p} \right) \widehat{w} \left(\frac{hA}{p}, \frac{kB}{p}\right).
\end{equation*}
The summation over $\alpha$ and $\beta$ is evaluated in Section \ref{section:charactersums} (the evaluation is completely straightforward).
Using Lemma \ref{lem:2a} we continue, obtaining
\begin{equation}
\label{eq:goodform}
-\frac{AB}{\log{X}} \sum_{p >3}  \varepsilon_p  \frac{2 \log{p}}{p^{3/2}}\widehat{\phi}\left( \frac{\log p}{\log X} \right) \sum_h \sum_k \left( \frac{k}{p} \right) e\left( \frac{-h^3 \bar{k}^2}{p} \right) \widehat{w} \left(\frac{hA}{p}, \frac{kB}{p} \right).
\end{equation}
We remark at this point that if we estimate this sum trivially we get a bound of the order (summing $p$ up to $P$, say)
\begin{equation*}
\frac{AB}{\log{X}} \sum_{p \leq P}  \frac{ \log{p}}{p^{3/2}} \left(1 + \frac{p}{A} \right)\left(1 + \frac{p}{B} \right) \ll (AB + BP^{1/2} + P^{3/2}) (\log{X})^{-1},
\end{equation*}
which is $O\left(AB (\log{X})^{-1} \right)$ when $P \leq X^{5/9}$.  Brumer essentially obtained this result \cite{Brumer}.  To get larger support we need to prove there is quite a lot of cancellation in the three variable character sum (\ref{eq:goodform}).  By exploiting some cancellation in this sum Heath-Brown \cite{Heath-Brown} has improved Brumer's result to the support range $(-2/3, 2/3)$.  Note that any improvement on Heath-Brown's result shows that the average rank is strictly less than two. 

The first step is to eliminate the variation in $\varepsilon_p$.  To do so we sum separately over the progressions $p \equiv 1 \pmod{4}$ and $p \equiv 3 \pmod{4}$ .  Effectively, it suffices to replace $\varepsilon_p$ by $\psi_4(p)$, a Dirichlet character $\mymod{4}$.  Now break up the summation in (\ref{eq:goodform}) into dyadic segments using a smooth partition of unity.  It suffices to consider sums of the form
\begin{equation}
\label{eq:dyadic}
\mathop{\mathop{ \mathop{\sum \sum \sum}_{H \leq h < 2H} }_{K \leq k < 2K} }_{P \leq p < 2P}
\frac{\log{p}}{p^{3/2}} \psi_4(p) \left( \frac{k}{p} \right) e\left( \frac{-h^3 \bar{k}^2}{p} \right) \widehat{\phi}\left( \frac{\log p}{\log X} \right) \widehat{w} \left(\frac{hA}{p}, \frac{kB}{p} \right) g(h, k, p),
\end{equation}
where $g$ is a smooth compactly supported function arising from the partition of unity.  We assume that the restrictions on $p$, $h$, and $k$ are redundant, following from the support of $g$.  It suffices to show that every sum of type (\ref{eq:dyadic}) is $\ll X^{-\varepsilon}$, with the implied constant depending only on $w, ||g||_{\infty}, ||g'||_{\infty},$ etc.  In addition we have to account for the contribution to (\ref{eq:goodform}) of $h=0$, but this contribution is negligible by trivial estimations.

Let $S(H, K, P)$ be the sum given by (\ref{eq:dyadic}) (in the notation we suppress the dependence on the test functions).



Using the bound $\widehat{w}(x, y) \ll_M (1 + |x|)^{-M} (1 + |y|)^{-M}$ we may assume that $H \ll_\varepsilon (P/A)^{1 + \varepsilon}$ and $K \ll_\varepsilon (P/B)^{1 + \varepsilon}$.  

It will be necessary to use different techniques of estimation in different ranges.  As a first step, we have the bound
\begin{equation}
\label{eq:Weylsbound}
\sum_{P \leq p < 2P} \left| \sum_{H \leq h < 2H} e \left( \frac{h^3 \bar{k}^2}{p} \right) \right| \ll \left(H^{3/4} P + H P^{3/4} +  H^{1/4} P^{5/4}\right) (HKP)^{\varepsilon}.
\end{equation}
The proof is standard by Weyl's method.
If we apply this to (\ref{eq:dyadic}) (after using partial summation to separate the variables) we get the bound
\begin{equation}
\label{eq:Weylsbound2}
S(H, K, P) \ll \left(H^{3/4} K P^{-1/2} + HK P^{-3/4} + H^{1/4} K P^{-1/4}  \right) X^{\varepsilon},
\end{equation}
which is useful in some ranges.

To cover the ranges where (\ref{eq:Weylsbound2}) is insufficient we continue with (\ref{eq:dyadic}).  

We prefer to sum over relatively prime $h$ and $k$ so we define $d = (h^3, k^2)$ and let $d_0$ to be the least positive integer such that $d | d_0^3$.  Since $d_0 | h$ we may set $h = d_0 h_0$.  
The condition $(h^3, k^2)  = d$ is equivalent to $(h_0, k^2/d) = 1$ and $(d_0^3/d, k^2/d) = 1$.  Then
\begin{equation*}
\sum_{h} \sum_{k} e\left( \frac{-h^3 \bar{k}^2}{p} \right) = \sum_k \mathop{\mathop{\sum_{d | k^2} \sum_{h_0}}_{(h_0, k^2/d) = 1}}_{(d_0^3/d, k^2/d) = 1}  e\left( \frac{-h_0^3 (d_0^3/d) \overline{(k^2/d)}}{p} \right).
\end{equation*}

The point is that everything in the exponential is coprime with $k^2/d$ (besides possibly $\overline{k^2/d}$) .  Now we may employ the following elementary reciprocity formula
\begin{equation}
\label{eq:elementaryreciprocity}
\frac{\bar{u}}{v} + \frac{\bar{v}}{u} \equiv \frac{1}{uv} \pmod{1},
\end{equation}
where $u, v, \bar{u}, \bar{v}$ are integers such that $(u,v) =1, u\bar{u} \equiv 1 \pmod{v}$, and $v\bar{v} \equiv 1 \pmod{u}$.  This reciprocity law was also employed by \cite{Heath-Brown} in his work on this problem.  In our application $u = k^2/d$ and $v = p$.  Now 
$S(H, K, P)$ has transformed into
\begin{equation}
\label{eq:goodestimate}
\mathop{ \mathop{\sum \sum}_{K \leq k < 2K}}_{P \leq p < 2P} \sum_{d | k^2}
\mathop{\mathop{  \sum_{\frac{H}{d_0} \leq h_0 < 2\frac{H}{d_0}}}_{  (h_0, \frac{k^2}{d}) =1} }_{(\frac{d_0^3}{d}, \frac{k^2}{d}) = 1}  \left( \frac{k}{p} \right) \psi_4(p) e\left( \frac{h_0^3 (d_0^3/d)\bar{p}}{k^2/d} \right) e\left( \frac{-h_0^3d_0^3}{p k^2} \right) U(h_0, d_0, k, p),
\end{equation} 
where
\begin{equation*}
U(h_0, d_0, k, p) =  g(h_0 d_0, k, p) \widehat{w} \left(\frac{h_0 d_0 A}{p}, \frac{kB}{p} \right) \frac{2  \log{p}}{p^{3/2}\log{X}} 
\widehat{\phi}\left( \frac{\log{p}}{\log X} \right).
\end{equation*}
To separate the variables we use the following expansion of additive characters into multiplicative characters via Gauss sums
\begin{equation}
\label{eq:gaussexpansion}
e\left(\frac{h_0^3 (d_0^3/d) \bar{p}}{k^2/d}\right) = \frac{1}{\phi(k^2/d)}\sum_{\chi \shortmod{k^2/d}} \tau(\chi) \overline{\chi}(h_0^3 (d_0^3/d) \bar{p}),
\end{equation}
(valid because of the coprimality conditions), obtaining the identity
\begin{equation}
\label{eq:expandedsum}
S(H, K, P) = \sum_{K \leq k < 2K}  \sum_{d | k^2} \frac{1}{\phi(k^2/d)} \sum_{\chi \shortmod{k^2/d}} \tau(\chi)  \overline{\chi}(d_0^3/d) Q(d, k, \chi)
\end{equation}
where
\begin{equation*}
Q(d, k, \chi) = \sum_{P \leq p <  2P}  \sum_{\frac{H}{d_0} \leq h_0 < 2\frac{H}{d_0}}  \psi_4(p) \chi(p) \left( \frac{k}{p}\right) \overline{\chi}^3(h_0) e\left( \frac{-h_0^3d_0^3}{p k^2} \right) U(h_0, d_0, k, p).
\end{equation*}
Our goal is to get a good bound for $Q(d, k, \chi)$ and estimate the rest trivially.

We wish to apply Lemma \ref{lem:4a} to $Q(d, k, \chi)$.  To do so we must define a number of parameters and check that the conditions of Lemma \ref{lem:4a} are satisfied.  We set
\begin{equation}
\label{eq:4test}
F(u, v) =  \left( \frac{v}{P} \right)^{-3/2} g(u d_0, k, v) \widehat{w} \left(\frac{u d_0 A}{v}, \frac{kB}{v} \right) \widehat{\phi}\left( \frac{\log v}{\log X} \right).
\end{equation}
The character $\chi(u)$ in Lemma \ref{lem:4a} is replaced by $\overline{\chi}^3(u)$ (with modulus $l_1 = \frac{k^2}{d}$), while the character $\psi(v)$ in Lemma \ref{lem:4a} is $\chi(v) \psi_4(v) (k/v)$ (with modulus $l_2=\text{lcm}(4,k^*, \frac{k^2}{d})$, where $k^*$ is the conductor of $\left(\frac{k}{v}\right)$), $q = -\frac{k^2}{d_0^3}$, $U = \text{min}\left\{\frac{H}{d_0}, \frac{P}{d_0A} \right\}$, and $V= P$.  We have
\begin{myclaim}
\label{lem:4g}
The test function $F(u, v)$ defined by (\ref{eq:4test}) satisfies the conditions of Lemma \ref{lem:4a} with $U = \text{min}\left\{H/d_0, P/d_0A\right\}$ and $V = P$.  
\end{myclaim}
\begin{proof}[Proof of claim]  The proof is a straightforward calculation.  The only slightly thorny issue is that differentiation of $\widehat{w}$ with respect to $v$ introduces a factor $(ud_0A/v)^\alpha$.  This can be absorbed (for $M$ large enough) by the bound of $(1 + u/(v/d_0A))^{-M}$ on any partial derivative of $\widehat{w}$.  
\end{proof}

To apply Lemma \ref{lem:4a} to $Q(d, k, \chi)$ we must first extend the summation over $p$ to prime powers instead of just primes.  It's easy to see by trivial estimations that we can extend the summation to prime powers without changing the bound for $Q(d, k, \chi)$.  Therefore we have
\begin{mycoro}  If $\chi \psi_4 (k/\cdot)$ and $\chi^3$ are nonprincipal then
\label{coro:Qbound}
\begin{equation*}
Q(d, k, \chi) \ll P^{-1} \left(\frac{H}{d_0}\right)^{1/2} \left( 1 + \frac{H^{3/2}}{P^{1/2} k} \right) X^\varepsilon.
\end{equation*}
If $\chi \psi_4 (k/\cdot)$ is principal but not $\chi^3$ then we lose a factor $P^{1/2}$.  If $\chi^3$ is principal but not $\chi \psi_4 (k/\cdot)$ then we lose a factor $(H/d_0)^{1/2}$.
\end{mycoro}



Having obtained the bound for $Q(d, k, \chi)$ we apply it to $S(H, K, P)$.  We have four cases according to which characters are principal.  Let $S = S_1 + S_2 + S_3 + S_4$, where $S_1$ corresponds to the terms where both characters are nonprincipal, $S_2$ corresponds to the terms where $\chi^3$ is nonprincipal but $\chi \psi_4 (k/\cdot)$ is not, $S_3$ corresponds to the terms where both characters are principal, and $S_4$ corresponds to the remaining terms where $\chi = \psi_4 (k/\cdot)$ is nonprincipal.

{\bf Case 1}.  To bound the sum $S_1$ we apply Corollary \ref{coro:Qbound} to (\ref{eq:expandedsum}), obtaining the bound 
\begin{equation*}
S_1 \ll \sum_{K \leq k < 2K} \sum_{d | k^2} \frac{1}{\phi(k^2/d)}  \sum_{\chi \shortmod{k^2/d}} P^{-1} \left(\frac{H}{d_0}\right)^{1/2} \left( 1 + \frac{H^{3/2}}{P^{1/2} k} \right) X^\varepsilon\left| \tau(\chi) \right|
\end{equation*}
\begin{equation*}
\ll H^{1/2}  P^{-1} \left(1 +  \frac{H^{3/2}}{P^{1/2} K} \right) X^{\varepsilon} \sum_{K \leq k < 2K} k \sum_{d | k^2} \frac{1}{(dd_0)^{1/2}}  
\end{equation*}
\begin{equation*}
\ll H^{1/2} P^{-1} \left( K^2 + H^{3/2}K P^{-1/2} \right) X^{\varepsilon}.
\end{equation*}
We require this bound to be $\ll X^{-\varepsilon}$.  Using $H \ll (P/A)^{1 + \varepsilon}$ and $K \ll (P/B)^{1 + \varepsilon}$ shows the requirement is $P \ll X^{7/9 - \varepsilon}$.  The existence of $7/9$ here exhibits the limit of our method.  The cases where one or both of the characters are principal are tedious to carry out but do not pose a significant barrier to obtaining larger support.

It remains to bound the sum (\ref{eq:expandedsum}) when one or both of the characters are principal. 
The loss of cancellation in these cases will be made up for by the rarity of principal characters.  First, $\overline{\chi}^3$ is trivial for $\ll k^\varepsilon$ characters $\chi$.  The character $\chi \psi_4 (k/\cdot)$ is trivial for only the character $\chi = \psi_4 (k/\cdot)\chi_0$.  Therefore the only way both characters are trivial is if $\chi$ is trivial and $(k/\cdot) \psi_4$ is trivial (which limits $k$ to be a square).  

{\bf Case 2}.  In case $\overline{\chi}^3$ is trivial but not $\chi \psi_4 (k/\cdot)$ we lose $(H/d_0)^{1/2}$ (from the loss of cancellation) but save $\phi(k^2/d)$ (from the rarity of such characters), which gives the bound
\begin{equation*}
S_2 \ll H P^{-1}\left(1 + \frac{H^{3/2}}{P^{1/2} K} \right)  X^{\varepsilon} \sum_k k^{-1}  \sum_{d | k^2} \frac{\sqrt{d}}{d_0}
\end{equation*}
\begin{equation*}
\ll H P^{-1} \left(1 + \frac{H^{3/2}}{P^{1/2} K} \right) K^{1/3} X^{\varepsilon} 
\end{equation*}
using only the obvious bound $d_0^{-1} \leq d^{-1/3}$ (from $d | d_0^3$).  Using $H \ll (P/A)^{1 + \varepsilon}$ and $1 \ll K \ll (P/B)^{1 + \varepsilon}$ and requiring this $\ll X^{-\varepsilon}$ means we must require $P \ll X^{5/6 - \varepsilon}$.

{\bf Case 3.}  In case both characters are trivial we know $k$ is a square and that $|\tau(\chi)| \leq 1$.  Using (\ref{eq:expandedsum}) and bounding $Q(d, k, \chi)$ trivially by $PH/d_0$ we easily get the bound of
\begin{equation*}
S_3 \ll H P^{-1/2} X^{\varepsilon} \mathop{\sum_{k = \square}}_{K \leq k < 2K} k^{-2} \sum_{d | k^2} \frac{d}{d_0} = HP^{-1/2} X^{\varepsilon} \sum_{(K)^{1/2} \leq l < (2K)^{1/2}} l^{-4} \sum_{d | l^4} \frac{d}{d_0}. 
\end{equation*}
Notice that if $a | d$ then $a/a_0 \leq d/d_0$ (look at each prime separately).  Therefore we have
\begin{eqnarray*}
S_3 & \ll & HP^{-1/2} X^{\varepsilon} \sum_{K^{1/2} \leq l < (2K)^{1/2}}  \prod_{p | l} p^{-2} \\
& \ll & H P^{-1/2} K^{-1/2} X^{\varepsilon},
\end{eqnarray*}
which is $\ll X^{-\varepsilon}$ when $K \gg P X^{-2/3 + \varepsilon}$.  To handle $K \ll P X^{-2/3 + \varepsilon}$ we simply apply Weyl's bound (\ref{eq:Weylsbound}) when $k = \square$, $K \leq k < 2K$.  We thus obtain the bound
\begin{equation*}
S_3 \ll \left(H^{3/4} P^{-1/2} + H P^{-3/4} + H^{1/4} P^{-1/4}  \right) K^{1/2} X^{\varepsilon},
\end{equation*}
which is $\ll X^{- \varepsilon}$ (using $H \ll (P/A)^{1 + \varepsilon}$ and $K \ll P X^{-2/3 + \varepsilon})$ when $P \ll X^{7/9 - \varepsilon}$.  We expect that use of current technology would allow us to take $P$ larger than $X^{7/9}$ here but since we are restricted to $7/9$ elsewhere we do not pursue such a result.

{\bf Case 4.}  In the last case where $\chi \psi_4 (k/\cdot)$ is trivial but $\overline{\chi}^3$ is not we lose $P^{1/2}$ but save $\phi(k^2/d)$.  We get further saving by noticing that $\psi_4 (k/\cdot)$ has conductor $k^*$ equal to the square-free part of $k$ (up to a factor $2$ or $4$).  Therefore $k^* | (k^2 d^{-1})$ and hence $d|(k^2 (k^*)^{-1})$.  We get the bound
\begin{equation*}
S_4 \ll H^{1/2} P^{-1/2} \left( 1+ \frac{H^{3/2}}{P^{1/2} K} \right) X^{\varepsilon}  \sum_{K \leq k < 2K}  k^{-1} \sum_{d |(k^2(k^*)^{-1})} \left(\frac{d}{d_0}\right)^{1/2}.
\end{equation*}
As before, $(d/d_0)^{1/2}$ is increasing with respect to divisibility so we get the bound
\begin{equation*}
S_4 \ll H^{1/2} P^{-1/2} \left( 1+ \frac{H^{3/2}}{P^{1/2} K} \right) X^{\varepsilon} \sum_{K \leq k < 2K}  k^{-1} \prod_{p^2 |k} p.
\end{equation*}
\begin{equation*}
\ll H^{1/2} P^{-1/2}X^{\varepsilon}  \left( 1+ \frac{H^{3/2}}{P^{1/2} K} \right).
\end{equation*}

For this bound to be $\ll X^{-\varepsilon}$ it is necessary and sufficient that 
\begin{equation}
\label{eq:equation18}
H^2 K^{-1} \ll P X^{-\varepsilon}.  
\end{equation}

We use a different method to estimate $Q(d, k, \chi)$ which works well for small $k$.  In our current case we have
\begin{equation*}
Q(d, k, \chi) = \sum_{P \leq p < 2P}{}^{'} \sum_{H/d_0 \leq h_0 < 2H/d_0} e \left(\frac{-h_0^3d_0^3}{pk^2} \right) \bar{\chi}^3(h_0) U(h_0, d_0, k, p).
\end{equation*}
We exploit cancellation in $\sum_{h_0} c_{h_{0}} e(-h_0^3 d_0^3/pk^2)$ with arbitrary complex coefficients $c_n$.  Precisely, we use the following
\begin{mylemma}  \label{lem:lemma4x} Let $c_n$ be complex numbers satisfying $|c_n| \leq 1$ and let
\begin{equation*}
R(N, P, d_0, k) = \sum_{P \leq p < 2P} \left| \sum_{N \leq n < 2N} e \left(\frac{n^3 d_0^3}{pk^2} \right) c_n \right|.
\end{equation*}
Then
\begin{equation*}
R(N, P, d_0, k) \ll N^{1/2}P + N^{1/4} P^{5/4} k^{1/2} d_0^{-3/4}.
\end{equation*}
\end{mylemma}
This is a special case of a more general result which is stated and proved in Appendix \ref{section:appendixb}; specifically, we apply Lemma \ref{lem:appendixb} with $f(x) = x^3$, $g(x) = x^{-1}$, and $Y = N^3 d_0^3 P^{-1} k^{-2}$.

Applying Lemma \ref{lem:lemma4x} to $Q(d, k, \chi)$ via partial summation gives (actually one must separate the variables $h_0$ and $p$ in $U(h_0d_0, k, p)$ before applying Lemma \ref{lem:lemma4x}, which can be done in any standard way with no cost)
\begin{mycoro} If $\chi \psi_4 (k/\cdot)$ is principal and $\overline{\chi}^3$ is non-principal then 
\label{cor:S}
\begin{equation*}
Q(d, k, \chi) \ll \left(\frac{H}{d_0}\right)^{1/2} P^{-1/2} +  P^{-1/4} H^{1/4} k^{1/2} d_0^{-1} X^{\varepsilon}. 
\end{equation*}
\end{mycoro}
Applying Corollary \ref{cor:S} to $S_4$ gives the bound 
\begin{equation*}
S_4 \ll \sum_{K \leq k < 2K} \sum_{d | k^2} d^{1/2} k^{-1} \left( \frac{H^{1/2}}{d_0^{1/2} P^{1/2}} + \frac{ H^{1/4} k^{1/2}   }{P^{1/4} d_0}X^{\varepsilon} \right).
\end{equation*}
Using the same techniques as before to estimate the sum over $d$ gives the bound
\begin{equation*}
S_4 \ll \frac{H^{1/2}K^{1/2} }{P^{1/2}} X^{\varepsilon} + \frac{H^{1/4}K^{1/2}}{P^{1/4}} X^{\varepsilon}. 
\end{equation*}
Using only $H \ll (P/A)^{1 + \varepsilon}$ and $K \ll (P/B)^{1 + \varepsilon}$ shows the first term is $\ll X^{-\varepsilon}$ when $P \ll X^{5/6 - \varepsilon}$.  
For the second term in the above bound to be sufficient we must have 
$K \ll P^{1/2} H^{-1/2} X^{-\varepsilon} (*)$.  Assuming (\ref{eq:equation18}) does not hold (i.e. $H^2 K^{-1} \gg P X^{- \varepsilon}$) and using $H \ll (P/A)^{1 + \varepsilon}$ we see that (*) 
holds when $P \ll X^{5/6 - \varepsilon}$.  Having considered all possible cases the proof is complete. 
\end{proof}

\subsection{A Complete Character Sum}
\label{section:charactersums}
In this section we evaluate a character sum which arose in our averaging.  Set
\begin{equation*}
\lambda_{a, b} (p) = - \sum_{x \shortmod{p}} \left( \frac{x^3 + ax + b}{p} \right).
\end{equation*}
Note that this is $\lambda_E(p)$ for the elliptic curve given by the equation $y^2 = x^3 + ax + b$.  Next, for any integers $h$ and $k$ define
\begin{equation*}
T(h, k; p) = \sum_{\alpha \shortmod{p}} \sum_{\beta \shortmod{p} } \lambda_{\alpha, \beta} (p) e\left( \frac{\alpha h + \beta k}{p} \right).
\end{equation*}

\begin{mylemma}  Let $p > 2$ be prime and $\bar{k}$ be defined by $k\bar{k} \equiv 1 \pmod{p}$ if $(k, p) = 1$ and $\bar{0} = 0$.  Then we have
\label{lem:2a}

\begin{equation*}
T(h, k;p) = - \varepsilon_p p^{3/2} \left( \frac{k}{p} \right) e\left( \frac{-h^3 \bar{k}^2}{p} \right), 
\end{equation*}
where $\varepsilon_p$ is the sign of the Gauss sum.  
\end{mylemma}

\begin{proof}  
By definition, 
\begin{equation*}
T(h, k;p) = - \sum_{x \shortmod{p}} \sum_{\alpha \shortmod{p}} \sum_{\beta \shortmod{p}} \left( \frac{x^3 + \alpha x + \beta}{p} \right) e\left( \frac{\alpha h + \beta k}{p} \right).
\end{equation*}
The change of variables $\beta \rightarrow \beta - x^3 - \alpha x$ gives
\begin{eqnarray*}
T(h, k;p) & = & \sum_{x \shortmod{p}}e\left( \frac{- x^3k }{p} \right)  \sum_{\alpha \shortmod{p}} e\left( \frac{\alpha(h - xk)}{p} \right) \sum_{\beta \shortmod{p}} \left(\frac{\beta}{p} \right) e\left( \frac{\beta k}{p} \right) \\
& = & - \varepsilon_p p^{3/2} \left( \frac{k}{p} \right) e\left( \frac{-h^3 \bar{k}^2}{p} \right). 
\end{eqnarray*}
\end{proof}

\subsection{An Estimate for an Incomplete Character Sum in Two Variables}
In this section we establish a general estimate for a character sum in two variables which will be applied in our work.  Precisely we have
\begin{mylemma}  Let $F(u, v)$ be a smooth function satisfying
\label{lem:4a}
\begin{equation}
\label{eq:smoothnesscondition}
F^{(\alpha_1, \alpha_2)}(u, v) u^{\alpha_1}v^{\alpha_2} \leq C(\alpha_1, \alpha_2) \left(1 + \frac{|u|}{U} \right)^{-2}\left(1 + \frac{|v|}{V} \right)^{-2} 
\end{equation}
for any $\alpha_1, \alpha_2 \geq 0$, the superscript on $F$ denoting partial differentiation.
Then
\begin{equation*}
\sum_u \sum_v \chi(u) \psi(v) \Lambda(v)e\left( \frac{u^c}{vq} \right)F(u, v)  \ll_\varepsilon 
U^{\onehalf} V^{\onehalf} \left(1 + \frac{U^c}{V|q|} \right)^{1/2} (l_1 l_2 UV)^{\varepsilon} 
\end{equation*}
where $\chi$ and $\psi$ are nonprincipal Dirichlet characters to the moduli $l_1$ and $l_2$, respectively, $\Lambda$ is the von Mangoldt function, $q$ is a nonzero rational number, $c$ is a positive integer, and $\varepsilon$ is any positive number, the implied constant depending only on $\varepsilon$,  $c$, and the numbers $C(\alpha_1, \alpha_2)$.  In case $\chi$ is principal but not $\psi$ the same bound holds but with $U^{1/2}$ lost.  In case $\psi$ is principal but not $\chi$ the same bound holds but with $V^{1/2}$ lost.
\end{mylemma}
Recall that we are assuming the Riemann Hypothesis for all Dirichlet L-functions; this lemma, of course, relies heavily on the GRH.

\begin{proof}  
First assume $\chi$ and $\psi$ are nonprincipal.  To begin, by Mellin inversion, 
\begin{equation*}
\sum_u \sum_v \chi(u) \psi(v) \Lambda(v)e\left( \frac{u^c}{vq} \right) F(u, v)  
\end{equation*}
\begin{equation*}
= \left( \frac{1}{2 \pi i} \right)^2 \int_{(1/2+ \varepsilon)} \int_{(1/2 + \varepsilon )} L(s_1, \chi) \frac{L'}{L}(s_2, \psi) H(s_1, s_2) ds_1 ds_2,
\end{equation*}
where
\begin{equation*}
H(s_1, s_2) = \int_0^\infty \int_0^\infty e\left( \frac{u^c}{vq} \right) F(u, v) u^{s_1} v^{s_2} \frac{du dv}{uv}.
\end{equation*}
Our goal is to obtain the bound 
\begin{equation*}
H(s_1, s_2) \ll \frac{U^{\sigma_1}V^{\sigma_2}}{|(s_1/c + s_2)(s_1/c + s_2 + 1) | \,|s_2|^{1 + \varepsilon}} \left(1 + \frac{U^c}{V|q|} \right)^{1/2} 
\end{equation*}
for $\text{Re}\: s_1 = \sigma_1$, $\text{Re} \: s_2 = \sigma_2$, $\sigma_j = 1/2 + \varepsilon$ or $1 + \varepsilon$, and use the GRH for the bound $L(s_1, \chi) \frac{L'}{L}(s_2, \psi) \ll |s_1 s_2 l_1 l_2|^{\varepsilon/2}$.  Putting these two estimates together will prove the desired result.  The full details of the proof of the bound for $H(s_1, s_2)$ are contained in Appendix \ref{section:appendixa}.

In the cases where one of the characters is principal we integrate over the line $\text{Re} \: s_j = 1 + \varepsilon$ instead of $\text{Re}\: s_j = 1/2 + \varepsilon$ for the appropriate variable.  The bounds on the L-functions are the same.  This accounts for the square-root loss.  The proof is now solely dependent on the details of Appendix \ref{section:appendixa}.
\end{proof}
\subsection{Proof of Corollary \ref{coro:rankbound}}
We follow an easy calculation in \cite{ILS}.  Take $w \geq 0$ and set
\begin{equation*}
p_m(X) = \frac{1}{W_X(\mathcal{F})} \mathop{\sum_{E \in \mathcal{F}}}_{\text{ord}_{s=1} L(s, f) = m} w_X(E).
\end{equation*}
We define the average analytic rank $r$ by
\begin{equation*}
r = \lim_{X \rightarrow \infty} \sum_{m = 1}^{\infty} m p_m(X),
\end{equation*}
if the limit exists.  Our method provides a bound on the limsup of the above quantity.
By taking a test function $\phi$ such that $\phi(x) \geq 0$, $\phi(0) = 1$, and support of $\widehat{\phi}$ contained in $[-\nu, \nu]$ we derive by Theorem \ref{thm:mainresult} and the Plancherel theorem that
\begin{equation*}
\sum_{m = 1}^{\infty} m p_m(X)  \leq g + o(X),
\end{equation*}
where
\begin{equation*}
g = \int_{- \infty}^{\infty} \widehat{\phi}(y) \widehat{\mathcal{W}}(O)(y) dy.
\end{equation*}
This uses GRH for all elliptic curve L-functions so that by positivity we can drop all zeros not at the central point.

By taking the Fourier pair
\begin{equation*}
\phi(t) = \left( \frac{ \sin(\pi \nu t)}{\pi \nu t} \right)^2, \; \; \widehat{\phi}(y) = \frac{1}{\nu} \left( 1 - \frac{|y|}{\nu} \right)
\end{equation*}
we obtain $g = \frac{1}{\nu} + \onehalf$.  Taking $\nu$ less than $7/9$ and letting $X \rightarrow \infty$ shows $r \leq 25/14 + \varepsilon$ for any $\varepsilon > 0$.  An identical calculation works for $\nu < 1 + \delta$ and completes the proof.

\subsection{A Note on Minimality}
\label{section:minimality}
In this section we investigate the family given in Theorem \ref{thm:mainresult} with the restriction that there are no primes $q$ such that both $q^4 | a$ and $q^6 | b$.  This condition ensures that the equation $E:y^2 = x^3 + ax + b$ is   minimal for all $p > 3$.  The result is the same as that stated in Theorem \ref{thm:mainresult}; we omit the precise statement for brevity.

The conductor condition (i.e. the analogue of Lemma \ref{lem:4conductor}) follows directly from Lemma \ref{lem:4conductor}.

The sum over the Fourier coefficients is hardly complicated by the divisibility restrictions.  We easily have
\begin{equation*}
\mathop{\sum_a \sum_b}_{q^4 | a \Rightarrow q^6 \notdiv b} P(E;\phi) w \left( \frac{a}{A}, \frac{b}{B} \right) = \sum_{p > 3} \frac{2 \log{p}}{p \log{X}} \widehat{\phi}\left(\frac{\log{p}}{\log{X}} \right) \mathop{\sum_a \sum_b}_{q^4 | a \Rightarrow q^6 \notdiv b} \lambda_{a, b} (p) w \left( \frac{a}{A}, \frac{b}{B} \right) 
\end{equation*} 
\begin{equation*}
= \sum_{p > 3} \frac{2 \log{p}}{p \log{X}} \widehat{\phi}\left(\frac{\log{p}}{\log{X}} \right) \sum_d \mu(d) \sum_a \sum_b \lambda_{ad^4, bd^6}(p) w \left( \frac{ad^4}{A}, \frac{bd^6}{B} \right) 
\end{equation*}
\begin{equation*}
= \sum_{p > 3} \frac{2 \log{p}}{p \log{X}} \widehat{\phi}\left(\frac{\log{p}}{\log{X}} \right) \mathop{\sum_{d \leq \log{X}}}_{(d, p) = 1} \mu(d) \sum_a \sum_b \lambda_{ad^4, bd^6}(p) w \left( \frac{ad^4}{A}, \frac{bd^6}{B} \right) + o(AB).
\end{equation*}
Completing the sum over $a$ and $b \pmod{p}$ leads to the sum
\begin{equation*}
AB \sum_{p > 3}  \mathop{\sum_{d \leq \log{X}}}_{(d, p) = 1} \mu(d) 
\sum_h \sum_k \sum_\alpha \sum_\beta \lambda_{\alpha d^4, \beta d^6}(p) e\left(\frac{\alpha h + \beta k}{p} \right) 
\Phi(h, k, p) + o(AB),
\end{equation*}
where
\begin{equation*}
\Phi(h, k, p) = \frac{2 \log{p}}{p \log{X}} \widehat{\phi}\left(\frac{\log{p}}{\log{X}} \right) \widehat{w} \left( \frac{hA}{d^4 p}, \frac{kB}{d^6 p} \right).
\end{equation*}
Applying the change of variables $\alpha \rightarrow \alpha d^{-4}$, $\beta \rightarrow \beta d^{-6}$ gives
\begin{equation*}
 AB \sum_{p > 3} \frac{1}{p^2} \mathop{\sum_{d \leq \log{X}}}_{(d, p) = 1} \frac{\mu(d)}{d^{10}} 
\sum_h \sum_k T(h \overline{d}^{4}, k \overline{d}^{6};p) \Phi(h, k, p) + o(AB),
\end{equation*}
where $T$ is given by Lemma \ref{lem:2a}.  
Applying Lemma \ref{lem:2a} we obtain the sum
\begin{equation*}
- \frac{AB}{\log{X}} \sum_{p >3}  \varepsilon_p  \frac{2 \log{p}}{p^{3/2}}\widehat{\phi}\left( \frac{\log p}{\log X} \right) \mathop{\sum_{d \leq \log{X}}}_{(d, p) = 1} \frac{\mu(d)}{d^{10}} \sum_h \sum_k \left( \frac{k}{p} \right) e\left( \frac{-h^3 \bar{k}^2}{p} \right) \widehat{w} \left(\frac{hA}{d^4 p}, \frac{kB}{d^6 p} \right).
\end{equation*}
Now simply remove the restriction $(d, p) = 1$ and move the summation over $d$ to the outside.  We are back to the problem of estimating (\ref{eq:goodform}) except with slightly smaller $A$ and $B$.  We can simply use the same bounds as before; the presence of $d$ contributes at worst bounded powers of $\log{X}$.  Since we showed that each sum of type (\ref{eq:dyadic}) is $\ll X^{-\varepsilon}$ we can use the same proof as for Lemma \ref{lem:5b}.

It is possible to impose minimality restrictions on the Weierstrass equations for other families.  Since no interesting features arise we have limited our discussions to this brief note.

\section{Proof of the Density Theorem for the Main Rank One Family}
\label{section:proof2}
In this section we prove Theorem \ref{thm:rankone}.  The proof will proceed much the same as Theorem \ref{thm:mainresult} but the arguments will be more intricate.

\subsection{The Conductor Condition}
To show (\ref{eq:approximateconductor}) holds we will apply Lemma \ref{lem:lemma3a}.  We take $R(d) = 16(4a^3 + 27b^4)$.  Then $R(d)=|\Delta(d)| \asymp X$.  All we need to show is
\begin{myprop}
\label{prop:conductorrank1}
\begin{equation*}
\sum_{a} \sum_{b} \sum_{p^\alpha || \Delta} \frac{\log{p^{\alpha - 1}}}{\log X} w\left(\frac{a}{A}, \frac{b}{B}\right)  \ll \frac{X^{7/12}}{\log{X}}
\end{equation*}
where $\Delta = -16(4a^3 + 27b^4)$.
\end{myprop}
\begin{proof}  We may assume $\alpha = 2$ because a similar argument to that used in the proof of Lemma \ref{lem:4conductor} with the roles of $a$ and $b$ switched will provide the necessary estimation for $\alpha \geq 3$.  We may further assume $p > 3$ and that $(a, p) = (b, p) = 1$.  
We now make two separate arguments to handle $p$ relatively small and $p$ relatively large.  For the former, we have
\begin{mylemma} For $P = X^{1/3}$
\label{lem:7smallp}
\begin{equation*} 
\sum_{3 < p \leq P} \mathop{\sum_{(a, p) = 1} \sum_{(b, p) = 1}}_{\Delta\equiv 0 \shortmod{p^2}} \frac{\log{p}}{\log{X}} w\left(\frac{a}{A}, \frac{b}{B}\right) \ll \frac{X^{7/12}}{\log{X}},
\end{equation*}
the implied constant depending only on $w$.
\end{mylemma}
For the latter we have
\begin{mylemma} For $P = X^{11/36 + \varepsilon}$
\label{lem:7bigp}
\begin{equation*} 
\sum_{p \geq P} \mathop{\sum_a \sum_b}_{p^2 || \Delta} \frac{\log{p}}{\log{X}} w\left(\frac{a}{A}, \frac{b}{B}\right) \ll \frac{X^{7/12}}{\log{X}},
\end{equation*}
the implied constant depending only on $\varepsilon$ and $w$.
\end{mylemma}

Since $11/36  < 1/3$ these two lemmas will allow us to take $\varepsilon$ small enought to close the gap and complete the proof of Proposition \ref{prop:conductorrank1}.  \end{proof}

\begin{proof}[Proof of Lemma \ref{lem:7smallp}]  By breaking the summation over $a$ into congruence classes $\mymod{p^2}$ we easily obtain the bound
\begin{equation*}
\ll \sum_{p \leq P} \sum_{b \ll B} \left(1 + \frac{A}{p^2} \right) \frac{\log{p}}{\log{X}} \ll \frac{AB}{\log{X}} \left(1 + \frac{P}{A} \right),
\end{equation*}
which is sufficient provided $P \ll A = X^{1/3}$, as claimed. \end{proof}

\begin{proof}[Proof of Lemma \ref{lem:7bigp}]
By majorizing $|w|$ by a smooth non-negative function with slightly larger support we may assume $w \geq 0$.
The conditions on $p$ imply that if we write $|\Delta|$ as $d^2 l$ with $l$ squarefree then $d \geq P$.  Therefore we have
\begin{eqnarray*}
 \sum_{p \geq P} \mathop{\mathop{\sum_a \sum_b}}_{p^2 || \Delta} \frac{\log p}{\log{X}} w \left( \frac{a}{A}, \frac{b}{B} \right) 
 & \leq &  \negthickspace \sum_{d \geq P} \mathop{\mathop{\mathop{\sum_{a} \sum_{b}} }_{|\Delta| = d^2 l}}_{l \text{ squarefree}} w \left( \frac{a}{A}, \frac{b}{B} \right)  \\
& \leq & \mathop{\sum_{l \ll XP^{-2}}}_{l \text{ squarefree}} \mathop{\sum_{a} \sum_{b}  }_{|\Delta| = d^2 l} w \left( \frac{a}{A}, \frac{b}{B} \right).
\end{eqnarray*}
Now define $s = (a, l)$.  Since $s$ is squarefree we must have $s | b$.  Define $l = sl'$, $a = sa'$, and $b = sb'$.  Then the condition on the discriminant is $4 s^2 a'^3 + 27s^3 b'^4 = (d/4)^2 l'$.  Since $l$ is squarefree we have $(s, l') = 1$, which implies $s | 4^{-1}d$.  Thus the condition is reduced to $4 a'^3 + 27sb'^4 = (d/4s)^2 l'$.  The important feature is that $(a', l') = (s, l') = 1$ and $(b', l') = 1 $ or $2$.  Set $l' = (2, l') (3, l') l''$, $a' = (3, a') a''$, $b' = (2, b') b''$, and $s = (2, s) s'$.  Then the discriminant condition becomes 
\begin{equation*}
\frac{4}{(2, l')} \frac{(3, a')^3}{(3, l')} (a'')^3 + \frac{27}{(3, l')} \frac{(2, s) (2, b')^4}{(2, l')} s' (b'')^4 = (d/4s)^2 l''.
\end{equation*}
By taking all the possible combinations of values for $(2, l'), (3, l'), (3, a'), \ldots $, we are left with finitely many equations of the form $c_1 a^3 + c_2 s b^4 = u^2 l$ with the condition $(c_1 c_2 abs, l) = 1$.
Therefore we are reduced to estimating sums of the type
\begin{equation*}
\sum_{s \leq L} \mathop{\sum_{l \leq L/s}{}^{'}}_{l \text{ squarefree}}  \mathop{ \sum_{a}{}^{'} \sum_{b}{}^{'}}_{c_1 a^3 + c_2 s b^4 = u^2 l} w \left( \frac{as}{A}, \frac{bs}{B} \right),
\end{equation*}
where $L \asymp XP^{-2}$ and the primes indicate the summation is restricted by $(c_1c_2abs, l) = 1$.  For squarefree $l$ set
\begin{equation*}
S(s, l) =  \mathop{\sum_{a}{}^{'} \sum_{b}{}^{'}}_{c_1a^3 + c_2sb^4 = u^2 l} w \left( \frac{as}{A}, \frac{bs}{B} \right).
\end{equation*}
Let $\mathcal{Q}$ be a set of primes $q$ of size $q \asymp Q$ (some fixed proportion), with $Q$ at our disposal (it will be chosen to be $X^\eta$ for $\eta$ small).  We detect the condition that an integer $m$ is a square by evaluating the Legendre symbol $\left(\frac{m}{q}\right)$ for $q \in \mathcal{Q}$ (a kind of amplification technique).  We obtain
\begin{equation*}
S(s, l)  =  \mathop{ \mathop{\sum{}^{'} \sum{}^{'}}_{c_1a^3 + c_2 sb^4 = cl}}_{c = \square} w \left( \frac{as}{A}, \frac{bs}{B} \right)
\end{equation*}
\begin{eqnarray*}
 & =  & \frac{1}{|\mathcal{Q}|^2} \mathop{\mathop{\sum {}^{'} \sum{}^{'}}_{c_1 a^3 + c_2 sb^4 = cl}}_{c = \square} \left| \sum_{q \in \mathcal{Q}}  \left( \frac{c}{q} \right) \right|^2 w \left( \frac{as}{A}, \frac{bs}{B} \right) \\
 &  & +{} \; \frac{1}{|\mathcal{Q}|^2}\mathop{ \mathop{\sum{}^{'} \sum{}^{'}}_{c_1 a^3 + c_2 sb^4 = cl}}_{c = \square} \left| \{ q \in \mathcal{Q} : q | c \} \right|^2 w \left( \frac{as}{A}, \frac{bs}{B} \right).
\end{eqnarray*}
Simplifying it becomes
\begin{equation*}
\leq  \frac{1}{|\mathcal{Q}|^2} \mathop{\sum_{a}{}^{'} \sum_{b}{}^{'}}_{c_1a^3 + c_2sb^4 \equiv 0 \shortmod{l}} \left| \sum_{q \in \mathcal{Q}}  \left( \frac{(c_1 a^3 + c_2 sb^4)l}{q} \right) \right|^2 w \left( \frac{as}{A}, \frac{bs}{B} \right) \end{equation*}
\begin{equation*}
+ \frac{1}{|\mathcal{Q}|^2}\mathop{ \mathop{\sum{}^{'} \sum{}^{'}}_{c_1 a^3 + c_2 sb^4 = cl}}_{c = \square} \left| \{ q \in \mathcal{Q} : q | cl \} \right|^2 w \left( \frac{as}{A}, \frac{bs}{B} \right),
\end{equation*}
since in the first sum we have relaxed the condition that $c$ is a square.  For the second sum, notice that since $q \asymp Q$ and $cl \ll X$
\begin{equation*}
\left| \{ q \in \mathcal{Q} : q | cl \} \right| \ll \frac{\log{X}}{\log{Q}}.
\end{equation*}
Therefore the second sum is $\ll (\log{X}/ \log{Q})^2 |\mathcal{Q}|^{-2} S(s, l)$.  If we let $S'(s, l)$ be the first sum above, then
\begin{equation*}
S(s, l) \leq \left(1 - \frac{\log^2{X}}{|\mathcal{Q}|^2 \log^2{Q}} \right)^{-1} S'(s, l).
\end{equation*}
Our choices of $\mathcal{Q}$ and $Q$ will show $S(s, l) \ll S'(s, l)$ so we consider $S'(s, l)$.
We expand the summation over $q$ and obtain
\begin{equation*}
S'(s, l) = \frac{1}{|\mathcal{Q}|^2}  \sum_{q_1 \in \mathcal{Q}} \sum_{q_2 \in \mathcal{Q}}  \mathop{\sum{}^{'} \sum{}^{'}}_{c_1 a^3 + c_2 sb^4 \equiv 0 \shortmod{l}}\left( \frac{(c_1 a^3 + c_2 sb^4)l}{q_1 q_2} \right) w \left( \frac{as}{A}, \frac{bs}{B} \right).
\end{equation*}
Let $S_2 = S_2(s, l, q_1, q_2)$ be the above summation over $a$ and $b$.  Set $r= q_1 q_2$ and apply Poisson summation in $a$ and $b \pmod{lr}$.  We get
\begin{equation*}
S_2 
 =  \frac{AB}{l^2 r^2 s^2} \sum_h \sum_k \mathop{\sum_{u \shortmod{lr}}{}^{'} \sum_{v \shortmod{lr}}{}^{'}}_{c_1u^3 + c_2sv^4 \equiv 0 \shortmod{l}}\left( \frac{(c_1u^3 + c_2 sv^4)l}{r} \right) e\left( \frac{hu + kv}{rl} \right) \widehat{w} \left( \frac{hA}{lrs}, \frac{kB}{lrs} \right). \\
\end{equation*}
Using the Chinese remainder theorem we write $u = u_1r + u_2l$ with $u_1$ given $\mymod{l}$ and $u_2$ given $\mymod{r}$, and similarly for $v$.  Then we get
\begin{equation*}
S_2  =   \frac{AB}{l^2 r^2s^2} \sum_h \sum_k U(l, r, s) \widehat{w} \left( \frac{hA}{lrs}, \frac{kB}{lrs} \right) \mathop{\sum_{u_1}{}^{'}  \sum_{v_1}{}^{'}}_{c_1 u_1^3  + c_2r s v_1^4 \equiv 0 \shortmod{l}} e\left( \frac{hu_1 + kv_1}{l} \right), 
\end{equation*}
where
\begin{equation*}
U(l, r, s) = \sum_{u_2 \shortmod{r}} \sum_{v_2 \shortmod{r}}  \left( \frac{c_1u_2^3 + c_2 lsv_2^4}{r} \right)e\left( \frac{hu_2 + kv_2}{r} \right).
\end{equation*}
Apply the change of variables $u_1 \rightarrow v_1 u_1$ and obtain
\begin{eqnarray*}
S_2 & =   & \frac{AB}{l^2 r^2 s^2} \sum_h \sum_k U(l, r, s) \widehat{w} \left( \frac{hA}{lrs}, \frac{kB}{lrs} \right) \mathop{\sum_{u_1}{}^{'} \sum_{v_1}{}^{'}}_{c_1 u_1^3 + c_2rs v_1 \equiv 0 \shortmod{l}} e\left( \frac{v_1(hu_1 + k)}{l} \right) \\
& = &  \frac{AB}{l^2 r^2s^2} \sum_h \sum_k U(l, r, s) \widehat{w} \left( \frac{hA}{lrs}, \frac{kB}{lrs} \right) \sum_{u_1}{}^{'} e\left( \frac{-c_1\overline{c_2 r s}u_1^3(hu_1 + k)}{l} \right), 
\end{eqnarray*}
since $(c_1c_2rs, l) = 1$.  It's easy to see that the exponential sum over $l$ factors into exponential sums of the form
\begin{equation*}
\sum_{x \shortmod{p}}{}^{'} e\left( \frac{\alpha_p x^3(hx + k)}{p} \right),
\end{equation*}
where $l = \prod_p p$ and $(\alpha_p, p) = 1$.  A corollary of the Riemann Hypothesis for curves (cf. \cite{Schmidt}, Corollary 2F, e.g.) implies that the summation over $x$ is $O(p^{1/2})$ (the implied constant absolute), unless both $h$ and $k$ are zero $\mymod{p}$, in which case the sum is exactly $p - 1$.  Therefore the summation over $u_1$ is $\ll l^{1/2} (h, k, l)^{1/2}\tau(l)$, unless $h = k = 0$, in which case the bound is $l$.
Clearly $U(l, r, s) \ll r^{3/2} + r^2 \delta_{q_1 q_2}$ so we get the bound
\begin{equation*}
S_2 \ll \frac{AB}{l^2 r^2} (r^{\frac{3}{2}} + r^2 \delta_{q_1 q_2} ) \left(l^{\frac{1}{2} + \varepsilon} \mathop{\sum_h \sum_k}_{(h,k) \neq (0, 0)} (h, k, l)^{\frac{1}{2}} \left(1 + \frac{hA}{lrs}\right)^{-2} \left(1 + \frac{kB}{lrs}\right)^{-2} + l \right)
\end{equation*}
\begin{equation*}
\ll \frac{AB}{l^2 r^2 s^2} (r^{\frac{3}{2}} + r^2 \delta_{q_1 q_2} ) \left(l^{\frac{1}{2} + \varepsilon} \frac{l^2 r^2s^2}{AB} + l^{\frac{1}{2} + \varepsilon} \frac{lrs}{B} + l \right) \end{equation*}
\begin{equation*}
= (r^{\frac{3}{2}} + r^2 \delta_{q_1 q_2} ) \left(l^{\frac{1}{2} + \varepsilon} + l^{-\frac{1}{2} + \varepsilon} \frac{A}{rs} +  \frac{AB}{lr^2s^2} \right).
\end{equation*}
Therefore
\begin{equation*}
S'(s, l) \ll 
Q^3 \left(1 + \frac{Q}{|\mathcal{Q}|}\right) \left( l^{1/2 + \varepsilon} + l^{-1/2 + \varepsilon} \frac{A}{sQ^2} + \frac{AB}{l s^2Q^4} \right). 
\end{equation*}
On summation over $l \ll Ls^{-1}$ and $s \leq L$ 
we get 
\begin{equation*}
\sum_{s \leq L} \sum_{l \leq Ls^{-1}} S'(s, l) \ll  
\left(1 + \frac{Q}{|\mathcal{Q}|}\right) \left(Q^3 L^{3/2 + \varepsilon} + AQ L^{1/2 + \varepsilon} + \frac{AB \log{L}}{Q} \right).
\end{equation*}
We take $Q = X^{\varepsilon}$ and $|\mathcal{Q}| \gg \sqrt{Q}$ (which implies $S(s, l) \ll S'(s, l)$).  The necessary bound on this sum is $\ll X^{7/12 - \varepsilon}$, which means the requirement on $L$ is $L \ll X^{7/18 - \varepsilon}.$  Since $L \asymp X/P^2$ the requirement on $P$ is $P \geq X^{11/36 + \varepsilon}$. Now the proof of Lemma \ref{lem:7bigp} is complete. \end{proof}

\subsection{Estimating the Sum of the Fourier Coefficients}
In this section we evaluate $\mathcal{P}(\mathcal{F}; \phi, w_X)$ for the family given in Theorem \ref{thm:rankone}.  This is accomplished with
\begin{mylemma} \label{lem:6b} Set $A=X^{1/3}$ and $B = X^{1/4}$.  Then
\begin{equation*}
\sum_{a} \sum_{b} -P(E;\phi) w\left(\frac{a}{A}, \frac{b}{B}\right) = \phi(0) W_X(\mathcal{F}) + O\left(\frac{X^{7/12}}{\log{X}}\right)
\end{equation*}
provided supp ${\widehat{\phi}} \subset (- \frac{23}{48}, \frac{23}{48})$.
\end{mylemma}
We mimic the proof of Lemma \ref{lem:5b}.  The details are similar so we often condense our arguments.

\begin{proof}  
On Poisson summation,
\begin{equation*}
 \sum_{a} \sum_{b} P(E;\phi) w\left(\frac{a}{A}, \frac{b}{B}\right)
\end{equation*}
\begin{equation*}
 =  AB \sum_{p > 3}  \frac{2 \log{p}}{p^3\log{X}}  \widehat{\phi}\left( \frac{\log p}{\log X} \right) \sum_h \sum_k T'(h, k;p) \widehat{w} \left(\frac{hA}{p}, \frac{kB}{p} \right),
\end{equation*}
where
\begin{equation*}
T'(h, k; p) = \sum_{\alpha \shortmod{p}} \sum_{\beta \shortmod{p} } \lambda_{\alpha, \beta^2} (p) e\left( \frac{\alpha h + \beta k}{p} \right).
\end{equation*}
The complete sum $T'$ can be evaluated explicitly; the calculation is made with Lemma \ref{lem:6}.  
The $\delta(h)\delta(k) p^2$ term in $T'$ gives the extra $\phi(0)W_X(\mathcal{F})$ (by the Prime Number Theorem).  The $p$ term is negligible via trivial estimations.  We are left with estimating the sum
\begin{equation*}
AB \sum_{p > 3} \sum_h \sum_k   \varepsilon_p \frac{2 \log{p}}{p^{3/2}\log{X}}  \widehat{\phi}\left( \frac{\log p}{\log X} \right) \left( \frac{-h}{p} \right) e\left( \frac{\bar{h}^3 k^4\bar{2}^6}{p} \right)  \widehat{w} \left(\frac{hA}{p}, \frac{kB}{p} \right).
\end{equation*}
Estimating this sum trivially obtains our result for support up to $7/18$.


To get larger support we need to show there is cancellation in the sum.  We estimate it in exactly the same way we did in the proof of Lemma \ref{lem:5b}.  First replace $\varepsilon_p$ by a character $\psi_4 \pmod{4}$.  Then break up the sum into dyadic segments using partitions of unity.  It suffices to consider sums of the type
\begin{equation}
\label{eq:rank1dyadic}
\mathop{\mathop{ \mathop{\sum \sum \sum}_{H \leq h < 2H} }_{K \leq k < 2K} }_{P \leq p < 2P}
\psi_4(p) \left( \frac{-h}{p} \right) e\left( \frac{\bar{h}^3 k^4\bar{2}^6}{p} \right) \frac{ \log{p}}{p^{3/2}} \widehat{\phi}\left( \frac{\log p}{\log X} \right) \widehat{w} \left(\frac{hA}{p}, \frac{kB}{p} \right) g(h, k, p),
\end{equation}
where $g$ is a function arising from the partitions of unity.  The contribution from $k=0$ is negligible by trivial estimations.  Let $S(H, K, P)$ be the sum given by (\ref{eq:rank1dyadic}).  It suffices to show $S(H, K, P) \ll X^{-\varepsilon}$.

Using the bound $\widehat{w}(x, y) \ll (1 + |x|)^{-M} (1 + |y|)^{-M}$ we may assume that $H \ll (P/A)^{1 + \varepsilon}$ and $K \ll (P/B)^{1 + \varepsilon}$.

In order to sum over coprime integers we set $d = (k^4, 2^6h^3)$ and define $d_0$ to be the least positive integer such that $d | d_0^4$.  Since $d_0 | k$ we may set $k = d_0 k_0$.  The condition $(k^4, 2^6h^3) = d$ is equivalent to the two conditions $(d_0^4/d, 2^6h^3/d) = 1$ and $(k_0, 2^6h^3/d) = 1$.  Therefore $S(H, K, P)$ is
\begin{equation*}
\mathop{ \mathop{ \sum \sum}_{H \leq h < 2H} }_{P \leq p < 2P}
 \mathop{\sum_{d| 2^6 h^3}}_{\left(\frac{d_0^4}{d}, \frac{2^6 h^3}{d}\right) = 1}  \mathop{\sum_{\frac{K}{d_0} \leq k_0 < 2\frac{K}{d_0}}}_{\left(k_0, 2^6\frac{h^3}{d}\right) = 1}   \psi_4(p) \left( \frac{-h}{p} \right) e\left( \frac{(d_0^4/d) k_0^4 \overline{(2^6h^3/d)} }{p} \right)  U(k_0, d, h, p),
\end{equation*}
where
\begin{equation*}
U(k_0, d, h, p) = g(h, k_0 d_0, p) \widehat{w} \left(\frac{hA}{p}, \frac{d_0k_0B}{p} \right) \frac{ \log{p}}{p^{3/2}} \widehat{\phi}\left( \frac{\log p}{\log X} \right). 
\end{equation*}
We apply the elementary reciprocity formula and the expansion of the exponential into multiplicative characters as in the proof of Lemma \ref{lem:5b} and obtain 
\begin{equation}
\label{eq:expandedrank1}
S(H, K, P) = \sum_{H \leq h < 2H} \sum_{d | 2^6 h^3} \frac{1}{\phi(2^6h^3/d)} \sum_{\chi \shortmod{2^6h^3/d}} \tau(\chi) \bar{\chi}(d_0^4/d) Q(d, h, \chi) 
\end{equation}
where 
\begin{equation*}
Q(d, h, \chi) = \sum_{P \leq p < 2P}  \sum_{K/d_0 \leq k_0 < 2K/d_0} \psi_4(p) \chi(p) \left( \frac{-h}{p} \right) \bar{\chi}^4(k_0) e\left( \frac{d_0^4 k_0^4}{2^6 p h^3} \right)  U(k_0, d, h, p).
\end{equation*}

Since the arguments are now extremely similar to those used in the proof of Lemma \ref{lem:5b} we will be brief.  We apply Lemma \ref{lem:4a} to $Q(d, h, \chi)$ and obtain
\begin{mylemma} \label{lem:Qboundrank1}
If $\chi \psi_4 (-h/\cdot)$ and $\chi^4$ are nonprincipal then
\begin{equation*}
Q(d, h, \chi) \ll P^{-1} K^{1/2} d_0^{-1/2} \left( 1 + \frac{K^4}{P H^3} \right)^{1/2} X^{\varepsilon}.
\end{equation*}
If $\chi^4$ is principal but $\chi \psi_4 (-h/\cdot)$ is not principal we lose a factor $K^{1/2} d_0^{-1/2}$.
\end{mylemma}

Let $S = S_1 + S_2 + S_3$, where $S_1$ corresponds to the terms where both characters are nonprincipal, $S_2$ corresponds to the terms where $\chi^4$ is principal but $\chi \psi_4 (-h/\cdot)$ is not principal, and $S_3$ corresponds to the remaining terms where $\chi \psi_4 (-h/\cdot)$ is principal (and necessarily $\chi^4$ is principal). 

{\bf Case 1.}  We apply Lemma \ref{lem:Qboundrank1} to $S_1$ and obtain the bound
\begin{eqnarray*}
S_1 & \ll & P^{-1} K^{1/2}  \left( 1 + \frac{K^4}{P H^3} \right)^{1/2} X^{\varepsilon} \sum_{H \leq h < 2H} \sum_{d | 2^6 h^3} h^{3/2} (dd_0)^{-1/2}  \\
& \ll & P^{-1} H^{5/2} K^{1/2}  \left( 1 + \frac{K^2}{P^{1/2} H^{3/2}} \right) X^{\varepsilon},
\end{eqnarray*}
which is $\ll X^{-\varepsilon}$ when $P \ll X^{23/48 - \varepsilon}$.

{\bf Case 2.} Applying Lemma \ref{lem:Qboundrank1} to $S_2$ gives
\begin{eqnarray*}
S_2 & \ll & P^{-1}  K \left( 1 + \frac{K^4}{P H^3} \right)^{1/2} X^{\varepsilon} \sum_{H \leq h < 2H} h^{-3/2} \sum_{d | 2^6 h^3}  d^{1/2}d_0^{-1}  \\
& \ll & P^{-1} K \left( 1 + \frac{K^2}{P^{1/2} H^{3/2}} \right) X^{\varepsilon}.
\end{eqnarray*}
Using only $K \ll (P/B)^{1 + \varepsilon}$ and $H \gg 1$ shows this is $\ll X^{-\varepsilon}$ when $P \ll X^{1/2 - \varepsilon}$.

{\bf Case 3.} Now consider the case where $\chi \psi_4 (-h/\cdot)$ (and hence $\chi^4$) are principal.  Since $\psi_4 (-h/\cdot)$ has conductor $h^*$ equal to the square free part of $h$ (up to a factor of $2$ or $4$) we get extra saving from the restriction $d | 2^6 h^3(h^*)^{-1}$.  Thus $\chi$ has conductor of size $h^* \ll h$ so $|\tau(\chi)| \ll h^{1/2}$.  We therefore have the bound
\begin{eqnarray*}
S_3 & \ll & P^{-1/2} H^{-5/2} K X^{\varepsilon} \sum_{H \leq h < 2H}  \sum_{d | 2^6 h^3/h^*} \frac{d}{d_0} \\
& \ll & P^{-1/2} K  X^{\varepsilon} \prod_{p} \left(1 + p^{-3/2} + \ldots \right) 
\end{eqnarray*}
which is $\ll X^{-\varepsilon}$ when $P \ll X^{1/2 - \varepsilon}$ since the infinite product converges.  Having considered all possible cases the proof is complete. \end{proof}
\subsection{A Complete Character Sum}
As in the proof of Theorem \ref{thm:mainresult} we need to evaluate a complete character sum.  We do this now.  Set
\begin{equation*}
T'(h, k; p) = \sum_{\alpha \shortmod{p}} \sum_{\beta \shortmod{p} } \lambda_{\alpha, \beta^2} (p) e\left( \frac{\alpha h + \beta k}{p} \right).
\end{equation*}
We have
\begin{mylemma}  Let $p > 2$ be prime and $\bar{h}$ be defined by $h\bar{h} \equiv 1 \pmod{p}$ if $(h, p) = 1$ and $\bar{0} = 0$.  Then we have
\label{lem:6}
\begin{equation*}
T'(h, k;p) = -p^2 \delta(h)\delta(k) - \varepsilon_p p^{3/2} \left( \frac{-h}{p} \right) e\left( \frac{k^4 \bar{h}^3 \bar{2}^6}{p} \right) + p,
\end{equation*}
where $\delta$ is the Kronecker delta function $\mymod{p}$.
\end{mylemma}

\begin{proof}  
By definition, 
\begin{equation*}
T'(h, k;p) = - \sum_{x \shortmod{p}} \sum_{\alpha \shortmod{p}} \sum_{\beta \shortmod{p}} \left( \frac{x^3 + \alpha x + \beta^2}{p} \right) e\left( \frac{\alpha h + \beta k}{p} \right).
\end{equation*}
\begin{equation*}
= - \sum_{\alpha} \sum_{\beta} \left( \frac{\beta^2}{p} \right) e\left( \frac{\alpha h + \beta k}{p} \right) - \sum_{x \neq 0} \sum_{\alpha} \sum_{\beta} \left( \frac{x^3 + \alpha x + \beta^2}{p} \right) e\left( \frac{\alpha h + \beta k}{p} \right).
\end{equation*}
\begin{equation*}
= T_1' + T_2',
\end{equation*}
say.  We easily have
\begin{equation*}
T_1' = \begin{cases}
		-p(p-1) & \text{if $h\equiv k \equiv 0 \mymod{p}$},
\\
		p & \text{if $h\equiv 0, k \not \equiv 0 \mymod{p}$},
\\
		0 & \text{if $h \not \equiv 0 \mymod{p}$}.
	\end{cases}
\end{equation*}

The sum $T_2'$ is, after the linear change of variables $\alpha \rightarrow \alpha - x^2 - \beta^2 \bar{x}$, given by
\begin{eqnarray*}
 T_2' & =  & - \sum_{x \neq 0}\left(\frac{x}{p} \right)\sum_{\beta}   e\left( \frac{ -h\beta^2 \bar{x} + \beta k - hx^2 }{p} \right) \sum_{\alpha} \left(\frac{\alpha}{p} \right) e\left( \frac{\alpha h}{p} \right) \\
& = & - \varepsilon_p p^{1/2}\left(\frac{h}{p} \right) \sum_{x \neq 0}\left(\frac{x}{p} \right) e\left(\frac{-hx^2}{p}\right) \sum_{\beta}   e\left( \frac{-h\beta^2 \bar{x} + \beta k }{p} \right) \\
& = & - \varepsilon_p p^{1/2}\left(\frac{h}{p} \right) \sum_{x \neq 0}\left(\frac{x}{p} \right)  e\left(\frac{-hx^2}{p}\right) \sum_{\beta}  e\left( \frac{-h x \beta^2 +  xk\beta}{p} \right) \; \text{(from $\beta \rightarrow x\beta$)}.
\end{eqnarray*}
To evaluate the summation over $\beta$ we apply the formula
\begin{equation}
\label{eq:quadraticsum}
\sum_{x \shortmod{p}} e\left( \frac{ax^2 + bx}{p} \right) = \begin{cases}
		\varepsilon_p \sqrt{p} \left(\frac{a}{p}\right) e\left( \frac{-\bar{a}b^2 \bar{4}}{p} \right)   & \text{if $(a, p) = 1$},
\\
		p & \text{if $a\equiv b \equiv 0 \mymod{p}$},
\\
		0 & \text{otherwise}.
	\end{cases}
\end{equation}

We obtain
\begin{equation*}
T_2' = - \varepsilon_p^2 p  \left(\frac{h^2}{p} \right) \left(\frac{-1}{p} \right) \sum_{x \neq 0} e\left( \frac{-hx^2 + k^2\bar{h}\bar{4}x}{p} \right).  
\end{equation*}
Applying (\ref{eq:quadraticsum}) again we obtain
\begin{equation*}
T_2' = - p \left(\frac{h^2}{p} \right) \left( \varepsilon_p \sqrt{p} \left(\frac{-h}{p} \right) e\left( \frac{\bar{h}^3k^4 \bar{2}^6}{p}\right) -1 \right).
\end{equation*}
Gathering terms and simplifying we finish the proof of the lemma. \end{proof}

\section{Conjecturally Enlarging the Support for the Family of all Elliptic Curves}
\label{section:conjecture}
\subsection{A Conjecture on the Size of a Character Sum}
In this section we investigate heuristically the behavior of $\mathcal{D}(\mathcal{F}; \phi, w_X)$ for the family given in Theorem \ref{thm:mainresult} when we extend the support of $\widehat{\phi}$ outside the range $(-1, 1)$.  Recall that this is the splitting point for the symmetry types $O$, $SO(\text{even})$, and $SO(\text{odd})$.  We predict that the symmetry type is $O$.  To provide evidence for this conjecture, we need to argue that
\begin{equation*}
\sum_a \sum_b P(E; \phi) w \left(\frac{a}{A}, \frac{b}{B} \right) = o(AB)
\end{equation*}
for $\widehat{\phi}$ with large support.  From (\ref{eq:goodform}) and (\ref{eq:goodestimate}) the problem is basically reduced to estimating the following character sum
\begin{equation}
\label{eq:S}
\sum_{k \leq K} \sum_{h \leq H} \sum_{p \leq P}  \left( \frac{k}{p} \right) e\left( \frac{h^3}{pk^2} \right) e\left( \frac{h^3 \bar{p}}{k^2} \right),
\end{equation}
with certain relations on the sizes of $H, K$, and $P$.  One of the relations is $H^3/PK^2 \asymp 1$ so for purposes of testing we ignore the $e(h^3/pk^2)$ term. 
\begin{myconj}  
\label{conj:mainconjecture}
There exists $\delta > 0$ and $\varepsilon > 0$ such that if 
$k \leq P$, $k$ is not a square (i.e. the character $(k/\cdot)$ is nonprincipal), and $H = P^{2/3 + \delta}$,  then 
\begin{equation*}
\sum_{h \leq H} \sum_{p \leq P} \left( \frac{k}{p} \right) e\left( \frac{h^3 \bar{p}}{k^2} \right) \ll P^{1 - 3\delta/2 - \varepsilon}.
\end{equation*}
\label{conj:primesum}
\end{myconj}
On summation over $k \leq P^{1/2 + 3\delta/2}$ this conjecture would indicate that Theorem \ref{thm:mainresult} remains true with test functions whose Fourier transforms $\widehat{\phi}$ have support outside of $(-1, 1)$.  The extent to which the support could exceed $(-1, 1)$ would depend on the value of $\delta$.  Precisely, we would obtain support up to $(1 - 3\delta)^{-1}$.  The heuristic used in Section \ref{section:evidence} lends support to the value 
$\delta = 1/48$, which would give support up to 
$16/15$.

\subsection{Evidence for Conjecture \ref{conj:mainconjecture}} 
\label{section:evidence}
To lend support for Conjecture \ref{conj:mainconjecture} we investigate the same sum but with $p$ ranging over positive integers coprime with $k$ instead of primes.  We have
\begin{mytheo} 
\label{thm:integersum}
Let $H = P^{2/3 + \delta}$, $k\asymp K = H^{3/2}P^{-1/2}$.  Then for any $\varepsilon > 0$ we have
\begin{equation*}
\sum_h{}^{'} \sum_m{} \left( \frac{m}{k} \right) e\left( \frac{h^3 \bar{m}}{k^2} \right) w\left( \frac{h}{H}, \frac{m}{P} \right) \ll_\varepsilon P^{5/6 + c\delta + \varepsilon},
\end{equation*}
where the prime indicates the summation is restricted to $(h, k) = 1$ and $c$ is a positive constant (the proof gives $c = 13/2$).
\end{mytheo}
Taking $\delta < \frac{1}{6}(c + \frac{3}{2})^{-1}$ and $\varepsilon$ small enough will show that the sum is $\ll P^{1 - 3\delta/2 - \varepsilon}$ for positive $\delta$ and $\varepsilon$.  The value $c = 13/2$ allows us to take any $\delta < 1/48$.

\begin{proof}  Let $S = S(k)$ be the sum to be estimated.  By Poisson summation in $m \pmod{k^2}$,
\begin{eqnarray*}
S & = & \sum_h{}^{'} \sum_m{} \left( \frac{m}{k} \right) e\left( \frac{h^3 \bar{m}}{k^2} \right) w\left( \frac{h}{H}, \frac{m}{P} \right) \\
& = & \frac{P}{k^2} \sum_h{}^{'}  \sum_{l}  \sum_{x \shortmod{k^2}} \left( \frac{x}{k} \right) e\left( \frac{h^3 \bar{x} + lx}{k^2} \right) w\left( \frac{h}{H}, \widehat{\frac{lP}{k^2}} \right), 
\end{eqnarray*}
where the hat over the second variable indicates we have taken the Fourier transform in that variable only.
Recall $K^2 = H^3/P = P^{1 + 3\delta}$ so $k^2/P \asymp P^{3\delta}$.  Now write $x = y(1 + kz)$ where $y$ and $z$ range over representatives $\mymod{k}$.
Then $\bar{x} = \bar{y} (1 - kz)$ where $\bar{y}$ is the multiplicative inverse of $y \mymod{k^2}$.  We obtain
\begin{eqnarray*}
S & = & \frac{P}{k^2} \sum_h{}^{'}  \sum_{l}  w\left( \frac{h}{H}, \widehat{\frac{lP}{k^2}} \right) \sum_{y \shortmod{k}} \sum_{z \shortmod{k}} \left( \frac{y}{k} \right) e\left( \frac{h^3  \bar{y}(1 - kz) + ly(1 + kz)}{k^2} \right) \\
& = & \frac{P}{k^2} \sum_h{}^{'}  \sum_{l}  w\left( \frac{h}{H}, \widehat{\frac{lP}{k^2}} \right) \sum_{y \shortmod{k}} \left( \frac{y}{k} \right) e\left( \frac{h^3 \bar{y} + l y}{k^2} \right) \sum_{z \shortmod{k}}  e\left( \frac{z(-h^3 \bar{y} + l y)}{k} \right) \\
& = & \frac{P}{k} \mathop{\sum_h{}^{'}  \sum_{l} \sum_{y \shortmod{k}}}_{ly^2 \equiv h^3 \shortmod{k}} \left( \frac{y}{k} \right) e\left( \frac{h^3 \bar{y} + ly}{k^2} \right)  w\left( \frac{h}{H}, \widehat{\frac{lP}{k^2}} \right).
\end{eqnarray*}
Write $h = h_0 + h_1 k$ where $h_0 \asymp H$, $(h_0, k) = 1$,  $h_1 \asymp H/k$, $h_0$ takes values in an interval of length $k$, and $h_1$ takes values in an interval of length $\asymp H/k$.  Now extend the summation over $h_0$ to an interval of length $\asymp H$, so the sum is repeated $\asymp H/k$ times.  We obtain
\begin{equation*}
S \ll \frac{P}{k} \frac{k}{H} \mathop{\sum_{h_0}{}^{'}  \sum_{l} \sum_{y \shortmod{k}}}_{ly^2 \equiv {h_0}^3 \shortmod{k}} \sum_{h_1} \left( \frac{y}{k} \right) e\left( \frac{h_0^3 \bar{y} + 3h_0^2 h_1 k \bar{y}}{k^2} \right) e\left( \frac{ly}{k^2}\right) w\left( \frac{h_0 + h_1k}{H}, \widehat{\frac{lP}{k^2}} \right).
\end{equation*}
Applying the change of variables $y \rightarrow h_0 ^2  \bar{y}$ gives
\begin{equation*}
S \ll \frac{P}{k} \frac{k}{H} \mathop{\sum_{h_0}{}^{'} \sum_{l} \sum_{y \shortmod{k}}}_{y^2 \equiv l h_0 \shortmod{k}} \sum_{h_1} \left( \frac{y}{k} \right) e\left( \frac{ h_0^2 l \bar{y}}{k^2} \right) e\left( \frac{(h_0 +3h_1k) y}{k^2} \right) w\left( \frac{h_0 + h_1k}{H}, \widehat{\frac{lP}{k^2}} \right).
\end{equation*}
There will be virtually no oscillation in $e( (h_0 + 3 h_1 k)y/k^2)$ if $y \ll P^{\varepsilon} k^2/H \ll P^{1/3 + 2\delta + \varepsilon}$.  On the other hand, if $y \gg P^{1/3 + 2\delta + \varepsilon}$ then it is easily shown that the summation over $h_1$ is $\ll_{\varepsilon, M} P^{-M}$.  Therefore we have fixed $0 < y \ll P^{\varepsilon} K^2/H$.  Thus the sum is reduced to
\begin{equation*}
S \ll \frac{P}{H} \sum_{|h_1| \ll H/k} \mathop{\sum_{h_0}  \sum_{l} \sum_{0 < y \ll P^{\varepsilon} K^2/H}}_{y^2 \equiv l h_0 \shortmod{k}} \left( \frac{y}{k} \right) e\left( \frac{h_0^2 l \bar{y}}{k^2} \right) w_1(h_0, h_1, k, l, y),
\end{equation*}
where $w_1$ is the new test function obtained by absorbing the non-oscillatory exponential factor into $w$ (any process of differentiation of $w_1$ with respect to $y$ introduces only factors of size $P^{\varepsilon}$).
Now we apply the elementary reciprocity law (\ref{eq:elementaryreciprocity}), obtaining
\begin{equation*}
S \ll \frac{P}{H} \sum_{|h_1| \ll H/k} \mathop{\sum_{h_0}  \sum_{l} \sum_{0 < y \ll P^{\varepsilon}K^2/H}}_{y^2 \equiv l h_0 \shortmod{k}}  \left( \frac{y}{k} \right) e\left( - \frac{h_0^2 \bar{k}^2 l}{y} \right) e \left( \frac{h_0^2 l}{yk^2} \right) w_1(h_0, h_1, k, l, y).
\end{equation*}
The equality $y^2 = l h_0 + sk$ means $|s| \ll (P^{2/3 + 4\delta}/k)^{1 + \varepsilon} \ll  P^{1/6 + 5\delta/2 + \varepsilon}$.  Then the sum is reduced to
\begin{eqnarray*}
S & \ll & \frac{P}{H} \sum_{|h_1| \ll \frac{H}{k}} \mathop{\sum_{h_0} \sum_{l} \sum_{0 < y \ll P^{\varepsilon} \frac{K^2}{H} }}_{y^2 = l h_0 + sk}  \left( \frac{y}{k} \right) e\left(- \frac{s^2 \bar{l}}{y} \right) e \left( \frac{(y^2 - sk)^2 }{lyk^2} \right) w_1(h_0, h_1, k, l, y) \\
& = & \frac{P}{H} \sum_{|h_1| \ll H/k} \mathop{\sum_{s}  \sum_{l} \sum_{0 < y \ll P^{\varepsilon} K^2/H}}_{y^2 \equiv sk \shortmod{l}} \left( \frac{y}{k} \right) 
e\left(  \frac{s^2\bar{y}}{l} \right) e\left( - \frac{s^2}{ly} + \frac{(y^2 - sk)^2 }{lyk^2} \right) w_2
\end{eqnarray*}
where $w_2 = w_2(s, h_1, k, l, y),$ is the replacement of $w_1$ after the change of variables which eliminates $h_0$ and introduces $s$.
Now we break the summation over $y$ into progressions $\mymod{l}$ and apply the P\'{o}lya-Vinogradov bound of
\begin{equation*}
\mathop{\sum_{y \ll Y}}_{y \equiv \lambda \mymod{l}} \left( \frac{y}{k} \right) \ll (lY)^{1/2} \log{lY}
\end{equation*}
to $S$ (via partial summation).  This gives the final bound of $P^{5/6 + c \delta + \varepsilon}$ where $c = 13/2$. \end{proof}



\section{Curves with Torsion $\mz/2\mz \times \mz/2\mz$}
\label{section:Z2Z2torsion}
In the following sections we study some interesting families of elliptic curves that have torsion points.  Each family has its own interesting features.  One common feature is that the square divisors of the conductor are generally rather easy to control.  The source of this ease is that the discriminants factor into polynomials of smaller degree.  In addition, the conductors are often much smaller than the discriminant because of high multiplicity in one or more of these polynomial factors.  This fact causes these torsion families to have a rather large number of curves with conductor $N \leq X$.

We refer to the paper of Kubert \cite{Kubert} as a reference for these torsion families.  In particular Table 3 contains essentially all the information we use.

In this section we investigate a particularly interesting family of elliptic curves, given in Weierstrass form by
\begin{equation*}
E: y^2 = x(x-a)(x+b).
\end{equation*}
$E$ has discriminant
\begin{equation*}
\Delta = 16a^2 b^2 (a+b)^2.
\end{equation*}
The torsion group is generated by the points $(0, 0)$ and $(-a, 0)$ and
\begin{equation*}
\lambda(p) = - \sum_{x \shortmod{p}} \left( \frac{x(x-a)(x+b)}{p} \right).
\end{equation*}
Helfgott has shown that the root number in this family is equidistributed \cite{Helfgott}.

We have the following
\begin{mytheo}
\label{thm:z2z2torsion}
Let $\mathcal{F}$ be the family of elliptic curves given by the Weierstrass equations $E_{a,b}: y^2 = x(x-a)(x+b)$ with $a$ and $b$ positive integers.  Set $A = B = X^{1/3}$, let $w$ be a smooth compactly supported function on $\mr^{+} \times \mr^+$, and set $w_X(E_{a, b}) = w\left(\frac{a}{A}, \frac{b}{B}\right)$.
Then
\begin{equation*}  
\mathcal{D}(\mathcal{F}; \phi, w_X) \sim [\widehat{\phi}(0) + \thalf \phi(0)]W_X(\mathcal{F}) \; \; \text{as } X \rightarrow \infty,
\end{equation*}
for $\phi$ with $\text{supp } \widehat{\phi} \subset (- \frac{2}{3}, \frac{2}{3})$.
\end{mytheo}
This family has particular interest because we can sum primes up to the size of the family (we are taking $\asymp AB = X^{2/3}$ curves), a natural barrier for any family (since square-root cancellation coming solely from averaging over the family gives us this support; to go further requires additional cancellation in the $\lambda(p)$'s as $p$ varies, at least on average).  To prove Theorem \ref{thm:z2z2torsion} we need two lemmas.  For the conductor condition we have
\begin{mylemma} Let $\mathcal{F}$ be the family given in Theorem \ref{thm:z2z2torsion}.  Then
\begin{equation*}
\sum_a \sum_b \frac{\log{N}}{\log{X}} w \left( \frac{a}{A}, \frac{b}{B} \right) = W_X(\mathcal{F}) + O\left(\frac{AB}{\log{X}}\right).
\end{equation*}
\end{mylemma}
\begin{proof}  We will apply Lemma \ref{lem:lemma3a}.  We take $R(d) = ab(a + b)$.  Then it's clear that $R(d) \asymp X$ and that the irreducible factors of $R(d)$ all divide $\Delta(d)$.
We first consider
\begin{equation*}
\sum_a \sum_b \mathop{\sum_{p || ab(a + b)}}_{p^2 | N, \; p > 3}  w \left( \frac{a}{A}, \frac{b}{B} \right) \log{p}.
\end{equation*}
This sum is empty because $p^2 | N$ implies $p|(a, b)$ (since $x(x-a)(x + b)$ has a triple root $\mymod{p}$).  To estimate the sum
\begin{equation*}
\sum_a \sum_b \sum_{p^\alpha || ab(a + b)} w \left( \frac{a}{A}, \frac{b}{B} \right) \log{p^{\alpha - 1}} 
\end{equation*}
we suppose $p^\gamma || a$ and $p^\delta || b$.  Thus $p^{\alpha - \gamma - \delta} || a + b$.  We first suppose $\gamma \neq \delta$.  By symmetry we may assume $\gamma < \delta$ (which implies $p^\gamma || a + b$ and  $\gamma < \alpha/3$).  Since $a \asymp A = B = X^{1/3}$ we always have $p \ll X^{1/3}$.  Then we get the bound
\begin{equation*}
\mathop{\mathop{\sum_{p^\alpha \ll X}}_{\alpha >0} }_{p \ll X^{1/3}}  \sum_{\gamma \leq \alpha/3} \left(1 + \frac{A}{p^\gamma} \right) \left(1 + \frac{B}{p^{\alpha - 2\gamma}} \right) \log{p^{\alpha - 1}} \ll X^{2/3}.
\end{equation*}
In case $\gamma = \delta$ we apply the change of variables $a \rightarrow a' p^\gamma, b \rightarrow  b' p^\gamma$ and get the bound
\begin{equation*}
\mathop{\mathop{\sum_{p^\alpha \ll X}}_{\alpha >0} }_{p \ll X^{1/3}} \sum_{\gamma \leq \alpha/3} \mathop{\sum_{a'} \sum_{b'}}_{p^{\alpha - 3\gamma} || a' + b'} w \left( \frac{a'p^\gamma}{A}, \frac{b'p^\gamma}{B} \right) \log{p^{\alpha - 1}} .
\end{equation*}
For fixed $a'$, $b'$ is determined $\mymod{p^{\alpha - 3\gamma}}$ so we get the bound
\begin{equation*}
\mathop{\mathop{\sum_{p^\alpha \ll X}}_{\alpha >0} }_{p \ll X^{1/3}}  \sum_{\gamma \leq \alpha/3}  \left(1 + \frac{A}{p^\gamma} \right) \left(1 + \frac{B}{p^{\alpha - 2\gamma}}\right) \log{p^{\alpha - 1}},
\end{equation*}
which is the same bound we have for $\gamma \neq \delta$, so the proof is complete. \end{proof}

The sum over the Fourier coefficients is handled with
\begin{mylemma} 
\label{lem:z2z2torsionsum}
Let $\mathcal{F}$ be the family given in Theorem \ref{thm:z2z2torsion}.  Then
\begin{equation*}
\sum_a \sum_b P(E; \phi) w \left( \frac{a}{A}, \frac{b}{B} \right) \ll \frac{X^{2/3}}{\log{X}}
\end{equation*}
provided supp $\widehat{\phi} \subset (- \frac{2}{3}, \frac{2}{3})$.
\end{mylemma}
\begin{proof}  Set
\begin{equation*}
S(p) = \sum_a \sum_b \lambda(p) w \left( \frac{a}{A}, \frac{b}{B} \right)
\end{equation*}
so that
\begin{equation*}
\sum_a \sum_b P(E; \phi) w \left( \frac{a}{A}, \frac{b}{B} \right) = \sum_{p > 3} S(p) \frac{2 \log{p}}{p \log{X}}  \widehat{\phi}\left( \frac{\log{p}}{\log{X}}\right).
\end{equation*}
Then
\begin{eqnarray*}
S(p) & = & - \sum_{x \shortmod{p}} \sum_a \sum_b \left( \frac{x(x-a)(x+b)}{p} \right) w \left( \frac{a}{A}, \frac{b}{B} \right) \\
& = & - \frac{AB}{p^2} \sum_h \sum_k \mathop{\sum \sum \sum}_{x, \rho, \sigma \shortmod{p}} \left( \frac{x(x-\rho)(x+\sigma)}{p} \right) e \left( \frac{\rho h + \sigma k}{p} \right) \widehat{w} \left( \frac{hA}{p}, \frac{kB}{p} \right) \\
& = & - \frac{AB}{p^2} \sum_h \sum_k \mathop{\sum \sum \sum}_{x, \rho, \sigma \shortmod{p}}    \left( \frac{x \rho \sigma}{p} \right) e \left( \frac{ - \rho h + \sigma k + x(h-k)}{p} \right) \widehat{w} \left( \frac{hA}{p}, \frac{kB}{p} \right) \\
& = & - \frac{AB}{p^{1/2}} \varepsilon_p \sum_h \sum_k \left( \frac{hk(h-k)}{p} \right) \widehat{w} \left( \frac{hA}{p}, \frac{kB}{p} \right). 
\end{eqnarray*}
To get further cancellation we sum over $p$.  To handle the variation of $\varepsilon_p$ we introduce a character $\psi_4 \pmod{4}$ evaluated at $p$.  For those terms with $(hk(h-k)/\cdot) \psi_4$ nontrivial we appeal to the Riemann Hypothesis for Dirichlet L-functions to obtain the bound (summing $p \leq P$)
\begin{equation*}
\frac{AB}{P} P^{\varepsilon} \frac{P^2}{AB} = P^{1 + \varepsilon},
\end{equation*}
which is $\ll X^{2/3 - \varepsilon}$ when $P \ll X^{2/3 - \varepsilon}$, i.e. the restriction on the support of $\widehat{\phi}$ is $2/3$.

When $(hk(h-k)/\cdot) \psi_4$ is trivial we do not obtain any saving in the summation over $p$.  The character is trivial only if $hk(h-k) = \pm \square$ so  the problem amounts to estimating the number of such solutions.  This is carried out by the following
\begin{mylemma}
\label{lem:pythagoreantriples}
Let
\begin{equation*} 
C(Y) = \{(h,k) \in \mz^2: |h| + |k| \leq Y, \; hk(h-k) = \pm \square \}.
\end{equation*}
Then $|C(Y)| \ll_\varepsilon Y^{1 + \varepsilon}$.
\end{mylemma}

Before proving the lemma we apply it to our sum.  By the rapid decay of the Fourier transform we may assume $|h|, |k| \ll (PA^{-1})^{1 + \varepsilon}$.  Using the value $Y = (PA^{-1})^{1 + \varepsilon}$ in the lemma we obtain the bound
\begin{equation*}
\frac{AB}{P^{1/2}} \sum_{p \leq P} |C(Y)| \ll B P^{1/2 + \varepsilon}
\end{equation*}
on the contribution of the terms with trivial character $(hk(h-k)/\cdot) \psi_4$.  The contribution is $\mathop{\ll X^{2/3 - \varepsilon}}$ when $P \ll X^{2/3 - \varepsilon}$, as desired.  This will complete the proof of Lemma \ref{lem:z2z2torsionsum} (and hence Theorem \ref{thm:z2z2torsion}) once we prove Lemma \ref{lem:pythagoreantriples}. \end{proof}

\begin{proof}[Proof of Lemma \ref{lem:pythagoreantriples}]  
Suppose we have a solution $hk(h-k) = \pm \square$ where $|h| + |k| \leq Y$.  By possibly changing the signs of $h$ and $k$ and switching the values of $h$ and $k$ we may assume $h > 0$, $k >0$, and $h > k$ (the cases where $hk(h-k) = 0$ are trivial).  Set $g = (h, k)$ and let $g = d^2l$ where $l$ is squarefree.  Set $h = gh'$ and $k = gk'$.  Then we have $l h'k'(h' - k') = \square$.  Since $(h', k') = 1$ and $l$ is squarefree we must have $l = l_1 l_2 l_3$ where $l_1 | h'$, $l_2 | k'$, $l_3 | (h' - k')$ and hence each of $l_1 h'$, $l_2 k'$, and $l_3 (h' - k')$ must be a square.  Thus $h' = l_1 \cdot \square$, $k' = l_2 \cdot \square$, and $h' - k' = l_3 \cdot \square$.  Thus every solution to $hk(h -k) = \square$ is given by
\begin{equation*}
h = d^2 l_1^2 l_2 l_3 x^2, \; \; \; k= d^2 l_1 l_2^2 l_3 y^2,
\end{equation*}
where
\begin{equation}
\label{eq:quadraticformsolution}
l_1 x^2 = l_2 y^2 + l_3 z^2,
\end{equation}
with the restrictions $(x, y) = (x, z) = (y, z) = 1$, $l=l_1 l_2 l_3$ is squarefree, and each $l_i$ is positive.  But now notice that the solution (\ref{eq:quadraticformsolution}) gives rise to the factorization
\begin{equation*}
l_1 l_2 x^2 = (l_2y + \sqrt{-l_2 l_3} z)(l_2 y - \sqrt{-l_2 l_3} z ).
\end{equation*}
The number of such factorizations is clearly bounded by $d_{\mq(\sqrt{-l_1 l_2})} (l x^2)$, the divisor function in the ring of integers of the field $\mq(\sqrt{-l_2 l_3})$.  It's well-known that $d_{K}(n) \leq c(\varepsilon) (|N_K(n)|)^{\varepsilon}$ where $N_K(n)$ is the norm in the number field $K$ and $c(\varepsilon)$ does not depend on the field $K$.  Using $h = d^2 l_1 l x^2 \geq d^2 l x^2$ we easily get the bound $|C(Y)| \ll Y^{1 + \varepsilon}$. \end{proof}

\section{Curves with Three-Torsion}
In this section we investigate another interesting family of elliptic curves.  Consider the elliptic curve $E = E(a, b)$ given by
\begin{equation*}
E: y^2 + a x y - b y = x^3.
\end{equation*}
$E$ has discriminant
\begin{equation*}
\Delta = - (a^3 + 27 b) b^3.
\end{equation*}
The point $(0, 0)$ is non-singular and has order three.  In fact, all curves over $\mq$ with three-torsion are \mq-isomorphic to one in our family.

For $p \neq 2$ our curve $\mymod{p}$ is equivalent to
\begin{equation*}
(2y + ax - b)^2 \equiv 4x^3 + (ax - b)^2 \mymod{p}.
\end{equation*}
Therefore
\begin{equation*}
\lambda(p) = - \sum_{x \shortmod{p}} \left( \frac{4 x^3 + (a x - b)^2}{p} \right).
\end{equation*}

We have the following density theorem for this family.
\begin{mytheo}
\label{thm:threetorsion}
Let $\mathcal{F}$ be the family of elliptic curves given by the Weierstrass equations $E_{a,b}: y^2 +axy - by = x^3$ with $a$ and $b$ positive integers.  Set $A = X^{1/6}$ and $B = X^{1/2}$.  Let $w$ be a smooth compactly supported function on $\mr^{+} \times \mr^+$ and set $w_X(E_{a, b})= w\left(\frac{a}{A}, \frac{b}{B}\right)$.  Then we have
\begin{equation*}  
\mathcal{D}(\mathcal{F}; \phi, w_X) \sim [\widehat{\phi}(0) + \thalf \phi(0)]W_X(\mathcal{F}) \; \; \text{as } X \rightarrow \infty,
\end{equation*}
for $\phi$ with $\text{supp } \widehat{\phi} \subset (- \frac{1}{2}, \frac{1}{2})$.
\end{mytheo}

Notice that $W_X(\mathcal{F}) = \widehat{w}(0, 0) AB+ o(AB)$, so we are taking $\asymp AB = X^{2/3}$ curves from our family.  Conjecture \ref{conj:fullfamily} would allow us to extend the support to 2/3, a natural barrier for the family.

We first prove that the conductor condition (\ref{eq:approximateconductor}) holds.
\begin{mylemma} Let $\mathcal{F}, A, B$ and $w$ be as in Theorem \ref{thm:threetorsion}.  Then
\label{lem:9b} 
\begin{equation*}
\sum_{E \in \mathcal{F}} \frac{\log N_E}{\log X}w_X(E) = W_X(\mathcal{F}) + O\left(\frac{AB}{\log{X}}\right).
\end{equation*}
\end{mylemma}
\begin{proof}  We will apply Lemma \ref{lem:lemma3a}.  We take $R(d) = b(a^3 + 27b)$.  Then it's clear that $R(d) \asymp X$ and that the irreducible factors of $R(d)$ all divide $\Delta(d)$.
We first consider
\begin{eqnarray*}
&  & \sum_{a} \sum_{b}  \mathop{\sum_{p || b(a^3 + 27b)}}_{p^2 | N_E}  w_X(E) \log{p}.
\end{eqnarray*}
We need the fact that for this family $c_4 = a(a^3 + 24b)$ so that the conditions $p|| b(a^3 + 27b)$ and $p| c_4$ cannot be simultaneously satisfied.  To handle the sum
\begin{equation*} 
\sum_{a} \sum_{b} \mathop{\sum_{p^\alpha || b(a^3 + 27b) }}_{\alpha > 0} w\left(\frac{a}{A}, \frac{b}{B} \right) \log p^{\alpha - 1} 
\end{equation*}
we suppose $p^\gamma || b$ (assume $p \neq 3$ for now) and sum over $\gamma \leq \alpha$.  If $\gamma > \alpha/2$ then $p^{\alpha - \gamma} || a^3$, so these terms give
\begin{equation*} 
\mathop{\sum_{p^\alpha \ll X}}_{\alpha > 0} \mathop{\sum_{\alpha/2 < \gamma \leq \alpha}}_{\gamma \equiv \alpha \shortmod{3}} \sum_{a'} \sum_{b'} w\left(\frac{a' p^{(\alpha - \gamma)/3} }{A}, \frac{b' p^\gamma}{B} \right) \log p^{\alpha - 1} 
\end{equation*}
\begin{equation*} 
\ll \mathop{\sum_{p^\alpha \ll X}}_{\alpha > 0} \mathop{\sum_{\alpha/2 < \gamma \leq \alpha}}_{\gamma \equiv \alpha \shortmod{3}}  \left( 1 + \frac{A}{p^{(\alpha - \gamma)/3}} \right) \left( 1 + \frac{B}{p^{ \gamma}} \right) \log p^{\alpha - 1}
\end{equation*}
\begin{equation*} 
\ll \mathop{\sum_{p^\alpha \ll X}}_{\alpha > 0} \mathop{\sum_{\alpha/2 < \gamma \leq \alpha}}_{\gamma \equiv \alpha \shortmod{3}}   \left(1 + \frac{A}{p^{(\alpha - \gamma)/3}} + \frac{B}{p^{ \gamma}}+ \frac{AB}{p^{\alpha + 2\gamma/3}} \right) \log p^{\alpha - 1}
\end{equation*}
\begin{equation*} 
\ll \mathop{\sum_{p^\alpha \ll X}}_{\alpha > 0}\left(\alpha + A+ \frac{B}{p^{[\alpha/2] + 1}}+ \frac{AB}{p^{4\alpha/3}} \right) \log p^{\alpha - 1} \ll X^{2/3}.
\end{equation*}
The terms with $\gamma \leq \alpha/2$ are treated similarly.  In this case $p^\gamma || a^3$ (so $\gamma$ is divisible by 3) and $p^{\alpha - 2\gamma} || a^3/p^\gamma + 27b/p^\gamma$.  Applying the change of variables $a \rightarrow a' p^{\gamma/3}, b \rightarrow b' p^\gamma$ gives (for $\gamma = \alpha/2$ the divisibility condition $p^0 || a'^3 + 27b'$ should be interpreted to mean $p \notdivtext a'^3 + 27b'$)
\begin{equation*} 
\mathop{\sum_{p^\alpha \ll X}}_{\alpha > 0} \mathop{\sum_{\gamma \leq \alpha/2}}_{\gamma \equiv 0 \shortmod{3}} \mathop{\sum_{a'} \sum_{b'}}_{p^{\alpha - 2 \gamma} || a'^3 + 27 b'}  w\left(\frac{a' p^{\gamma/3} }{A}, \frac{b' p^\gamma}{B} \right) \log p^{\alpha - 1} . 
\end{equation*}
For fixed $a'$, $b'$ is determined $\mymod{p^{\alpha - 2\gamma}}$, so we get the bound
\begin{equation*} 
\mathop{\sum_{p^\alpha \ll X}}_{\alpha > 0} \mathop{\sum_{\gamma \leq \alpha/2}}_{\gamma \equiv 0 \shortmod{3}}  \left(1 + \frac{A}{p^{\gamma/3}} \right) \left(1 + \frac{B}{p^{\alpha - \gamma}} \right) \log p^{\alpha - 1}
\end{equation*}
\begin{equation*} 
\ll \mathop{\sum_{p^\alpha \ll X}}_{\alpha > 0}\left(\alpha + A + \frac{B}{p^{\alpha - [\alpha/2]}} + \frac{AB}{p^{\alpha - 2[\alpha/2]/3}}\right) \log p^{\alpha - 1},
\end{equation*}
which is bounded by $X^{2/3}$.  The case $p = 3$ is similar to the above with only superficial changes.
Now the proof of Lemma \ref{lem:9b} is complete. \end{proof}

For the sum of the Fourier coefficients we have 
\begin{mylemma}  For $\mathcal{F}$ and $w$ as in Theorem \ref{thm:threetorsion} we have
\label{lem:threetorsionlemma}
\begin{equation*}
\sum_{a} \sum_{b} P (f ; \phi) w\left( \frac{a}{A}, \frac{b}{B} \right) \ll O\left(\frac{AB}{ \log{X}}\right)
\end{equation*}
for $\phi$ with $\text{supp } \widehat{\phi} \subset (- \frac{1}{2}, \frac{1}{2})$.
\end{mylemma}
\begin{proof}
Let
\begin{equation*}
S(p) = \sum_{a} \sum_{b} \lambda(p) w\left(\frac{a}{A}, \frac{b}{B} \right) = - {\sum_{x \shortmod{p}}} \sum_{a} \sum_{b} \left( \frac{ 4x^3 + (ax + b)^2}{p} \right) w\left( \frac{a}{A}, \frac{b}{B} \right),
\end{equation*}
so that
\begin{equation*}
\sum_{a} \sum_{b}  P (f ; \phi) w\left( \frac{a}{A}, \frac{b}{B} \right) = -\sum_{p } S(p)  \widehat{\phi}\left(\frac{\log p}{\log R}\right) \frac{2 \log p}{p \log X}.
\end{equation*}

To compute $S$, we employ Poisson summation $\mymod{p}$ in $b$, obtaining
\begin{eqnarray*}
S(p) & = & - \frac{B}{p}{\sum_{x \shortmod{p}}}  \sum_{k} \sum_{a} \sum_{\beta \shortmod{p}} \left( \frac{ 4x^3 + (a x + \beta)^2}{p} \right) e \left( \frac{ \beta k}{p} \right) w\left( \frac{a}{A}, \widehat{\frac{kB}{p}} \right),
\end{eqnarray*}
where (in abuse of notation) we have placed a hat over the 2nd variable to indicate that we have applied Poisson summation in that variable only.  Breaking the sum up depending on whether $\beta \equiv 0$ or not and applying the change of variables $x \rightarrow \beta x$ when $\beta \not \equiv 0$, we obtain
\begin{eqnarray*}
S(p) & = &  - \frac{B}{p}{\sum_{x \shortmod{p}}}  \sum_{k} \sum_{a} \sum_{\beta \not \equiv 0} \left( \frac{ 4\beta x^3 + (ax + 1)^2}{p} \right) e \left( \frac{ \beta k }{p} \right) w\left( \frac{a}{A}, \widehat{\frac{kB}{p}} \right) \\
&  & - \frac{B}{p}{\sum_{x \shortmod{p}}}  \sum_{k} \sum_{a} \left( \frac{ 4x^3 + (ax)^2}{p} \right) w\left( \frac{a}{A}, \widehat{\frac{kB}{p}} \right). 
\end{eqnarray*}
Separating the terms where $x = 0$ gives
\begin{eqnarray*}
S(p) & = & - \frac{B}{p}  \sum_k \sum_{a} \sum_{\beta \not \equiv 0} e \left( \frac{ \beta k }{p} \right) w\left( \frac{a}{A}, \widehat{\frac{kB}{p}} \right) \\
& & - \frac{B}{p}{\sum_{x \not \equiv 0}}   \sum_k  \sum_{a} \left( \frac{ 4x + a^2}{p} \right) w\left( \frac{a}{A}, \widehat{\frac{kB}{p}} \right)  \\ 
& & - \frac{B}{p}\sum_{x \not \equiv 0 }   \sum_{a} \sum_{\beta \not \equiv 0} \left( \frac{ 4\beta x^3 + (ax + 1)^2}{p} \right) e \left( \frac{ \beta k }{p} \right) w\left( \frac{a}{A}, \widehat{\frac{kB}{p}} \right).
\end{eqnarray*}
Upon simplification we obtain
\begin{eqnarray*}
S(p) & =  & - \frac{B}{p} \sum_{k} \sum_{a}  (p \delta(k) - 1) w\left( \frac{a}{A}, \widehat{\frac{kB}{p}} \right) + \frac{B}{p}\sum_{k} \sum_{a} \left( \frac{a^2}{p} \right) w\left( \frac{a}{A}, \widehat{\frac{kB}{p}} \right)  \\ 
& & - \frac{B}{p}\sum_{x \not \equiv 0 }  \sum_{k} \sum_{a} \sum_{\beta} \left( \frac{\beta}{p} \right) e \left( \frac{ \bar{4}\bar{x}^3k(\beta - (ax+1)^2) }{p} \right) 
w\left( \frac{a}{A}, \widehat{\frac{kB}{p}} \right) \\
& & + \frac{B}{p}\sum_{x \not \equiv 0 }  \sum_{k} \sum_{a}  \left( \frac{(ax+1)^2}{p} \right) w\left( \frac{a}{A}, \widehat{\frac{kB}{p}} \right).
\end{eqnarray*}
Further simplification gives
\begin{eqnarray*}
S(p) & =  & - B \sum_{k \equiv 0} \sum_{a} w\left( \frac{a}{A}, \widehat{\frac{kB}{p}} \right) + \frac{B}{p}\sum_{k} \sum_{a\not \equiv 0} w\left( \frac{a}{A}, \widehat{\frac{kB}{p}} \right) \\
& & +  \frac{B}{p} \sum_x \sum_k \sum_a \left( \frac{ax+1}{p} \right)^2 w\left( \frac{a}{A}, \widehat{\frac{kB}{p}} \right)  \\ 
& & - \frac{B\varepsilon_p}{\sqrt{p}}\sum_{x }  \sum_{k} \sum_{a} \left( \frac{xk}{p} \right) e \left( \frac{ -\bar{4}\bar{x}^3k (ax+1)^2 }{p} \right)  w\left( \frac{a}{A}, \widehat{\frac{kB}{p}} \right). 
\end{eqnarray*}
Evaluating the 3rd sum and grouping terms gives
\begin{eqnarray*}
S(p) & =  &  B \sum_{k \not \equiv 0} \sum_{a} w\left( \frac{a}{A}, \widehat{\frac{kB}{p}} \right)  \\
& & - \frac{B\varepsilon_p}{\sqrt{p}}\sum_{x }  \sum_{k} \sum_{a} \left( \frac{xk}{p} \right) e \left( \frac{ -\bar{4}\bar{x}^3k (ax+1)^2 }{p} \right)  w\left( \frac{a}{A}, \widehat{\frac{kB}{p}} \right)  \\
& = & - \frac{B\varepsilon_p}{\sqrt{p}} \sum_{x }  \sum_{k} \sum_{a} \left( \frac{xk}{p} \right) e \left( \frac{ -\bar{4}xk (x+a)^2 }{p} \right)  w\left( \frac{a}{A}, \widehat{\frac{kB}{p}} \right) + O(Ap).
\end{eqnarray*}
Bounding the sum over $x \pmod{p}$ using the Riemann Hypothesis for curves (cf. \cite{Schmidt}, Theorem 2G), we have proved
\begin{equation*}
S(p) \ll A p.
\end{equation*}

After summation over $p$ this translates to requiring $\text{supp}\; \widehat{\phi} \subset(-\onehalf , \onehalf)$ and completes the proof (of both Lemma \ref{lem:threetorsionlemma} and Theorem \ref{thm:threetorsion}). \end{proof}

\medskip
To investigate Theorem \ref{thm:threetorsion} with larger support of $\widehat{\phi}$ we use a different form of $S(p)$.  To do so we execute Poisson summation in $a$ also, obtaining
\begin{equation*}
S(p)  =  B \sum_{k \not \equiv 0} \sum_a \left( 1 - \frac{\varepsilon_p}{\sqrt{p}} \sum_{x } \left( \frac{xk}{p} \right) e \left( \frac{ -\bar{4}xk (x+a)^2 }{p} \right) \right) w\left( \frac{a}{A}, \widehat{\frac{kB}{p}} \right) 
\end{equation*}
\begin{eqnarray*} 
& = & \frac{AB}{p} \sum_h \sum_{k \not \equiv 0} \sum_{\alpha \shortmod{p}} e\left( \frac{\alpha h}{p} \right) \widehat{w}\left( \frac{hA}{p}, \frac{kB}{p}\right) \\
& - & \frac{AB}{p} \sum_h \sum_{k \not \equiv 0} \sum_{\alpha \shortmod{p}} \frac{\varepsilon_p}{\sqrt{p}} \sum_{x } \left( \frac{xk}{p} \right) e \left( \frac{ -\bar{4}xk (x+\alpha)^2 + \alpha h }{p} \right)  \widehat{w}\left( \frac{hA}{p}, \frac{kB}{p}\right). 
\end{eqnarray*}
Executing summation over $\alpha$ in the first sum gives
\begin{eqnarray*}
&  & \frac{AB}{p} \sum_{h \equiv 0 \shortmod{p}} \sum_{k \not \equiv 0} p \; \widehat{w}\left(\frac{hA}{p}, \frac{kB}{p}\right). \\
\end{eqnarray*}
Changing variables $\alpha \rightarrow \alpha - x$ and executing summation over $x$ in the second sum gives
\begin{eqnarray*}
&  & - \frac{AB}{p} \sum_h \sum_{k \not \equiv 0} {\varepsilon_p^2} \sum_{\alpha } \left( \frac{-k(\bar{4}k \alpha^2 + h)}{p} \right) e \left( \frac{\alpha h}{p} \right) \widehat{w}\left( \frac{hA}{p}, \frac{kB}{p}\right). 
\end{eqnarray*}
Apply the change of variables $\alpha \rightarrow 2\alpha/k$ (and use $\varepsilon_p^2 = (-1/p)$) and get (for the second sum)
\begin{equation*}
 -\frac{AB}{p} \sum_h \sum_{k \not \equiv 0}  \sum_{\alpha } \left( \frac{\alpha^2 + hk}{p} \right) e \left( \frac{ 2 \alpha \bar{k} h }{p} \right)  \widehat{w}\left( \frac{hA}{p}, \frac{kB}{p}\right).
\end{equation*}
The sum over $\alpha$ is a Kloosterman sum.  To see this, use the identity
\begin{equation*}
\left(\frac{y}{p}\right) = \frac{1}{\varepsilon_p \sqrt{p}} \sum_{\gamma \shortmod{p}} \left(\frac{\gamma}{p}\right) e\left(\frac{\gamma y}{p} \right)
\end{equation*}
so that for $h \not \equiv 0$ we have
\begin{eqnarray*}
& & \sum_{\alpha \shortmod{p}} \left(\frac{\alpha^2 + hk}{p} \right) e \left( \frac{2 h \bar{k} \alpha}{p} \right) \\
& = & \frac{1}{\varepsilon_p \sqrt{p}} \sum_{\gamma \shortmod{p}} \left(\frac{\gamma}{p}\right)  e\left(\frac{\gamma hk}{p} \right) \sum_{\alpha \shortmod{p}} e \left( \frac{ \gamma \alpha^2 + 2 h \bar{k} \alpha}{p} \right) \\
& = & \sum_{\gamma \not \equiv 0} e \left(\frac{hk \gamma - h^2 \bar{k}^2 \bar{\gamma}}{p} \right) 
 =  S(-h^3\bar{k}, 1;p).
\end{eqnarray*}
Thus
\begin{equation*}
S(p) = -\frac{AB}{p} \sum_{h \not \equiv 0} \sum_{k \not \equiv 0}  S(-h^3\bar{k}, 1;p)  \widehat{w}\left( \frac{hA}{p}, \frac{kB}{p}\right) + \frac{AB}{p} \sum_{h \equiv 0} \sum_{k \not \equiv 0} \widehat{w}\left( \frac{hA}{p}, \frac{kB}{p}\right).
\end{equation*}
Now that we have $S(p)$ as a sum of Kloosterman sums (plus a small remainder) we make the following
\begin{myconj} \label{conj:fullfamily} For any $X \geq 1$ and $\varepsilon > 0$ there exists $\delta > 0$ such that we have
\begin{equation*}
\sum_{h \leq X^{1/2}} \sum_{k \leq X^{1/6}} \sum_{p \leq X^{2/3 - \varepsilon}} \frac{S(-h^3\bar{k}, 1;p)}{p} \log{p} \ll X^{2/3 - \delta}.
\end{equation*}
\end{myconj}
This conjecture amounts to saying that on average there is square root cancellation in the summation over $p$, in the ranges of $h$ and $k$ that are of interest to us.  Since it is not relevant for our purposes we do not attempt to formulate a precise conjecture for more general ranges of $h$ and $k$.

If we apply this conjecture to the averaging of $P(E; \phi)$ we get that
\begin{equation*}
\sum_{a} \sum_{b} P(E; \phi) w\left( \frac{a}{A}, \frac{b}{B} \right) \ll \frac{AB}{ \log{X}}
\end{equation*}
for $\text{supp } \widehat{\phi} \subset (-\frac{2}{3}, \frac{2}{3})$.

\section{Curves with Two-Torsion}
In this section we investigate the family given by
\begin{equation*}
E: y^2 = x(x^2 + ax - b).
\end{equation*}
$E$ has discriminant
\begin{equation*}
\Delta = 16b^2( a^2 + 4b ).
\end{equation*}
We have the following
\begin{mytheo}
\label{thm:z2torsion}
Let $\mathcal{F}$ be the family of elliptic curves given by the Weierstrass equations $E_{a,b}: y^2 = x(x^2 + ax - b)$ with $a$ and $b$ positive integers.  Let $w$ be a smooth compactly supported function on $\mr^+ \times \mr^+$.  Set $A = X^{1/4}$, $B = X^{1/2}$, and $w_X(E_{a, b}) = w\left(\frac{a}{A}, \frac{b}{B}\right)$.  Then we have
\begin{equation*}  
\mathcal{D}(\mathcal{F}; \phi, w_X) \sim [\widehat{\phi}(0) + \thalf \phi(0)]W_X(\mathcal{F}) \; \; \text{as } X \rightarrow \infty,
\end{equation*}
for $\phi$ with $\text{supp } \widehat{\phi} \subset (- \frac{1}{2}, \frac{1}{2})$.
\end{mytheo}

To handle the conductor condition we have
\begin{mylemma} Let $\mathcal{F}$ be the family given in Theorem \ref{thm:z2torsion}.  Then
\begin{equation*}
\sum_a \sum_b \frac{\log{N}}{\log{X}} w \left( \frac{a}{A}, \frac{b}{B} \right) = W_X(\mathcal{F}) + O\left(\frac{AB}{\log{X}}\right).
\end{equation*}
\end{mylemma}
\begin{proof}  We will apply Lemma \ref{lem:lemma3a}.  We take  $R(d) = 16b(a^2 + 4b)$.  Then it's clear that $R(d) \asymp X$ and that the irreducible factors of $R(d)$ all divide $\Delta(d)$.
We first consider
\begin{equation*}
\sum_a \sum_b \mathop{\sum_{p || 16b(a^2 + 4b)}}_{p^2 | N, \; p > 3}  w \left( \frac{a}{A}, \frac{b}{B} \right) \log{p}.
\end{equation*}
This sum is empty because $p^2 | N$ implies $p|(a, b)$ (since $x(x^2 +ax - b)$ has a triple root $\mymod{p}$).  It remains to bound the sum
\begin{equation*}
\sum_a \sum_b \mathop{\sum_{p^\alpha || 16b(a^2 + 4b)}}_{\alpha > 0} w \left( \frac{a}{A}, \frac{b}{B} \right) \log{p^{\alpha - 1}}.
\end{equation*}
For $p \neq 2$ the details are so similar to those in the proof of Lemma \ref{lem:9b} that we omit them.  The case $p = 2$ is also simple.  When $a$ is odd the desired bound is immediate.  When $a$ is even we simply change variables $a \rightarrow 2a$ and reduce the divisibility condition to 
$2^{\alpha'} || b (a^2 + b)$, which is virtually identical to the case $p \neq 2$, so again we omit the details.
\end{proof}

For the sum of the Fourier coefficients we have
\begin{mylemma} Let $\mathcal{F}$ be the family given in Theorem \ref{thm:z2torsion}.  Then
\begin{equation*}
\sum_a \sum_b P(E; \phi) w \left( \frac{a}{A}, \frac{b}{B} \right) \ll \frac{X^{3/4}}{\log{X}} 
\end{equation*}
for supp $\widehat{\phi} \subset (- \frac{1}{2}, \frac{1}{2})$.
\end{mylemma}
\begin{proof}  Set
\begin{equation*}
S(p) = \sum_a \sum_b \lambda(p) w \left( \frac{a}{A}, \frac{b}{B} \right).
\end{equation*}
Then by definition,
\begin{eqnarray*}
S(p) & = & - \sum_{x \shortmod{p}} \sum_a \sum_b \left( \frac{x(x^2 + ax -b)}{p} \right) w \left( \frac{a}{A}, \frac{b}{B} \right). 
\end{eqnarray*}
Applying Poisson summation $\mymod{p}$ gives
\begin{eqnarray*}
&  & - \frac{AB}{p^2} \sum_h \sum_k \mathop{\sum \sum \sum}_{x, \alpha, \beta \shortmod{p}} \left( \frac{x(x^2 + \alpha x - \beta )}{p} \right) e \left( \frac{\alpha h + \beta k}{p} \right) \widehat{w} \left( \frac{hA}{p}, \frac{kB}{p} \right). 
\end{eqnarray*}
After the change of variables $\beta \rightarrow -\beta + x^2 + \alpha x$ we have
\begin{equation*}
S(p)  =   - \frac{AB}{p^2} \sum_h \sum_k \mathop{\sum \sum \sum}_{x, \alpha, \beta \shortmod{p}} \left( \frac{x}{p} \right) \left( \frac{\beta}{p} \right) e \left( \frac{\alpha h - \beta k + \alpha x k + x^2 k}{p} \right) \widehat{w} \left( \frac{hA}{p}, \frac{kB}{p} \right). 
\end{equation*}
Evaluating the summation in $\beta$ (it is a Gauss sum) gives
\begin{equation*}
S(p)  =   - \frac{AB}{p^{3/2}} \varepsilon_p \sum_h \sum_k \mathop{\sum \sum}_{x, \alpha \shortmod{p}}\left( \frac{x}{p} \right) \left( \frac{-k}{p} \right) e \left( \frac{\alpha h + \alpha x k + x^2 k}{p} \right) \widehat{w} \left( \frac{hA}{p}, \frac{kB}{p} \right). 
\end{equation*}
We evalute the summation over $\alpha$ and obtain
\begin{eqnarray*}
S(p) & =  & - \frac{AB}{p^{1/2}} \varepsilon_p \sum_h \sum_k \sum_{x \shortmod{p}} \left( \frac{x}{p} \right) \left( \frac{-k}{p} \right) \delta(h + xk) e \left( \frac{ x^2 k}{p} \right) \widehat{w} \left( \frac{hA}{p}, \frac{kB}{p} \right) \\
& = & - \frac{AB}{p^{1/2}} \varepsilon_p \sum_h \sum_k{}^{'} \left( \frac{h}{p} \right) e \left( \frac{ h^2 \bar{k}}{p} \right) \widehat{w} \left( \frac{hA}{p}, \frac{kB}{p} \right),
\end{eqnarray*}
the prime indicating $(k, p) = 1$.
Trivially estimating this sum gives
\begin{equation*}
S(p) \ll p^{3/2}. 
\end{equation*}
Using similar methods to the estimation of (\ref{eq:goodform}) (i.e. the reciprocity law which changes the modulus in the exponential to $k$, the separation of variables using the expansion of the exponential into multiplicative characters, and the Riemann Hypothesis to bound a character sum over primes) we should be able to obtain larger support.  The limit of the method is likely to be $2/3$.  We have not executed the details because the harder (and more interesting) case was already performed in the proof of Theorem \ref{thm:mainresult}.

Using the trivial estimation we obtain
\begin{equation*}
\mathcal{P}(\mathcal{F};\phi, w_X) \ll \sum_{p \leq P} |S(p)| \frac{\log{p}}{p\log{X}} \ll \frac{P^{3/2}}{\log{X}}.
\end{equation*}
which is required to be $\ll X^{3/4}(\log{X})^{-1}$, i.e. we require supp $\widehat{\phi} \subset (- \onehalf, \onehalf)$. \end{proof} 

\section{Curves with Four-Torsion}
In this section we consider the family of curves $E=E(b)$ given by
\begin{equation*}
E: y^2 + xy - by = x^3 - bx^2.
\end{equation*}
This curve has discriminant $\Delta = \Delta(b)$ given by
\begin{equation*}
\Delta = b^4(1 + 16b).
\end{equation*}
The parameter $c_4$ is
\begin{equation*}
c_4 = 16b^2 + 16b + 1.
\end{equation*}
The torsion group is cyclic of order 4 and is generated by the point $(0, 0)$.  We wish to write $\lambda(p)$ as a character sum so we complete the square:  
\begin{equation*}
4\left(y + \frac{x-b}{2}\right)^2 = 4x^3 - 4bx^2 + (x-b)^2 = (x-b) (4x^2 + x - b) 
\end{equation*}
so
\begin{equation*}
\lambda(p) = - \sum_{x \shortmod{p}} \left( \frac{(x-b) (4x^2 + x - b)}{p} \right)
\end{equation*}

We have the following density theorem.
\begin{mytheo}
\label{thm:z4torsion}
Let $\mathcal{F}$ be the family of elliptic curves given by the Weierstrass equations $E_{b}: y^2 + xy - by = x^3 - bx^2$ with $b$ a positive integer.  Let $w$ be a smooth compactly supported function on $\mr^{+}$ and set $B = X^{1/2}$.  Define $w_X(E_b) = w\left(\frac{b}{B}\right)$.  Then we have
\begin{equation*}  
\mathcal{D}(\mathcal{F}; \phi, w_X) \sim [\widehat{\phi}(0) + \thalf \phi(0)]W_X(\mathcal{F}) \; \; \text{as } X \rightarrow \infty,
\end{equation*}
for $\phi$ with $\text{supp } \widehat{\phi} \subset (- \frac{1}{2}, \frac{1}{2})$.
\end{mytheo}
Note that we can sum primes as large as the size of the family here.  

For the conductor condition we have
\begin{mylemma} Let $\mathcal{F}$ be the family given in Theorem \ref{thm:z4torsion}.  Then
\begin{equation*}
\sum_b \frac{\log{N}}{\log{X}} w \left( \frac{b}{B} \right) = W_X(\mathcal{F}) + O\left(\frac{B}{\log{X}}\right).
\end{equation*}
\end{mylemma}
\begin{proof}  
As usual, we will apply Lemma \ref{lem:lemma3a}.  We take $R(d) = b(1 + 16b)$.  Then it's clear that $R(d) \asymp X$ and that the irreducible factors of $R(d)$ all divide $\Delta(d)$.
We first consider
\begin{equation*}
\sum_b \mathop{\sum_{p || b(1 + 16b)}}_{p^2 | N, \; p > 3} w \left( \frac{b}{B} \right) \log{p}.
\end{equation*}
This sum is empty because $(\Delta, c_4) = 1$.  The sum
\begin{equation*}
\sum_b \mathop{\sum_{p^\alpha || b(1 + 16b)}}_{\alpha \geq 2} w \left( \frac{b}{B} \right) \log{p^{\alpha - 1}}
\end{equation*}
is trivially bounded by $X^{1/2}$ since $(b, 1+ 16b) = 1$.
\end{proof}

\begin{mylemma} 
Let $\mathcal{F}$ be the family given in Theorem \ref{thm:z4torsion}.  Then
\begin{equation*}
\sum_b P(E; \phi) w \left( \frac{b}{B} \right) \ll \frac{X^{1/2}}{\log{X}}
\end{equation*}
provided supp $\widehat{\phi} \subset (- \frac{1}{2}, \frac{1}{2})$.
\end{mylemma}
\begin{proof}  
Set
\begin{equation*}
S(p) = \sum_b  \lambda(p) w \left( \frac{b}{B} \right).
\end{equation*}
Then
\begin{eqnarray*}
S(p) & = & - \sum_{x \shortmod{p}} \sum_b \left( \frac{(x-b)(4x^2 + x - b)}{p} \right) w \left( \frac{b}{B} \right) \\
& = & - \frac{B}{p} \sum_h \sum_{x \shortmod{p}} \sum_{\beta \shortmod{p}} \left( \frac{(x- \beta)(4x^2 + x - \beta)}{p} \right) e \left( \frac{\beta h}{p} \right) \widehat{w} \left( \frac{hB}{p} \right),
\end{eqnarray*}
by Poisson summation.  Now separate the terms where $\beta \not \equiv 0$ and apply $x \rightarrow \beta x$ for these terms, getting
\begin{eqnarray*}
S(p) & = &   - \frac{B}{p} \sum_h \sum_{x \shortmod{p}}\left( \frac{ x- 1}{p} \right)  \sum_{\beta \not \equiv 0} \left( \frac{4\beta x^2 + x - 1}{p} \right) e \left( \frac{\beta h}{p} \right) \widehat{w} \left( \frac{hB}{p} \right) \\
& & + \; \frac{B}{p} \sum_h \widehat{w} \left( \frac{hB}{p} \right). 
\end{eqnarray*}
Extending the summation back to $\beta \equiv 0$ and simplifying gives
\begin{eqnarray*}
S(p)  & = &   - \frac{B}{p} \sum_h \sum_{x \shortmod{p}}\left( \frac{ x- 1}{p} \right)  \sum_{\beta } \left( \frac{4\beta x^2 + x - 1}{p} \right) e \left( \frac{\beta h}{p} \right) \widehat{w} \left( \frac{hB}{p} \right) \\
& & + \;  B \sum_h \widehat{w} \left( \frac{hB}{p} \right). 
\end{eqnarray*}
Now separate the terms where $x \not \equiv 0$, getting
\begin{eqnarray*}
S(p) & =  & - \frac{B}{p} \sum_h \sum_{x \not \equiv 0} \left( \frac{ x- 1}{p} \right)  \sum_{\beta } \left( \frac{4\beta x^2 + x - 1}{p} \right) e \left( \frac{\beta h}{p} \right) \widehat{w} \left( \frac{hB}{p} \right) \\
& & - \; B \sum_h \widehat{w}(hB) + B \sum_h \widehat{w} \left( \frac{hB}{p} \right). 
\end{eqnarray*}
Let $S_1$ be the first sum above.  After the change of variables $\beta \rightarrow \bar{4}\bar{x}^2 (\beta - (x - 1))$ we have
\begin{eqnarray*} 
S_1 = - \frac{B}{p} \sum_h \sum_{x \not \equiv 0} \left( \frac{ x- 1}{p} \right)  \sum_{\beta } \left( \frac{\beta}{p} \right) e \left( \frac{\bar{4}h \bar{x}^2 (\beta - x + 1)}{p} \right) \widehat{w} \left( \frac{hB}{p} \right).  
\end{eqnarray*}
After evaluating the Gauss sum we get 
\begin{eqnarray*} 
S_1 & =  & - \frac{B}{p^{1/2}}\varepsilon_p \sum_h \sum_{x \not \equiv 0} \left( \frac{ (x - 1)(\bar{4} h \bar{x}^2)}{p} \right) e \left( \frac{\bar{4}h \bar{x}^2 ( - x + 1)}{p} \right) \widehat{w} \left( \frac{hB}{p} \right)  \\ 
& = & - \frac{B}{p^{1/2}}\varepsilon_p \sum_h \sum_{x} \left( \frac{ (x - x^2)h}{p} \right) e \left( \frac{\bar{4}h ( x^2 - x)}{p} \right) \widehat{w} \left( \frac{hB}{p} \right) \; \; (\text{from } x \rightarrow \bar{x}). \\ 
\end{eqnarray*}
Applying RH for curves (cf. \cite{Schmidt}, Theorem 2G, p. 45) to the summation over $x$ gives $S(p) \ll p$.  In order to apply Weil's bound the requirement is that the polynomials $Y^2 - X^2 + X$ and $Y^p - X^2 + X$ are absolutely irreducible.  We use the condition that $Y^d - f(X)$ is absolutely irreducible if and only if  $(d, d_1, d_2, \ldots, d_s) = 1$, where $f(X) = a (X-x_1)^{d_1} \ldots (X - x_s)^{d_s}$ is the factorization of $f$ in the algebraic closure of $\mf_p$ (assume that the $x_i$'s are distinct) (\cite{Schmidt}, Lemma 2C, p. 11).  This condition is clearly satisfied for $Y^d - X(X-1)$ for any $d$.

Then we easily have
\begin{equation*}
\sum_b P(E;\phi) \ll \sum_{p \leq P} |S(p)| \frac{\log{p}}{p\log{X}} \ll P^{1 + \varepsilon}
\end{equation*}
which is required to be $\ll X^{1/2}(\log{X})^{-1}$, i.e. supp $\widehat{\phi} \subset (- \onehalf, \onehalf)$.
\end{proof}

\section{Curves with Five-Torsion} 
In this section we consider the family of curves $E=E(b)$ given by
\begin{equation*}
E: y^2 + xy - bxy - by = x^3 - bx^2.
\end{equation*}
This curve has discriminant
\begin{equation*}
\Delta = b^5(b^2 - 11b - 1).
\end{equation*}
The parameter $c_4$ is
\begin{equation*}
c_4 = b^4 - 12 b^3 + 14b^2 + 12b + 1.
\end{equation*}
The torsion group has order 5 and is generated by the point $(0, 0)$.  Now we wish to write $\lambda(p)$ as a character sum so we complete the square:  
\begin{equation*}
4\left(y + \frac{(1-b)x-b}{2}\right)^2 = 4x^3 - 4bx^2 + ((1-b)x-b)^2 
\end{equation*}
so for $p > 2$
\begin{equation*}
\lambda(p) = - \sum_{x \shortmod{p}} \left( \frac{ 4x^3 - 4bx^2 + ((1-b)x-b)^2  }{p} \right).
\end{equation*}

We have the following
\begin{mytheo}
\label{thm:z5torsion}
Let $\mathcal{F}$ be the family of elliptic curves given by the Weierstrass equations $E_b:y^2 + xy - bxy - by = x^3 - bx^2$ with $b$ a positive integer.  Let $w$ be a smooth compactly supported function on $\mr^{+} \times \mr^+$ and set $B = X^{1/3}$.  Define $w_X (E_{a, b}) = w\left(\frac{a}{A}, \frac{b}{B}\right)$.  Then we have
\begin{equation*}  
\mathcal{D}(\mathcal{F}; \phi, w_X) \sim [\widehat{\phi}(0) + \thalf \phi(0)]W_X(\mathcal{F}) \; \; \text{as } X \rightarrow \infty,
\end{equation*}
for $\phi$ with $\text{supp } \widehat{\phi} \subset (- \frac{2}{9}, \frac{2}{9})$.
\end{mytheo}

For the conductor condition we have
\begin{mylemma} Let $\mathcal{F}$ be the family given in Theorem \ref{thm:z5torsion}.  Then
\begin{equation*}
\sum_b \frac{\log{N}}{\log{X}} w \left( \frac{b}{B} \right) = W_X(\mathcal{F}) + O\left(\frac{B}{\log{X}}\right).
\end{equation*}
\end{mylemma}
\begin{proof}  As usual, we will apply Lemma \ref{lem:lemma3a}.  We take $R(d) = b(b^2 - 11b - 1)$.  Then it's clear that $R(d) \asymp X$ and that the irreducible factors of $R(d)$ all divide $\Delta(d)$.
We first consider
\begin{equation*}
\sum_b \mathop{\sum_{p || b(b^2 - 11b - 1)}}_{p^2 | N, \; p > 3}  w \left( \frac{b}{B} \right) \log{p}.
\end{equation*}
An easy calculation shows the only possible prime divisor of $(c_4, \Delta) $ is $5$.  Therefore the contribution from these terms is $\ll B$. We split the sum
\begin{equation*}
\sum_b \mathop{\sum_{p^\alpha || b(b^2 - 11b - 1)}}_{\alpha > 0}  w \left( \frac{b}{B} \right) \log{p^{\alpha - 1}}
\end{equation*}
into two parts depending on whether $p^\alpha || b$ or $p^\alpha || b^2 - 11b -1$.  In the second case we have the bound
\begin{equation*}
\ll \sum_{p^\alpha \ll B^2}  \left(1 + \frac{B}{p^\alpha} \right) \log{p^{\alpha - 1}} \ll B.
\end{equation*}
The first case is even simpler, so the proof is complete. \end{proof}

\begin{mylemma} Let $\mathcal{F}$ be the family given in Theorem \ref{thm:z5torsion}.  Then
\begin{equation*}
\sum_b P(E; \phi) w \left( \frac{b}{B} \right) \ll \frac{X^{1/3}}{\log{X}}
\end{equation*}
for supp $\widehat{\phi} \subset (- \frac{2}{9}, \frac{2}{9})$.
\end{mylemma}
\begin{proof}  For $p > 2$ set
\begin{equation*}
S(p) = \sum_b  \lambda(p) w \left( \frac{b}{B} \right).
\end{equation*}
Then
\begin{eqnarray*}
S(p) & = & - \sum_{x \shortmod{p}} \sum_b \left( \frac{ 4x^3 - 4bx^2 + ((1-b)x-b)^2 }{p} \right) w \left( \frac{b}{B} \right) \\
& = & - \frac{B}{p} \sum_k \mathop{ \mathop{ \sum \sum}_{x \shortmod{p}} }_{\beta \shortmod{p}} \left( \frac{ 4x^3 - 4 \beta x^2 + ((1- \beta)x - \beta)^2 }{p} \right) e \left( \frac{\beta k}{p} \right) \widehat{w} \left( \frac{kB}{p} \right), 
\end{eqnarray*}
by Poisson summation.  Expanding the polynomial, extracting the term $k = 0$, and estimating the terms with $k \neq 0$ with the Weil's bound (RH for curves) gives
\begin{eqnarray*}
S(p) = - \frac{B}{p} \widehat{w} (0)  \mathop{ \mathop{ \sum \sum}_{x \shortmod{p}} }_{\beta \shortmod{p}} \left( \frac{ \beta^2(1+x)^2 + \beta(-2x - 6x^2) + 4x^3 + x^2 }{p} \right) + O(p^{3/2}).
\end{eqnarray*}
Separating the terms with $x = 0$ and executing the change of variables $\beta \rightarrow x \beta$ for $x \not \equiv 0$ gives
\begin{eqnarray*}
S(p) & = &   - \frac{B}{p} \widehat{w} (0)  \sum_{x \not = 0 } \sum_{\beta \shortmod{p}} \left( \frac{ \beta^2(1+x)^2 + \beta (-2 - 6x) + (4x + 1) }{p} \right)  \\
& & - \frac{B}{p} \widehat{w} (0) (p-1) + O(p^{3/2}).
\end{eqnarray*}
After evaluation of the sum over $\beta$ we get
\begin{eqnarray*}
S(p) & =  & - \frac{B}{p} \widehat{w} (0)  \sum_{x \not = 0 } p \delta((-2 - 6x)^2 - 4(1 + x)^2(4x + 1))  + \frac{B}{p} \widehat{w} (0) + O(p^{3/2}) \\
& = & - \frac{B}{p} \widehat{w} (0)  \sum_{x \not = 0} p \delta(x^3)  + \frac{B}{p} \widehat{w} (0) + O(p^{3/2}) \\
& \ll & p^{3/2} + \frac{B}{p}.
\end{eqnarray*}

Then
\begin{equation*}
\sum_b P(E;\phi) \ll \sum_{p \leq P} |S(p)| \frac{\log{p}}{p\log{X}} \ll \frac{P^{3/2}}{\log{X}} + \frac{B}{\log{X}}
\end{equation*}
which is required to be $\ll X^{1/3}(\log{X})^{-1}$, i.e. we require $\text{supp } \widehat{\phi} \subset (-\frac{2}{9}, \frac{2}{9})$.
\end{proof}

\section{Some Families with Complex Multiplication}
\label{section:CM}
\subsection{The Family $y^2 = x^3 + b$}
\begin{mytheo}
\label{thm:CM1} 
Let $\mathcal{F}$ be the family of elliptic curves given by the Weierstrass equations $E_b:y^2 = x^3 + b$ with $b$ a positive integer.  Let $w$ be a smooth compactly supported function on $\mr^+ \times \mr^+$ and set $B = X^{1/2}$.  Define $w_X(E_{b}) = w\left(\frac{b}{B} \right) $.  Then we have
\begin{equation*}  
\mathcal{D}(\mathcal{F}; \phi, w_X) \sim [\widehat{\phi}(0) + \thalf \phi(0)]W_X(\mathcal{F}) \; \; \text{as } X \rightarrow \infty,
\end{equation*}
for $\phi$ with $\text{supp } \widehat{\phi} \subset (- \frac{1}{2}, \frac{1}{2})$.
\end{mytheo}
For the conductor condition we have
\begin{mylemma} Let $N_b$ be the conductor of the curve $y^2 = x^3 + b$.  Let $B \geq 1$ and set
\begin{equation*}
E(B) = \sum_b \frac{\log{(X/N_b)}}{\log{X}} w \left( \frac{b}{B} \right).
\end{equation*}
Then $E(B) \ll B(\log{X})^{-1}$.
\end{mylemma}
\begin{proof}[Proof of Lemma]
\begin{equation*}
E(B)  =   \sum_b \frac{\log{(X/b^2)}}{\log{X}} w \left( \frac{b}{B} \right) + \sum_b \frac{\log{(b^2/N_b)}}{\log{X}} w \left( \frac{b}{B} \right).
\end{equation*}
The first term is trivially $\ll B (\log{X})^{-1}$.  The second term is treated as usual.  The crucial fact is that $p | N \Rightarrow p^2 | N$. 
\end{proof}

\begin{proof}[Proof of Theorem]
We now need to evaluate $\mathcal{P}(\mathcal{F}; \phi, w_X)$.  As usual, we define
\begin{equation*}
S(p) =  \sum_{b} \lambda(p) w \left( \frac{b}{B} \right).
\end{equation*}
Then
\begin{eqnarray*}
S(p) & = &  - \sum_{x \shortmod{p}} \sum_{b} \left( \frac{x^3 + b}{p} \right) w \left( \frac{b}{B} \right) \\
 & = & - \frac{B}{p} \sum_k \sum_{x \shortmod{p}} \sum_{\beta \shortmod{p}} \left( \frac{x^3 + \beta}{p} \right) e \left( \frac{\beta k}{p} \right) \widehat{w} \left( \frac{kB}{p} \right) \\
& = & - \frac{B}{p} \sum_k \sum_{x \shortmod{p}} \sum_{\beta \shortmod{p}} \left( \frac{\beta}{p} \right) e \left( \frac{\beta k - x^3 k}{p} \right) \widehat{w} \left( \frac{kB}{p} \right) \\
& = & - \frac{B}{p^{1/2}} \varepsilon_p \sum_k \sum_{x \shortmod{p}} \left( \frac{k}{p} \right) e \left( \frac{ x^3 k}{p} \right) \widehat{w} \left( \frac{kB}{p} \right). 
\end{eqnarray*}
The next step is to recognize that the summation over $x$ can be written in terms of (twisted) cubic Gauss sums.  In case $p \equiv 2 \pmod{3}$ it's easily seen the sum is zero (since then every element of $\mf_p$ is a cube).  If $p \equiv 1 \pmod{3}$ then $p = \pi \bar{\pi}$ splits in $\mz[e^{2 \pi i/3}]$.  Letting $\chi_{\pi}$ denote the cubic residue character $\negthickspace \pmod{\pi}$ and letting $g(\chi_{\pi})$ denote the corresponding Gauss sum, we then have
\begin{equation*}
\sum_{x \shortmod{p}} e \left( \frac{ x^3 k}{p} \right) = \overline{\chi_{\pi}}(k) g(\chi_{\pi}) + \chi_{\pi}(k) g(\overline{\chi_{\pi}}) + p\delta(k).
\end{equation*}

Using only the bound $|g(\chi_{\pi})| \leq p^{1/2}$ we obtain $S(p) \ll p$ and therefore
\begin{equation*}
\sum_b P(E; \phi) w \left( \frac{b}{B} \right) \ll \frac{P}{\log{X}}.
\end{equation*}
Requiring this $\ll B$ gives support up to $1/2$.  We mention that the presence of cubic Gauss sums should not be completely surprising.  Each curve of the form $y^2 = x^3 + b$ has complex multiplication by $\mq(e^{2 \pi i/3})$.  The L-function can be written in terms of the Gr\"{o}ssencharacter $(4b/\pi)_6$.   Here, then, is the source of the cubic and quadratic characters in $S(p)$.  We have used the more elementary approach for simplicity.
\end{proof}

There is an intriguing possibility of obtaining larger support by virtue of oscillation in the Gauss sums as $p$ varies.  Heath-Brown and Patterson have proved that (for $l \in \mz$) 
\begin{equation*}
\sum_{N(\pi) \leq P} \frac{g(\chi_{\pi})}{N(\pi)^{1/2}} \left( \frac{\bar{\pi}}{|\pi|} \right)^l \ll P^{\varepsilon}(P^{30/31} + |l| P^{29/31}).
\end{equation*}
They conjecture that the sum above is $cP^{5/6} + O(P^{1/2 + \varepsilon})$ ($c \neq 0$) when $l = 0$ and is $O(P^{1/2 + \varepsilon})$ when $l \neq 0$
\cite{Heath-BrownPatterson}.  We therefore feel quite safe to suppose
\begin{equation}
\label{eq:weakcubicgaussbound}
\sum_{P < p \leq 2P}  \sum_{x \shortmod{p}} \left(\frac{k}{p} \right) e \left( \frac{ x^3 k}{p} \right) \ll (kP)^{\varepsilon} P^{4/3} 
\end{equation}
In this case we have
\begin{mytheo}
Assuming (\ref{eq:weakcubicgaussbound}) holds, then Theorem \ref{thm:CM1} holds for $\phi$ with $\text{supp } \widehat{\phi} \subset (- \frac{3}{5}, \frac{3}{5})$.
\end{mytheo}
We mention that if we assume the improved bound of $(kP)^{\varepsilon} P$ for (\ref{eq:weakcubicgaussbound}) (assuming $k$ is not a sixth power) then we would obtain support $(-1, 1)$.

\subsection{The Family $y^2 = x^3 + ax$}
\begin{mytheo} 
Let $\mathcal{F}$ be the family of elliptic curves given by the Weierstrass equations $E_a:y^2 = x^3 + ax$ with $a$ a positive integer.  Let $w$ be a smooth compactly supported function on $\mr^+ \times \mr^+$ and set $A = X^{1/2}$.  Define $w_X(E_a) = w\left(\frac{a}{A} \right)$.  Then we have
\begin{equation*}  
\mathcal{D}(\mathcal{F}; \phi, w_X) \sim [\widehat{\phi}(0) + \thalf \phi(0)]W_X(\mathcal{F}) \; \; \text{as } X \rightarrow \infty,
\end{equation*}
for $\phi$ with $\text{supp } \widehat{\phi} \subset (- \frac{1}{2}, \frac{1}{2})$.
\end{mytheo}
This family has discriminant $-64a^3$ and conductor usually of size about $a^2$.  Precisely we have
\begin{mylemma} Let $N_a$ be the conductor of the curve $y^2 = x^3 + ax$.  Set
\begin{equation*}
B(A) = \sum_a \frac{\log{(X/N_a)}}{\log{X}} w \left( \frac{a}{A} \right).
\end{equation*}
Then $B(A) \ll A (\log{X})^{-1}$.
\end{mylemma}
\begin{proof}
\begin{equation*}
B(A)  =   \sum_a \frac{\log{(X/a^2)}}{\log{X}} w \left( \frac{a}{A} \right) + \sum_a \frac{\log{(a^2/N)}}{\log{X}} w \left( \frac{a}{A} \right).
\end{equation*}
The first term is trivially $\ll A(\log{X})^{-1}$.  The second term is treated as usual. 
\end{proof}

To prove the Theorem we need to compute $\mathcal{P}(\mathcal{F}; \phi, w_X)$.  Since this family has complex multiplication over the Gaussian integers we expect (in analogy to the family $y^2 = x^3 + b$) an appearance of a quartic Gauss sum.  As in the previous section, we proceed elementarily.  Set
\begin{equation*}
S(p) =  \sum_{a} \lambda(p) w \left( \frac{a}{A} \right).
\end{equation*}
Then
\begin{eqnarray*}
S(p) & = &  - \sum_{x \shortmod{p}} \sum_{a} \left( \frac{x^3 + ax}{p} \right) w \left( \frac{a}{A} \right) \\
 & = & - \frac{A}{p} \sum_h \sum_{x \shortmod{p}} \sum_{\alpha \shortmod{p}} \left( \frac{x^3 + \alpha x}{p} \right) e \left( \frac{\alpha h}{p} \right) \widehat{w} \left( \frac{hA}{p} \right) \\
& = & - \frac{A}{p} \sum_h \sum_{x \shortmod{p}} \sum_{\alpha \shortmod{p}} \left( \frac{\alpha x}{p} \right) e \left( \frac{\alpha h - x^2 h}{p} \right) \widehat{w} \left( \frac{hA}{p} \right) \\
& = & - \frac{A}{p^{1/2}} \varepsilon_p \sum_h \left( \frac{h}{p} \right) \sum_{x \shortmod{p}} \left( \frac{x}{p} \right) e \left( \frac{- x^2 h}{p} \right) \widehat{w} \left( \frac{hA}{p} \right). 
\end{eqnarray*}
When $p \equiv 3 \pmod{4}$ it's easy to see that the summation over $x$ is zero (by $x \rightarrow - x$).
For $\mathop{p \equiv 1 \pmod{4}}$ we have
\begin{equation*}
\sum_{x \shortmod{p}} \left( \frac{x}{p} \right) e \left( \frac{- x^2 h}{p} \right) =  - \varepsilon_p \sqrt{p} \left(\frac{-h}{p} \right) + \sum_{x \shortmod{p}} e \left( \frac{- x^4 h}{p} \right)
\end{equation*}
\begin{equation*}
= 
\overline{\chi_{\pi}}(-h) g(\chi_{\pi}) + \chi_{\pi}(-h) g(\overline{\chi_{\pi}}),  
\end{equation*}
where $p = \pi \bar{\pi}$ is the prime decomposition of $p$ in $\mz[i]$ when $p \equiv 1 \pmod{4}$, $\chi_{\pi}$ is the quartic residue character $\mymod{\pi}$, and $g(\chi_{\pi})$ is the Gauss sum of $\chi_{\pi}$.
Now we easily see $S(p) \ll p$ and 
therefore
\begin{equation*}
\sum_b P(E; \phi) w \left( \frac{a}{A} \right) \ll P.
\end{equation*}
Requiring this $\ll A = X^{1/2}$ gives support up to $1/2$.  At this point we could obtain larger support if we assume better bounds on the sum
\begin{equation*}
\sum_{N(\pi) \leq P} \overline{\chi_{\pi}}(h) g(\chi_{\pi}).
\end{equation*}
In principle this should be no different than the cubic case.

\section{A Positive Rank Family of Quadratic Twists}
\label{section:rankonequadratictwists}
Many people have studied families of quadratic twists, especially the aspect of rank frequencies amongst these families.  For instance, Stewart and Top \cite{ST} and Rubin and Silverberg \cite{RS} have produced examples of infinite families of twists with ranks as large as 3.  
In this section we investigate families composed of quadratic twists of a given fixed curve.  

We consider a fixed elliptic curve
\begin{equation*}
E: y^2 = x^3 + ax + b
\end{equation*}
and twist it by an integer $d$, giving the curve
\begin{equation*}
E_d: dy^2 = x^3 + ax + b.
\end{equation*}
Let $N$ be the conductor of $E$ and $N_d$ be the conductor of $E_d$.  If $d$ is squarefree and $(d, N) = 1$ the twisted conductor is given by
\begin{equation*}
N_d = d^2 N.
\end{equation*}
Provided $d \equiv 1 \pmod{4}$ and $(d, N) = 1$ the root number $w_d$ of $E_d$ satisfies
\begin{equation*}
w_d = \left( \frac{d}{-N} \right) w,
\end{equation*}
where $w$ is the root number of $E$.

Now let $\mathcal{A}$ be a set of integers.  We set
\begin{equation*}
\mathcal{D}(\mathcal{A}; \phi, w_X) = \sum_{d \in \mathcal{A}} D(E_d ; \phi) w_X(d).
\end{equation*}
We measure this sum against
\begin{equation*}
W_X(\mathcal{A}) = \sum_{d \in \mathcal{A}} w_X ( d).
\end{equation*}
Before specializing the family $\mathcal{A}$ we set up our machinery.  As always, we will employ the explicit formula (\ref{eq:explicitformula}).  To obtain a density we need the conductor condition
\begin{equation}
\label{eq:conditiona}
\sum_{d \in \mathcal{A}} \frac{\log{N_d}}{\log{X}} w_X (d) \sim \sum_{d \in \mathcal{A}} w_X(d) \; \; \text{as } X \rightarrow \infty
\end{equation}
and an evaluation of the sum of the Fourier coefficients
\begin{equation}
\sum_{d \in \mathcal{A}} P(E_d ; \phi) w_X (d).
\end{equation}
The sum of the Fourier coefficients may reveal extra zeros, depending on the family $\mathcal{A}$.  It is particularly simple to calculate $\lambda(p)$ for quadratic twists since $\lambda_{E_d}(p) = \left( \frac{d}{p} \right) \lambda_{E}(p)$.  Thus
\begin{equation}
\label{eq:conditionb}
\sum_{d \in \mathcal{A}} P(E_d ; \phi) w_X (d) =  \sum_{p} \lambda_{E}(p) \frac{2\log{p}}{p \log{X}} \widehat{\phi}\left( \frac{\log{p}}{\log{X}} \right) \sum_{d \in \mathcal{A}} \left( \frac{d}{p} \right) w_X (d).
\end{equation}


We consider a lacunary family of quadratic twists.  Instead of twisting by all squarefree integers we twist by integers $d$ of the form $d = u^3 + au + b$.  When we do so the twisted curve $E_d$ always has the point $(u, 1)$ so in general we expect the curve to have positive rank.  For this family we have the density theorem
\begin{mytheo} 
\label{thm:rank1quadratic}
Let $E: y^2 = x^3 + ax + b$ be a fixed elliptic curve.  Let $\mathcal{A}$ be the set of all integers of the form $u^3 + au + b$ for $u$ an integer.  Let $w$ be a smooth compactly supported function on $\mr^{+}$, set $U = X^{1/6}$, and take $w_X(d) = w\left(\frac{u}{U}\right)$.  Then we have
\begin{equation*}  
\mathcal{D}(\mathcal{A}; \phi, w_X) \sim [\widehat{\phi}(0) + {\textstyle \frac{3}{2}} \phi(0)] W_X (\mathcal{A}) \; \; \text{as } X \rightarrow \infty, 
\end{equation*}
for $\phi$ with $\text{supp } \widehat{\phi} \subset (- \frac{1}{6}, \frac{1}{6})$.
\end{mytheo}
Notice that we can sum primes up to the size of the family here.

For notational simplicity we sum over $u$ and take $w\left(\frac{u}{U}\right)$ as our test function.
We first treat the conductor condition (\ref{eq:conditiona}).  The result is
\begin{mylemma} 
\label{lem:twistconductor}
With notation as in Theorem \ref{thm:rank1quadratic}
\begin{equation*}
\sum_u \frac{\log{ N_{u^3 + au + b}}}{\log{X}} w \left( \frac{u}{U} \right) = X^{1/6} \widehat{w}(0) + O\left(\frac{X^{1/6}}{\log^{1/4}{X}}\right).
\end{equation*}
\end{mylemma}

\begin{proof}  It suffices to consider
\begin{equation*}
\frac{1}{\log{U}} \sum_u \log\left({\frac{ N_{u^3 + au + b} }{(u^3 + au + b)^2} }\right) = \frac{1}{\log{U}} \sum_u \mathop{\sum_{p^\beta || N_{u^3 + au + b}}}_{p^{2 \alpha} || (u^3 + au + b)^2} \log{p^{\beta - 2\alpha}} w \left( \frac{u}{U} \right) .
\end{equation*}
We may trivially sum over the terms with $\alpha = 0$ so we may assume $p | u^3 + au + b$ and therefore $\beta = 2$.  Thus we obtain the sum
\begin{equation*}
\frac{1}{\log{U}} \sum_u  \mathop{\sum_{p^\alpha || u^3 + au + b}}_{\alpha > 0} 2 \log{p^{\alpha - 1}} w \left( \frac{u}{U} \right).
\end{equation*}

By splitting the summation over $u$ into progressions $\mymod{ p^{\alpha}}$ we get the bound
\begin{equation*}
\frac{1}{\log{U}} \mathop{\sum_{p^\alpha \ll U^3 }}_{\alpha \geq 3} \log{p^{\alpha - 1}} \left(1 + U/p^{\alpha} \right) \ll \frac{U}{\log{U}},
\end{equation*} 
so we may assume $\alpha = 2$.  Using the same technique we can sum $p \ll U$ when $\alpha = 2$.  In case $u^3 + au +b$ is reducible (over \mz) we immediately get $p \ll U$, so we may assume $u^3 + au + b$ is irreducible.  
To handle large $p$ we use the same method as Lemma \ref{lem:7bigp} to obtain
\begin{equation*}
\sum_{p \geq P} \sum_{u^3 + au + b \equiv 0 \shortmod{p^2}} \frac{\log{p}}{\log{U}} \; w \left( \frac{u}{U} \right) \ll \sum_{d \geq P} \sum_{u^3 + au + b = d^2 l} \left|w \left( \frac{u}{U} \right) \right|.
\end{equation*}
We may bound $|w|$ by a smooth non-negative function with slightly larger support.  Abusing notation we rename $w$ to be this new function.  Setting
\begin{equation*}
S(l) = \sum_{u^3 + au + b = d^2l }  w \left( \frac{u}{U} \right) 
\end{equation*}
we get
\begin{equation*}
S(l) \left( 1 - \frac{\log^2 U}{|\mathcal{Q}|^2 \log^2 Q}\right) 
\end{equation*}
\begin{equation*}
 \leq  \frac{1}{|\mathcal{Q}|^2} \sum_{q_1} \sum_{q_2} \sum_{u^3 + au + b \equiv 0 \shortmod{l}} \left( \frac{l}{q_1q_2} \right) \left( \frac{u^3 + au + b}{q_1q_2} \right)w \left( \frac{u}{U} \right), 
\end{equation*}
where $\mathcal{Q}$ and $Q$ are defined as in Lemma \ref{lem:7bigp}.  We have compressed a number of steps because the calculations are identical to those used in the proof of Lemma \ref{lem:7bigp}.  Let $S_1$ be the inner sum over $u$, set $r = q_1 q_2$, and apply Poisson summation in $u \pmod{lr}$ to obtain 
\begin{equation*}
S_1 = \frac{U}{lr} \sum_h \mathop{\sum_{x \shortmod{ lr}}}_{x^3 + ax + b \equiv 0 \shortmod{l}} \left( \frac{l}{r} \right) \left( \frac{x^3 + ax + b}{r} \right) e \left( \frac{xh}{lr} \right) \widehat{w} \left( \frac{hU}{lr} \right).
\end{equation*}
Write $x = x_1 r + x_2 l$ where $x_1$ is given $\mymod{l}$ and $x_2$ is given $\mymod{r}$ and get
\begin{equation*}
S_1 = \frac{U}{lr} \left( \frac{l}{r} \right)  \sum_h \widehat{w} \left( \frac{hU}{lr} \right) Y(h, l, r),
\end{equation*} 
where
\begin{equation*}
Y(h, l, r) = \sum_{x_2 \shortmod{r}} \left( \frac{x_2^3 + ax_2  + b}{r} \right) 
e \left( \frac{x_2h \bar{l}}{r} \right) 
 \mathop{\sum_{x_1 \shortmod{ l}}}_{x_1^3 r^3 + ax_1 r + b \equiv 0 \shortmod{l}} e \left( \frac{x_1h}{l} \right).
\end{equation*}
If $q_1 \neq q_2$ the summation over $x_2$ is $O(r^{1/2})$ by the Riemann Hypothesis for curves (Weil's bound).  If $q_1 = q_2$ then the summation over $x_2$ is $r \delta_r(h) + O(1)$ (here $\delta_r(y)$ is the characteristic function of the arithmetic progression $y \equiv 0 \pmod{r}$).  The summation over $x_1$ is bounded by $\rho(l)$, where $\rho(l)$ is the number of solutions to $x^3 + ax + b \equiv 0 \pmod{l}$ (since we may assume $(l, r) = 1$).  Therefore we have the bound
\begin{equation*}
S_1 \ll \rho(l) \left( r^{1/2} +  \frac{U}{l r^{1/2}}+  \frac{U}{l} \delta_{q_1 = q_2} \right)
\end{equation*}
and hence
\begin{equation*}
S(l) \left( 1 - \frac{\log^2 U}{|\mathcal{Q}|^2 \log^2 Q}\right) \ll \frac{U}{|\mathcal{Q}|} \frac{\rho(l)}{l} + Q \rho(l).
\end{equation*}
Using the Prime Ideal Theorem we get the bounds $\sum_{l \leq L} \rho(l) \ll L$ and $\sum_{l \leq L} \rho(l) l^{-1} \ll 1$ (since we may assume $u^3 + au + b$ is irreducible), so 
\begin{equation*}
\sum_{l \leq L} S(l) \left( 1 - \frac{\log^2 U}{|\mathcal{Q}|^2 \log^2 Q}\right) \ll \frac{U}{|\mathcal{Q}|} + L Q.
\end{equation*}
Now we take $Q \asymp \log{U} \log \log{U}$, $|\mathcal{Q}| \gg \log U$, and $L = U/\log^{3/2}{U}$ to obtain
\begin{equation*}
\sum_{l \leq L} S(l) \ll \frac{U}{\sqrt{\log{U}}} 
\end{equation*} 
Converting $L$ to $P$ via $LP^2 \asymp U^3$ means we can sum over $p \geq U \log^{3/4}{U}$.  To close the gap we exploit the rarity of primes; for any $U_1 \geq U$ we have
\begin{equation*}
\sum_{U \leq p \leq U_1} \sum_{u^3 + au + b \equiv 0 \shortmod{p^2}} w \left( \frac{u}{U} \right) \ll \sum_{U \leq p \leq U_1} \left( 1 + \frac{U}{p} \right)
\end{equation*}
\begin{equation*}
\ll \frac{U_1}{\log{U_1}} \ll \frac{U}{\log^{1/4}{U}},
\end{equation*}
for the choice $U_1 = U \log^{3/4}{U}$.  The proof is now complete. \end{proof}

As for condition (\ref{eq:conditionb}), we calculate
\begin{equation*}
- \sum_u \left( \frac{u^3 + au + b}{p} \right) w \left(\frac{u}{U} \right) = \frac{U}{p} \lambda_E (p) \widehat{w}(0) + O(p^{1/2}).
\end{equation*}
Therefore we obtain
\begin{eqnarray*}
& & \sum_{p \leq P} \lambda_{E}(p) \frac{2\log{p}}{p \log{X}} \widehat{\phi}\left( \frac{\log{p}}{\log{X}} \right) \left( \frac{U}{p} \lambda_E (p) \widehat{w}(0) + O(p^{1/2}) \right) \\
& = & U \widehat{w}(0) \sum_{p \leq P} \lambda_{E}^2(p) \frac{2\log{p}}{p^2 \log{X}} \widehat{\phi}\left( \frac{\log{p}}{\log{X}} \right) + O\left(\frac{P}{\log{X}}\right) \\
& = & U \widehat{w}(0) \sum_{p \leq P} \chi_0 (p) \frac{2 \log{p}}{p \log{X}}\widehat{\phi}\left( \frac{\log{p}}{\log{X}} \right)  \\
& + & U \widehat{w}(0) \sum_{p \leq P} \left[ \frac{\lambda_{f}^2(p)}{p} - \chi_0(p) \right] \frac{2\log{p}}{p \log{X}} \widehat{\phi}\left( \frac{\log{p}}{\log{X}} \right) + O\left(\frac{P}{\log{X}}\right) \\
& = & U \widehat{w}(0) \phi(0) +  O\left(\frac{U \log{\log{X}}}{\log{X}}\right) + O\left(\frac{P}{\log{X}}\right),
\end{eqnarray*}
by the Prime Number Theorem and the Riemann Hypothesis for the symmetric-square L-function associated to $E$.  This bound is sufficient provided $P \ll U$, i.e., $\text{supp } \widehat{\phi} \subset \left(- \frac{1}{6}, \frac{1}{6} \right)$. 

\appendix
\section{A Technical Character Sum}
\label{section:appendixa}
Recall our goal is to obtain the bound
\begin{equation*}
H(s_1, s_2) \ll \frac{U^{\sigma_1}V^{\sigma_2}}{|(s_1/c + s_2)(s_1/c + s_2 + 1) ||s_2|^{1 + \varepsilon}} \left(1 + \frac{U^c}{V|q|} \right)^{1/2}, 
\end{equation*}
where
\begin{equation*}
H(s_1, s_2) = \int_0^\infty \int_0^\infty e\left( \frac{u^c}{vq} \right) F(u, v) u^{s_1} v^{s_2} \frac{du dv}{uv}.
\end{equation*}

We note first that the trivial bound on $H(s_1, s_2)$ is $U^{\sigma_1}V^{\sigma_2}$.  We begin with the change of variables $t = u^c/{vq}$, obtaining
\begin{equation*}
H(s_1, s_2) = \int \int e(t) F((tvq)^{1/c}, v) (tvq)^{s_1/c} v^{s_2} \frac{dt dv}{c t v}  \end{equation*}
\begin{equation*}
= \int e(t) (tq)^{s_1/c} \left(\int F((tvq)^{1/c}, v) v^{s_1/c + s_2 -1} dv \right) \frac{dt}{ct}.
\end{equation*}

In the inner integral we apply integration by parts twice, obtaining
\begin{eqnarray*}
& & \frac{1}{(s_1/c + s_2)(s_1/c + s_2 + 1)} \int \frac{\partial^2 F((tvq)^{1/c}, v)  }{\partial v^2} v^{s_1/c + s_2 +1} dv \\
& = & \frac{1}{(s_1/c + s_2)(s_1/c + s_2 + 1)} \int G((tvq)^{1/c}, v) v^{s_1/c + s_2} \frac{dv}{v}
\end{eqnarray*}
where
\begin{equation*}
G(u, v)  = c_1 F^{(2, 0)}(u, v) u^{2} + c_2 F^{(1, 1)}(u, v) u v + c_3 F^{(1, 0)}(u, v) u + c_4 F^{(0, 2)}(u, v) v^2
\end{equation*}
and $c_1$ through $c_4$ are absolute constants.  Notice that $G$ satisfies (\ref{eq:smoothnesscondition}).  Therefore, at the expense of differentiating twice we have obtained the convergence factor $((s_1/c + s_2)(s_1/c + s_2 + 1))^{-1}$, i.e.,
\begin{equation*}
H(s_1, s_2) \ll  |(s_1/c + s_2)(s_1/c + s_2 + 1)|^{-1} \left| \int \int e(t) G((tvq)^{1/c}, v) (tvq)^{s_1/c} v^{s_2} \frac{dt dv}{tv} \right|.
\end{equation*}
Now we concentrate on the $t$-variable aspect, specifically the integral
\begin{equation*}
\int_0^\infty e(t) G((tvq)^{1/c}, v) t^{s_1/c - 1} dt
\end{equation*}
\begin{equation*}
=\int_0^\infty e^{i L(t)} G((tvq)^{1/c}, v) t^{-1 + \sigma_1/c} dt
\end{equation*}
where $s_1 = \sigma_1 + i r_1$ and $L(t) = 2\pi t + (r_1 \log{t})/c$.  We expect that there should be a lot of cancellation in this integral as long as $L(t)$ has some variation, i.e. away from points where $L'(t) = 0$.  $L'(t)$ has its only zero at $t_0 = -r_1/2c\pi$.  Suppose $t_0 > 0$.  Let $f_1 + f_2 + f_3$ be a partition of unity such that 
\begin{equation*}
f_1(t) = \begin{cases} 1 & \text{for } t < 1/2 \\   0 & \text{for } t > 3/4 \end{cases}, \; \; 
f_2(t) = \begin{cases} 0 & \text{for } |t - 1| > 1/2 \\   1 & \text{for } |t - 1| < 1/4 \end{cases}, \; \; 
f_3(t) = \begin{cases} 0 & \text{for } t < 5/4 \\   1 & \text{for } t > 3/2 \end{cases}.
\end{equation*}
Using the partition (with arguments scaled by $t_0$) break up the integral into three pieces
\begin{equation*}
\int_{0}^{3t_0/4} + \int_{t_0/2}^{3t_0/2} + \int_{5t_0/4}^{\infty}.
\end{equation*}

For the first integral, we compute
\begin{equation*}
\int_0^{3t_0/4} e^{i L(t)} G((tvq)^{1/c}, v) t^{-1 + \sigma_1/c} f_1(t/t_0) dt \end{equation*}
\begin{equation*}
= \int_0^{3t_0/4} iL'(t)e^{i L(t)} \frac{G((tvq)^{1/c}, v) t^{-1 + \sigma_1/c}f_1(t/t_0)}{iL'(t)} dt
\end{equation*}
\begin{equation*}
=  - \int_0^{3t_0/4} e^{i L(t)} \left[\frac{G((tvq)^{1/c}, v) t^{-1 + \sigma_1/c}f_1(t/t_0)}{iL'(t)}\right]' dt
\end{equation*}
\begin{equation*}
= -i\int_0^{3t_0/4} e^{i L(t)} \frac{G_1((tvq)^{1/c}, v) t^{-2 + \sigma_1/c} L'(t) +G_2((tvq)^{1/c}, v) t^{-1 + \sigma_1/c}  L''(t)}{L'(t)^2} dt 
\end{equation*}
where $G_1$ and $G_2$ satisfy (\ref{eq:smoothnesscondition}).  The boundary term for $t = 0$ is nonexistent because $L'(t) \asymp t^{-1}$ as $t \rightarrow 0$.  
The integral simplifies to
\begin{equation*}
\int_0^{3t_0/4} e^{i L(t)} G_3((tvq)^{1/c}, v) t^{-1 + \sigma_1/c} R(t) dt,
\end{equation*}
where $G_3$ satisfies (\ref{eq:smoothnesscondition}) and $R(t) = \frac{d_0t + d_1 r_1}{( t -t_0)^2}$ for $d_0$ and $d_1$ absolute constants.  
Apply integration by parts one more time in this integral and obtain
\begin{equation*}
- \int_0^{3t_0/4} e^{i L(t)} \left[\frac{G_3((tvq)^{1/c}, v) t^{-1 + \sigma_1/c} R(t)}{iL'(t)}\right]' dt.
\end{equation*}
The boundary term at $3t_0/4$ is zero because $f_1 \equiv 0$ for $t \geq3/4$.  The boundary term at $0$ is zero because $R(t) \asymp 1$ as $t \rightarrow 0$. 
Now
\begin{equation*}
\left[\frac{G_3((tvq)^{1/c}, v) t^{-1 + \sigma_1/c}R(t)}{iL'(t)}\right]' \end{equation*}
\begin{equation*}
= \left[\frac{G_3((tvq)^{1/c}, v) t^{-1 + \sigma_1/c}}{iL'(t)}\right]' R(t) + \left[\frac{G_3((tvq)^{1/c}, v) t^{-1 + \sigma_1/c}}{iL'(t)}\right] R'(t)
\end{equation*}
\begin{equation*}
= G_4 ((tvq)^{1/c}, v) t^{-1 + \sigma_1/c} \tilde{R}(t)
\end{equation*}
where $G_4(u, v) = c_1 G_3^{(1, 0)} u + c_2 G_3(u, v)$ for some pair of constants $c_1$ and $c_2$ and 
\begin{equation*}
\tilde{R}(t) = \frac{d_2 t^2 + d_3 t_0 t + d_4 {t_0}^2}{(t - t_0)^4}
\end{equation*}
for absolute constants $d_2, d_3$, and $d_4$.  Again, $G_4$ satisfies (\ref{eq:smoothnesscondition}).  The conclusion is that at the expense of differentiating twice more we have gained the convergence factor $\tilde{R}(t)$.  

Estimating the integral trivially gives the bound
\begin{equation*}
\frac{|q|^{\sigma_1/c}}{t_0^2} \int_0^\infty  \int_0^{3t_0/4} \left(1 + \frac{(tv|q|)^{1/c}}{U} \right)^{-2} \left(1 + \frac{v}{V} \right)^{-2}  t^{-1 + \sigma_1/c} v^{\sigma_1/c + \sigma_2 - 1} dt dv,
\end{equation*}
using the bound $|\tilde{R}(t)| \ll t_0^{-2}$.
Simplifying, it's
\begin{equation*}
\ll \frac{|q|^{\sigma_1/c}}{t_0^2} \int_0^\infty  \left(\text{min}\left\{t_0, \frac{U^c}{v |q|} \right\} \right)^{\sigma_1/c} \left(1 + \frac{v}{V} \right)^{-2}  v^{\sigma_1/c + \sigma_2 - 1} dv
\end{equation*}
\begin{equation*}
\ll t_0^{-2} U^{\sigma_1} V^{\sigma_2}.
\end{equation*}


Therefore we get the bound
\begin{equation*}
\int_0^{\infty} \int_0^{3t_0/4} \ll \frac{U^{\sigma_1} V^{\sigma_2}}{|s_1/c +s_2||s_1/c + s_2 + 1| |s_1|^{2}}. 
\end{equation*}
We obtain the same bound for the integration over $t \geq 5t_0/4$.  The reason is that the integration by parts gives the same convergence factor $\tilde{R}(t)$.  
The only change is that we apply the bound $\tilde{R}(t) \ll t^{-2}$ to get convergence at infinity.  The bound is the same.


In the case that $t_0 < 0$ there is no need to break up the integral at all and we gain the convergence factor $\tilde{R}(t)$, which immediately saves us $|t_0|^{-2}$.

It remains to bound the integral
\begin{equation*}
 \int_{t_0/2}^{3t_0/2} e^{i L(t)} G((tvq)^{1/c}, v) t^{-1 + \sigma_1/c} f_2(t/t_0) dt
\end{equation*}
\begin{equation*}
=  t_0^{\sigma_1/c}e^{i L(t_0)} \int_{-1/2}^{1/2} e^{i [L(t_0(1 + t)) - L(t_0)]} G((t_0(1 + t)vq)^{1/c}, v) (1 + t)^{-1 + \sigma_1/c} f_2(1 + t)dt.
\end{equation*}
Set $t_0 \Phi(t) = L(t_0(1 + t)) - L(t_0) = 2\pi t_0 (t - \log(1 + t))$ and define
\begin{equation*}
a(t, t_0) = G((t_0(1 + t)vq)^{1/c}, v) \left(1 + \frac{v}{V}\right)^2 \left( 1 + \frac{(t_0vq)^{1/c}}{U} \right)^2 (1 + t)^{-1 + \sigma_1/c} f_2(1 + t).
\end{equation*}
With these changes in notation the integral is
\begin{equation*}
  t_0^{\sigma_1/c}e^{i L(t_0)} \left(1 + \frac{v}{V}\right)^{-2} \left( 1 + \frac{(t_0vq)^{1/c}}{U} \right)^{-2} \int_{-\infty}^{\infty} e^{it_0 \Phi(t) } a(t, t_0) dt.
\end{equation*}
Then it's clear from (\ref{eq:smoothnesscondition}) that $a(t, t_0)$ satisfies (for $t_0 > 1$, as we may assume)
\begin{equation*}
\left| \left(\frac{\partial}{\partial t} \right)^\alpha \left(\frac{\partial}{\partial t_0} \right)^\gamma a(t, t_0) \right| \ll (1 + t_0)^{-\gamma} 
\end{equation*}
uniformly in $v$.  Since $\Phi(0) = \Phi'(0) = 0$ and $\Phi'(t) \neq 0$ for $t \neq 0$ we may apply the Van der Corput bound of $O\left(t_0^{-1/2}\right)$ to the integration over $t$ (cf. \cite{Sogge}, Lemma 1.1.2).  Continuing, we integrate over $v$ and get
\begin{equation*}
\ll  t_0^{-1/2 + \sigma_1/c} |q|^{\sigma_1/c} \int_0^{\infty} \left(1 + \frac{v}{V}\right)^{-2} \left( 1 + \frac{(t_0vq)^{1/c}}{U} \right)^{-2} v^{\sigma_1/c +\sigma_2  -1} dv
\end{equation*}
\begin{equation*}
\ll t_0^{-1 - \varepsilon} U^{\sigma_1} V^{\sigma_2} \left( 1 + \left(\frac{U^c}{V|q|}\right)^{1/2 + \varepsilon} \right).
\end{equation*}
Therefore
\begin{equation*}
H(s_1, s_2) \ll \frac{U^{\sigma_1} V^{\sigma_2} }{|(s_1/c + s_2)(s_1/c + s_2 + 1)| |s_1|^{1 + \varepsilon}} \left( 1 + \left(\frac{U^c}{V|q|}\right)^{1/2 + \varepsilon} \right)
\end{equation*}
and the proof of the main statement is complete.

\section{A Technical Exponential Sum}
\label{section:appendixb}
This appendix is dedicated to stating and proving the following general result.
\begin{mylemma}
\label{lem:appendixb}
Let $f$ and $g$ be smooth real-valued functions defined on an open interval containing $[1, 2]$.  Suppose $f(x) \asymp 1 \asymp f'(x)$ and $g^{(k)}(x) \asymp 1$ for $k = 0, 1, 2, 3$.
Let $c_n$, $n = 1, 2, \ldots$ be arbitrary complex numbers satisfying $|c_n| \leq 1$, and let $M \geq 1, N \geq 1$, and $Y$ be real numbers.  Consider the sum
\begin{equation*}
S = \sum_{M \leq m < 2M} \left| \sum_{N \leq n < 2N} c_n e\left( Y f(n/N) g(m/M) \right) \right|.
\end{equation*}
If $M \leq C|Y|$ ($C$ a positive real number) then the bound
\begin{equation*}
S \ll N^{1/2} M + \frac{NM}{1 +  |Y|^{1/2}} \log{N}
\end{equation*}
holds.
In case $M > C|Y| $ the bound
\begin{equation*}
S \ll N^{1/2} M + \frac{NM}{1 + |Y|^{1/4}}
\end{equation*}
holds.  The implied constants in these bounds depend only on $f$, $g$, and $C$.
\end{mylemma}
\begin{proof}  
Let $F_{\varepsilon}(x)$ be a smooth, non-negative function which takes the value $1$ for $1 \leq x \leq 2$ and has support in the interval $(1 - \varepsilon, 2 + \varepsilon)$.  Suppose $\varepsilon > 0$ is small enough so that $f$ and $g$ are defined on the interval $(1-\varepsilon, 2 + \varepsilon)$.
Then
\begin{equation*}
S \leq \sum_{m \in \mz} F_{\varepsilon}\left(\frac{m}{M}\right) \left| \sum_{N \leq n < 2N} e\left( Y f(n/N) g(m/M) \right) \right|.
\end{equation*}
Applying Cauchy's inequality, we obtain
\begin{equation*}
S^2 \ll M \sum_{N \leq n_1 < 2N}  \sum_{N \leq n_2 < 2N} c_{n_1} \overline{c_{n_2}} \sum_{m} F_{\varepsilon}\left(\frac{m}{M}\right) e\left( Y (f(n_1/N) - f(n_2/N)) g(m/M) \right). 
\end{equation*}
The diagonal terms contribute $O(NM^2)$.  

Let $T$ be any real number and consider the sum
\begin{equation}
S_1 (T) = \sum_{m \in \mz} F_{\varepsilon}\left(\frac{m}{M}\right) e\left( T g(m/M) \right).
\end{equation}
Our goal is to obtain the bound
\begin{equation}
\label{eq:S_1bound}
S_1(T) \ll
\begin{cases} M(1 + |T|)^{-1} & \text{when $M \gg |T|$} \\
M(1 + |T|^{1/2})^{-1} & \text{when $M \ll |T|$},
\end{cases}
\end{equation}
with implied constants depending on $f$, $g$, $F_\varepsilon$, and $C$.  Applying this bound to $S$ using $T = Y(f(n_1/N) - f(n_2/N)) \asymp Y (n_1 - n_2) N^{-1}$ will prove Lemma \ref{lem:appendixb} by noting that since $|T| \gg Y$ we can use the better bound for all pairs $n_1, n_2$ if  $M \gg |T|$. 
\end{proof}

\begin{proof}[Proof of (\ref{eq:S_1bound})]  
By Poisson summation,
\begin{equation}
\label{eq:sumtointegral}
\sum_{m \in \mz} F_{\varepsilon} \left(\frac{m}{M}\right) e\left( T g(m/M) \right)
= \sum_{r = - \infty }^{\infty} \int_{- \infty}^{\infty} F_{\varepsilon} \left( \frac{u}{M} \right) e(T g(u/M) - ru) du.
\end{equation}
Apply the change of variables $u \rightarrow Mu$ and obtain
\begin{equation}
\label{eq:appfreq}
\sum_{r = - \infty }^{\infty} M \int_{1 - \varepsilon }^{2 + \varepsilon} F_{\varepsilon}\left(u \right) e(T g(u) - ruM) du.
\end{equation}
Set $G(u) = G_{r, M, T} (u) = 2\pi(T g(u) - ruM)$ and integrate by parts to get
\begin{equation}
\label{eq:firstparts}
M \int_{1 - \varepsilon }^{2 + \varepsilon} F_{\varepsilon}\left(u \right) e(T g(u) - ruM) du = Mi \int e^{i G(u)} \frac{F_{\varepsilon}'(u) G'(u) - F_{\varepsilon}(u) G''(u)}{G'(u)^2} du.
\end{equation}
Integrate by parts again and obtain
\begin{equation*}
- M \int e^{i G(u)} \left(\frac{F_{\varepsilon}''(u) G'(u) - F_{\varepsilon}(u) G'''(u)}{G'(u)^3} + 3 \frac{(F_{\varepsilon}(u) G''(u) - F_{\varepsilon}'(u) G'(u))G''(u)}{G'(u)^4}\right) du.
\end{equation*}
Setting $l = l_{r, T} = \text{min}_u \{|TM^{-1} g'(u) - r| \}$ and estimating this integral trivially gives the bound (provided $l \neq 0$)
\begin{equation}
\label{eq:boundonl}
\ll M^{-1} l^{-2} + |T| M^{-2} l^{-3} + |T|^2 M^{-3} l^{-4},
\end{equation}
the implied constant depending on $f, g$, and $F_{\varepsilon}$ only.  If $|T| \leq \varepsilon M$ for $\varepsilon$ small enough with respect to the implied constant in $g'(x) \asymp 1$ then we can sum over all $r \neq 0$ in (\ref{eq:appfreq}) and get the bound $O(1)$ depending on $f$, $g$, and $F_\varepsilon$ only (but from now on $\varepsilon$ and $F_\varepsilon$ are fixed so in fact all implied constants depend on $f$ and $g$ only).  The case $r = 0$ is easily estimated by (\ref{eq:firstparts}) and contributes at most
\begin{equation*}
\ll \frac{M}{1 + |T|}.
\end{equation*}
In case $|T| \geq \varepsilon M$ we can sum over all $r$ except those of size $|r| \asymp |T| M^{-1}$ (using $|r| \asymp l$ and (\ref{eq:boundonl})) to get the bound $O(1)$.  Now assume we are in the range $|r| \asymp |T|M^{-1}$.  The contribution from those $r$ such that 
\begin{equation}
\label{eq:stationarypoints}
l M |T|^{-1} \geq \delta > 0
\end{equation}
is at most $O_{\delta}(M(1 + |T|^2)^{-1})$.  Thus we may take $\delta$ small enough (with respect to the implied constant in $g'(x) \asymp 1$) so that either (\ref{eq:stationarypoints}) holds or $l = 0$ (since $g'(u) \asymp 1$ and $|r|M|T|^{-1} \asymp 1$).  There are a bounded number (bounded in terms of the various implied constants already mentioned) of possible values of $r$ such that $l = 0$.
For each such $r$ there is one value of $u$, say $u_0$, such that $Tg'(u_0) -rM = 0$.  By taking $\Phi(u) = g(u) - g(u_0) + r MT^{-1} (u - u_0)$ it suffices to bound the integral
\begin{equation*}
M \int_{- \infty}^{\infty} F_{\varepsilon}(u) e\left(T\Phi(u) \right) du,
\end{equation*}
where $\Phi(0) = \Phi'(0) = 0$, $\Phi'(u) \neq 0$ for $u \neq 0$, and $\Phi''(0) \neq 0$ 
Using stationary phase estimates (cf. \cite{Sogge}, Theorem 1.1.1.) we get the bound of 
\begin{equation*}
\ll \frac{M}{1 + |T|^{1/2}}
\end{equation*} 
for this integral, the implied constant depending on $f$ and $g$ only.  Now the proof is complete. 
\end{proof}

\end{document}